\documentclass[fleqn]{cas-sc}

\usepackage[numbers,sort&compress]{natbib}

\def\tsc#1{\csdef{#1}{\textsc{\lowercase{#1}}\xspace}}
\tsc{WGM}
\tsc{QE}
\usepackage{cleveref}
\usepackage{enumerate}
\usepackage{lipsum}
\usepackage{mathtools}
\usepackage{amsmath,amsfonts,amsthm,bm,amsthm,amssymb} 
\usepackage{graphicx}
\usepackage{epstopdf}
\usepackage{algorithmic}
\ifpdf
  \DeclareGraphicsExtensions{.eps,.pdf,.png,.jpg}
\else
  \DeclareGraphicsExtensions{.eps}
\fi

\usepackage{mathtools}
\usepackage{comment}
\usepackage[labelfont=bf]{caption} 
\usepackage{subcaption}
\usepackage{tikz}

\newsavebox{\foobox}

\newcommand{\bunderline}[1]{\underline{#1}}
\usepackage{arydshln}

\usepackage{cleveref}
\usepackage{enumerate}
\usepackage{lipsum}
\usepackage{mathtools}
\usepackage{amsmath,amsfonts,amsthm,bm,amsthm,amssymb} 
\usepackage{graphicx}
\usepackage{epstopdf}
\usepackage{algorithmic}
\ifpdf
  \DeclareGraphicsExtensions{.eps,.pdf,.png,.jpg}
\else
  \DeclareGraphicsExtensions{.eps}
\fi
\usepackage{mathtools}
\usepackage{comment}
\usepackage[labelfont=bf]{caption} 
\usepackage{subcaption}
\usepackage{tikz}

\usepackage{arydshln}

\newtheorem{remark}{Remark}

\DeclareCaptionLabelSeparator{period}{.\ }
\renewcommand{\vec}[1]{{\bunderline{#1}}}

\newcommand{\mat}[1]{{\bunderline{\bunderline{#1}}}}

\newcommand{\reals}{\mathbb R}

\def\maxw{{\mathcal M}}
\def\gauss{{\mathcal G}}

\def\pr{{\mathbb P}}

\usepackage{tcolorbox}

\makeatletter
\renewcommand*\env@matrix[1][\arraystretch]{%
  \edef\arraystretch{#1}%
  \hskip -\arraycolsep
  \let\@ifnextchar\new@ifnextchar
  \array{*\c@MaxMatrixCols c}}
\makeatother

\usepackage{placeins}

\newcommand{\inred}[1]
        {\mbox{\color{red} #1}}
        
\begin{document}
\let\WriteBookmarks\relax
\def\floatpagepagefraction{1}
\def\textpagefraction{.001}

\shorttitle{Micro-macro kinetic flux-vector splitting schemes for the multidimensional Boltzmann-ES-BGK equation}    

\shortauthors{J.A. Rossmanith and P. Sar}  

\title [mode = title]{Micro-macro kinetic flux-vector splitting schemes for the multidimensional Boltzmann-ES-BGK equation}   


%

\author[1]{James A. Rossmanith}[type=editor,orcid=0000-0002-3629-8895]
\ead{rossmani@iastate.edu}
\cormark[1]

\author[2]{Preeti Sar}[orcid=0000-0002-4937-7023]
\ead{sarp@ornl.gov}



\affiliation[1]{organization={Deparment of Mathematics, Iowa State University},
            addressline={411 Morrill Road}, 
            city={Ames},
            postcode={50011}, 
            state={IA},
            country={USA}}

\affiliation[2]{organization={Fusion Energy Division, Oak Ridge National Laboratory},
    city={Oak Ridge},
    postcode={37831}, 
    state={TN},
    country={USA}}

\credit{}

\cortext[1]{Corresponding author}



\begin{abstract}
The kinetic Boltzmann equation models gas dynamics over a wide range of spatial and temporal scales. Simplified versions of the full Boltzmann collision operator, such as the classical Bhatnagar-Gross-Krook (BGK) and the closely related Ellipsoidal-Statistical-BGK (ES-BGK) operators, can dramatically reduce the computational cost of solving kinetic equations numerically. Classical BGK yields incorrect transport coefficients (relative to the full Boltzmann collision operator) at low Knudsen numbers, whereas ES-BGK captures them correctly. In this work, we develop a finite-volume method based on a micro-macro decomposition of the distribution function, which requires a smaller velocity mesh than direct kinetic methods for low and intermediate Knudsen numbers. The macro portion of the model is a fluid model with a moment closure derived from the heat-flux tensor calculated from the micro portion. The micro portion is obtained by applying to the original kinetic equation a projector into the orthogonal complement of the null space of the collision operator -- this projector depends on the macro portion. In particular, we extend the technique of Bennoune, Lemou, and Mieussens [{\it Uniformly stable schemes for the Boltzmann equation preserving the compressible Navier-Stokes asymptotics, J. Comput. Phys. (2008)}] to two-space dimensions, the ES-BGK collision operator, and problems with reflecting wall boundary conditions. The collision operator in the micro and macro equations is handled via L-stable implicit time discretizations, while the transport terms are computed via kinetic flux vector splitting (for the macro equations) and upwind differencing (for the micro equation). The resulting scheme is applied to various test cases in 1D and 2D. The 2D version of the code is parallelized using MPI, and we present weak- and strong-scaling studies with varying numbers of processors.
\end{abstract}


\maketitle

\tableofcontents


\section{Introduction}
\label{sec:intro}
A rarefied gas is a collection of unbound molecules in which the mean free path (i.e., the average distance a molecule will travel before colliding with another molecule) is comparable or larger than some characteristic length scale in the problem; the ratio of these two lengths (i.e., the mean free path divided by the characteristic length scale) is called the {\it Knudsen number}, which we denote by $\varepsilon>0$. Different flow regimes can be enumerated through the size of the Knudsen number; a common approach is to define the following four flow regimes: (1) continuum $\left(\varepsilon<10^{-2}\right)$, (2) slip $\left(10^{-2}<\varepsilon<10^{-1}\right)$, (3) transition $\left(10^{-1}<\varepsilon<10\right)$, and (4) free molecular $\left(\varepsilon>10\right)$. A gas is typically considered {\it rarefied} in all but the continuum regime. Applications in which rarefied gases are essential in the overall dynamics include flow around high-altitude aircraft, various propulsion mechanisms for spacecraft, and flow through microfluidic and nanofluidic devices \cite{book:karniadakis2005microflows}.

To accurately predict rarefied gas dynamics (RGD), continuum fluid models such as the Navier-Stokes equations are insufficient; instead, kinetic transport models, such as the Boltzmann equation, are typically required. Kinetic Boltzmann equations are ordinarily more challenging to solve computationally than continuum fluid models, since they require solving integro-differential equations over the full phase space and discretizing potentially complex collision operators. We briefly review some standard approaches for solving kinetic equations below. 
\begin{description}
\item {\bf Direct Simulation Monte Carlo (DSMC).} DSMC (e.g., see Bird \citep{Bird1970DirectEquation,Bird1994MolecularEdition}, Nanbu \cite{article:Nanbu1983}, and
Pareschi and Russo \cite{Pareschi2001AnEquation}) is a classical approach for solving the Boltzmann kinetic equation that uses a Lagrangian description of gas particles and models particle collisions using a probabilistic Monte Carlo method. DSMC is efficient for large Knudsen numbers (i.e., large mean-free paths), but the computational cost increases as the Knudsen number decreases. Another disadvantage of DSMC is that the numerical results will contain statistical noise resulting from the Monte Carlo treatment of the collision operator, which decreases relatively slowly with the number of particles, $N$: ${\mathcal O}\left(N^{-1/2}\right)$.
\item {\bf Direct Kinetic Methods.} An alternative to the DSMC approach is to {\it directly} discretize the kinetic equation, meaning that all terms in the equations are approximated deterministically using a variety of finite difference, finite volume, finite element, and spectral methods. Several such approaches exist in the literature, including Chorin \cite{article:Chorin1972}, Pareschi and Russo \citep{Pareschi2000NumericalOperator},
Pieraccini and Puppo \cite{Pieraccini2007Implicit-explicitEquations}, Xu and Huang \cite{article:XuHuang2010}, Filbet and Jin \cite{Filbet2011AnEquation},
G\"u\c{c}l\"u and Hitchon \cite{article:Guclu2012}, Kitzler and Sch\"oberl \cite{article:Kitzler2015}, Su et al. \cite{Su2020ImplicitEquation}, Li et al. \cite{article:LiFang2022}, as well as many others. Most direct kinetic methods are applied to simplified collision operators like BGK \cite{Bhatnagar1954ASystems} or ES-BGK \cite{article:Holway1965}. The advantage of direct kinetic models over DSMC is that they do not suffer from the statistical noise of Monte Carlo methods and can remain relatively efficient even for small Knudsen numbers with appropriate implicit or semi-implicit time-stepping schemes. However, the main drawback of such methods is that they require discretization of the full phase space, which, due to the curse of dimensionality (i.e., exponential growth in computational cost with the number of dimensions), can result in high computational costs. 
\item {\bf Moment-Closure Methods.} An alternative to direct kinetic methods is to replace the full particle density function (PDF) by a small number of moments of the PDF (e.g., mass, momentum, and energy). The dynamics of these moments occur on a lower-dimensional manifold than the full PDF (e.g., 3D rather than 6D), significantly reducing computational complexity. Furthermore, moments of the collision operator are typically simpler than the full collision operator. The difficulty with this approach is determining an accurate and efficient moment-closure that allows one to replace the full PDF with a small number of moments. The first systematic moment-closure approach was developed by Grad \cite{Grad1949OnGases} in 1949; in subsequent decades, many alternative moment-closures have been developed, including maximum entropy closures \cite{article:Dreyer87,Levermore1996MomentTheories,article:MuRu93}, quadrature-based moment closures \cite{article:desjardins2008,article:Fox08}, the R13 closure \cite{Struchtrup2003RegularizationAnalysis}, and many others. Torrilhon \cite{article:Torrilhon2009} and the references therein review many of these approaches. The advantages of these methods are that they greatly reduce computational complexity and directly capture the small-Knudsen-number limit. The disadvantage of these approaches is that they often produce accurate solutions only near thermodynamic equilibrium and therefore may not provide accurate solutions if the Knudsen number becomes too large.
\end{description}

\begin{figure}[!th]
\begin{center}
\begin{tikzpicture}[scale=.8, transform shape]
\draw [very thick] [-]  (-2.25,0) -- (3.25,0);
\draw [very thick] [-]  (3.25,1.5) -- (3.25,0);
\draw [very thick] [-]  (3.25,1.5) -- (-2.25,1.5);
\draw [very thick] [-]  (-2.25,1.5) -- (-2.25,0);
\small
\node[align=center] at (0.5,1.0) {Numerical method for fluid model:};
\node[align=center] at (0.5,.5) {$\vec{q}_{,t} + \nabla_{\vec{x}} \cdot \mat{\bf F}\left(\vec{q}\right) = \vec{0}$};
\normalsize
\draw [very thick] [->]  (3.5,0.75) -- (5.75,0.75);
\node[align=center, above] at (4.75,0.75) {$\begin{bmatrix} \Delta t \\ \Delta x  
\end{bmatrix}\rightarrow 0^+$};
\draw [very thick] [-]  (6,0) -- (11,0);
\draw [very thick] [-]  (11,1.5) -- (11,0);
\draw [very thick] [-]  (11,1.5) -- (6,1.5);
\draw [very thick] [-]  (6,1.5) -- (6,0);
\small
\node[align=center] at (8.5,1.0) {Exact solution of fluid model:};
\node[align=center] at (8.5,.5) {$\vec{q}_{,t} + \nabla_{\vec{x}} \cdot \mat{\bf F}\left(\vec{q}\right) = \vec{0}$};
\normalsize
\draw [very thick] [-]  (-2.25,3.25) -- (-2.25,4.75);
\draw [very thick] [-]  (-2.25,4.75) -- ( 3.25,4.75);
\draw [very thick] [-]  ( 3.25,4.75) -- ( 3.25,3.25);
\draw [very thick] [-]  ( 3.25,3.25) -- (-2.25,3.25);
\small
\node[align=center] at (0.5,4.25) {Numerical method for kinetic model:};
\node[align=center] at (0.5,3.75) {$f_{,t}+\vec{v} \cdot \nabla_{\vec{x}} f 
= \frac{1}{\varepsilon} C(f)$};
\normalsize
\draw [very thick] [->]  (3.5,4.0) -- (5.75,4.0);
\node[align=center, above] at (4.75,4.0) {$\begin{bmatrix} \Delta t \\ \Delta x \\ \Delta v 
\end{bmatrix}\rightarrow 0^+$};
\draw [very thick] [-]  (6,3.25) -- (6,4.75);
\draw [very thick] [-]  (6,4.75) -- (11,4.75);
\draw [very thick] [-]  (11,4.75) -- (11,3.25);
\draw [very thick] [-]  (11,3.25) -- (6,3.25);
\small
\node[align=center] at (8.5,4.25) {Exact solution of kinetic model:};
\node[align=center] at (8.5,3.75) {$f_{,t}+\vec{v} \cdot \nabla_{\vec{x}} f 
= \frac{1}{\varepsilon} C(f)$};
\draw [very thick] [->]  (0.5,3.0) -- (0.5,1.75);
\draw [very thick] [->]  (8.5,3.0) -- (8.5,1.75);
\node[align=center, right] at (0.75,2.5) {$\varepsilon \rightarrow 0^+$};
\node[align=center, left] at (8.25,2.5) {$\varepsilon \rightarrow 0^+$};
\normalsize
\end{tikzpicture}
\caption{(\S \ref{sec:intro}) Commutative diagram that defines the asymptotic-preserving (AP) property. If we start with some numerical method for the kinetic model (top left), then we can imagine taking two different limits: (1) the limit as the mesh parameters vanish ($\Delta t$, $\Delta x$, $\Delta v$ $\rightarrow 0^+$), and (2) the limit as the Knudsen number vanishes ($\varepsilon \rightarrow 0^+$). If these limits commute, we say the method is asymptotic preserving (AP). The significance of this is that if a numerical method for the kinetic model (top left) is AP, it produces stable and accurate solutions on a fixed mesh for any value of $\varepsilon>0$. \label{fig:apscheme}}
\end{center}
\end{figure}

\subsection{Asymptotic Preserving Property}
The so-called {\it fluid limit} of the kinetic Boltzmann occurs when the Knudsen number tends to zero (from the right): $\varepsilon \rightarrow 0^+$. This limit is of vital importance in gas dynamics, since it represents the instantaneous relaxation to thermodynamic equilibrium. The difficulty is that this limit is also a {\it singular limit} of the kinetic Boltzmann equation, since certain high-order derivative terms (i.e., viscosity and thermal diffusion) vanish in that limit. In complex flows, parts of the domain may be in or very near thermodynamic equilibrium, while other portions may not; thus, accurately and efficiently handling the full range of Knudsen numbers is critical.

The key property that allows a numerical scheme to remain stable and accurate in the singular limit, $\varepsilon \rightarrow 0^+$, is commonly referred to as the {\it asymptotic preserving} (AP) property (e.g., see
\cite{Caflisch1997UniformlyRelaxation,Coron1991NumericalEquations,Gabetta1997RelaxationEquations,article:Jin1995,Jin1999EfficientEquations,Jin2012AsymptoticReview,Jin1996NumericalTerms}). The basic idea is this: given a stable and consistent numerical method for the kinetic equation, it should be possible to fix the mesh parameters (i.e., the time-step and mesh spacing), and then take the limit $\varepsilon \rightarrow 0^+$ in the numerical method and arrive at a new stable and consistent numerical method for the system in thermodynamic equilibrium (i.e., compressible Euler equations). If this is possible, we refer to the original numerical method for the kinetic equation as an asymptotic-preserving (AP) method. The key mathematical property that must be met is this: the limit of $\varepsilon \rightarrow 0^+$ and the limit of the mesh parameters tending to zero (again, from the right) must {\it commute}. A commutative diagram illustrating this point is shown in  Figure \ref{fig:apscheme}. To achieve the AP property in practice, one must use an implicit or semi-implicit time discretization for the collision operator.

\subsection{Scope of this work}
One useful tool in solving the Boltzmann equation is to decompose the distribution function into two parts: (1) the equilibrium distribution (i.e., Maxwell-Boltzmann), plus (2) the non-equilibrium portion. Coupled evolution equations can be formulated for each portion. Liu and Hu \cite{article:LiuHu2004} introduced such an approach and called it a {\bf micro-macro decomposition}. Numerical methods based on this micro-macro decomposition idea were developed by Bennoune, Lemou, and Mieussens \cite{Bennoune2008UniformlyAsymptotics}, with additional modifications developed by several authors, including Crestetto, Crouseilles, and Lemou \cite{article:Crestetto2012}, Xiong, Jang, Li, and Qiu \cite{XIONG201570}, Gamba, Jin, and Liu \cite{GAMBA2019264}, and Endeve and Hauck \cite{ENDEVE2022111227}. 

The basic idea is to write the full particle distribution function (PDF) as a sum of its thermodynamic equilibrium (i.e., the macro portion) plus everything else (i.e., the micro portion). The micro-macro methods of \cite{Bennoune2008UniformlyAsymptotics,article:Crestetto2012,XIONG201570,GAMBA2019264,ENDEVE2022111227} can be thought of as direct kinetic methods, in that they solve for the entire distribution on a mesh; however, they also borrow from moment-closure methods in that they directly evolve fluid moments (i.e., the macro portion). Unlike ``standard'' direct kinetic methods, micro-macro methods discretize only the non-equilibrium portion of the PDF on the full phase-space mesh, allowing velocity grids with smaller velocity extents and lower resolutions than in comparable direct kinetic methods. The benefit of the micro-macro methods over moment-closure methods is that one can accurately simulate flows at any Knudsen number, whether arbitrarily close to or far from thermodynamic equilibrium.

This work aims to extend the micro-macro decomposition method of Bennoune, Lemou, and Mieussens \cite{Bennoune2008UniformlyAsymptotics} in several key directions. We show how to extend the micro-macro decomposition to multiple dimensions and to accurately handle the Ellipsoidal-Statistical-Bhatnagar-Gross-Krook (ES-BGK) \cite{article:Holway1965} collision operator. One disadvantage of the micro-macro decomposition approach is that boundary conditions can be much more difficult to derive and implement than in direct kinetic methods. In this work, we show how to implement reflecting wall boundary conditions in both the macro and micro portions of the numerical update to recover solutions comparable to those computed with direct kinetic methods. In the proposed numerical scheme, the collision operator, as it appears in both the micro and macro equations, is handled via L-stable implicit time discretizations. The remaining transport terms are computed via kinetic flux vector splitting (for the macro equations) and upwind differencing (for the micro equation). The semi-implicit nature of the full update makes the scheme asymptotic preserving (AP). The resulting scheme is applied to various test cases in 1D and 2D. The 2D version of the code is parallelized using MPI, and we present weak- and strong-scaling studies with varying numbers of processors. The resulting methods are written in an open-source C++ code with MPI for parallelization and Python plotting routines; this code can be freely downloaded as a {\tt bitbucket.org} repository \cite{code:RossmanithSar2025}.

The remainder of this paper is organized as follows. In \S\ref{sec:boltzmann} we describe the kinetic Boltzmann-BGK and Boltzmann-ES-BGK models and review their properties. We develop in \S\ref{sec:micmac1D} the one-dimensional version of the proposed micro-macro decomposition numerical method and provide detailed descriptions of the numerical scheme, both in the macro and micro portions. Numerical examples of the resulting scheme are presented in \S\ref{sec:numerical1D}, including a convergence study, a shock-tube problem, and a boundary-driven flow that shows the robustness of the scheme across a wide range of Knudsen numbers. In \S\ref{sec:micmac2D} we extend the proposed method to the 2D Boltzmann-ES-BGK model and again provide detailed descriptions of the numerical scheme's macro and micro portions. Numerical examples of the resulting scheme are presented in \S\ref{sec:numerical2D}, including a convergence study, a cylindrical shock problem, and the results of a lid-driven cavity example in the slip/transition regime that clearly shows dynamics beyond the Navier-Stokes-Fourier approximation. We also briefly describe a parallel MPI implementation in the 2D case and provide both weak- and strong-scaling analyses. Finally, we summarize our findings in \S\ref{sec:conclusions}.


\section{Boltzmann Equation and Approximate Collision Operators}
\label{sec:boltzmann}
We begin with a brief review of the kinetic Boltzmann equation for a monoatomic gas, the BGK and ES-BGK approximate
collision operators, as well as the fluid equations that result in the strongly collisional limit.

\subsection{Boltzmann Equation}
The dimensionless form of the single-species Boltzmann equation can be written as:
\begin{equation}
    \label{eq2-1}
    \frac{\partial f}{\partial t} + \vec{v}\cdot\nabla_{\vec{x}}  f = \frac{1}{\varepsilon} \, \mathcal{Q}(f),
\end{equation}
where $f\left(t,\vec{x},\vec{v}\right)\!: \reals_{\ge 0} \times \reals^{d_x} \times \reals^{d_v} \mapsto \reals_{\ge 0}$ is the particle density function, $t \in \reals_{\ge 0}$ is time, $\vec{x} \in \reals^{d_x}$ is the spatial coordinate, and $\vec{v} \in \reals^{d_v}$ is the velocity coordinate, and $d_x, d_v \in \left\{1,2,3\right\}$ are the spatial and velocity dimensions, respectively. The parameter $\varepsilon \in \reals_{>0}$ is the Knudsen number, a nondimensional number that compares the mean free path of particles to a characteristic length scale; small (large) Knudsen numbers mean high (low) rates of collisions. The collision operator, $\mathcal{Q}(f): \reals_{\ge 0} \mapsto \reals$, models inter-particle collisions.

We use the following shorthand to denote integration over velocity:
\begin{equation}
\left\langle g \right\rangle := \int_{\mathbb{R}^{d_v}} g\left( \vec{v} \right) \, d\vec{v}.
\end{equation}
The first $d_v+2$ moments of the distribution function, $f$, are then given by
\begin{equation}
\label{eqn2-3}
    \vec{U} := \left( \rho, \, \rho \vec{u} \, , \, {\mathcal E} \right) :=  \left\langle \vec{m}\left( \vec{v} \right) f \right\rangle \qquad \text{where} \qquad \vec{m}\left( \vec{v} \right) := \left(1, \, \vec{v} \, , \, \frac{\|\vec{v}\|^2}{2} \right),
\end{equation}
and $\rho\left(t,\vec{x} \right): \reals_{\ge 0} \times \reals^{d_x} \mapsto \reals_{\ge 0}$ is the mass density,
$\rho \vec{u} \left(t,\vec{x} \right): \reals_{\ge 0} \times \reals^{d_x} \mapsto \reals^{d_v}$ is the momentum density,
and ${\mathcal E}\left(t,\vec{x} \right): \reals_{\ge 0} \times \reals^{d_x} \mapsto \reals_{\ge 0}$ is the scalar energy density.
The following variables can be directly computed from the first $d_v+2$ moments:
\begin{align}
\text{macroscopic velocity:} & \quad \vec{u} = \frac{\rho \vec{u}}{\rho} = \frac{1}{\rho} \int_{\mathbb{R}^{d_v}} \vec{v} \, f \, d\vec{v} \, , \\
\label{eqn:scalar-pressure}
\text{scalar pressure:} & \quad  p = \frac{2 {\mathcal E} -   \rho \| \vec{u} \|^2}{d_v}  = 
\frac{1}{d_v} \int_{\mathbb{R}^{d_v}} \| \vec{v} - \vec{u} \|^2 \, f \, d\vec{v} \, , \\
\label{eqn:scalar-temperature}
 \text{scalar temperature:} & \quad  T = \frac{p}{\rho} = \frac{1}{\rho  d_v} \int_{\mathbb{R}^{d_v}} \| \vec{v} - \vec{u} \|^2 \, f \, d\vec{v} \, .
\end{align}

The collision operator, $\mathcal{Q}(f)$, can take different forms (e.g., see Villani \cite{article:Villani2002}), but in each version
should obey the following key properties:
\begin{itemize}
    \item Conservation of mass, momentum, and energy:
    \begin{equation}
    \label{eqn:collision-conservation}
       \left\langle \vec{m}\left(\vec{v} \right) \, \mathcal{Q}(f) \right\rangle = 0, \quad \forall f \geq 0;
    \end{equation}
    \item Entropy inequality:
    \begin{equation}
        \label{eqn:collision-entropy}
\left\langle \mathcal{Q}(f) \, \log(f) \right\rangle \leq 0, \quad \forall f \geq 0;
    \end{equation}
    \item The Maxwell-Boltzmann distribution is the equilibrium distribution:
    \begin{equation}
    \label{eqn:maxwboltz}
        \mathcal{Q}(f) \equiv 0 \quad \Longleftrightarrow \quad f \equiv \mathcal{M}\left[f\right] :=  \frac{\rho}{\left({2\pi T}\right)^{d_v/2}}\exp\left(-\frac{\|\vec{v}-\vec{u}\|^2}{2T}\right),
    \end{equation}
where $\rho, \vec{u}$, and $T$ are the density, macroscopic velocity, and scalar temperature associated with $f$:
\begin{equation}
\label{eq2-4}
\left\langle \vec{m}\left(\vec{v} \right) \, \mathcal{M}[f] \right\rangle \equiv \left\langle \vec{m}\left(\vec{v} \right) \, f \right\rangle.
\end{equation}
\end{itemize}

\subsection{Bhatnagar-Gross-Krook (BGK) Collision Operator} 
To avoid the complexity of the full Boltzmann collision operator, approximate collision operators are often employed in numerical simulations that mimic the behavior of the Boltzmann operator, but that are much less computationally expensive to evaluate. While there exist multiple possible approximations, the simplest approximate collision operator that satisfies all of the properties outlined above in \eqref{eqn:collision-conservation}, \eqref{eqn:collision-entropy}, and
\eqref{eqn:maxwboltz}, is the Bhatnagar-Gross-Krook (BGK) operator \citep{Bhatnagar1954ASystems}, which can be written as 
\begin{equation}
\label{eqn:collision-bgk}
\mathcal{Q}(f) := {\tau} \bigl( \mathcal{M}\left[f\right] - f \bigr),
\end{equation}
where $\mathcal{M}\left[f\right]$ is the Maxwell-Boltzmann distribution \eqref{eqn:maxwboltz},
and $\tau$ is a collision parameter that allows us to make the collision frequency solution-dependent
\citep{Coron1991NumericalEquations,Mieussens2000DiscreteDynamics}.

One way to quantify the effect of the BGK collision operator is to understand the dynamics of the Boltzmann-BGK equation in the limit $\varepsilon \rightarrow 0^+$. This can be done by applying a power series expansion for $f$ in the Knudsen number, $\varepsilon$, which is commonly referred to as a Chapman-Enskog expansion (e.g., see \cite{book:Cercignani2000}), of the following form:
\begin{equation}
\label{eqn:chapman-enskog}
   f = \mathcal{M}\left[f\right] \left( 1 + \varepsilon g_{1} + \varepsilon^2 g_{2} + \cdots \right),
\end{equation}
which results to ${\mathcal O}\left( \varepsilon^2 \right)$ in the Navier-Stokes-Fourier system:
\begin{equation}
\label{eqn:NSF}
    \partial_t\begin{pmatrix}
        \rho \\ \rho \vec{u} \\ {\mathcal E}
    \end{pmatrix} + \nabla_\vec{x}\cdot\begin{pmatrix}
        \rho \vec{u} \\ \rho \vec{u} \otimes \vec{u} + p \mat{\mathcal{I}} \\
        \left({\mathcal E}+p\right) \vec{u}
    \end{pmatrix} = \varepsilon \nabla_\vec{x} \cdot \begin{pmatrix}
        0 \\ \mu \, \mat{\sigma} \\ \, \mu \, \mat{\sigma}\;\vec{u} - \vec{h_1} \,
    \end{pmatrix} + {\mathcal O}\left( \varepsilon^2 \right),
\end{equation}
where $\mat{\mathcal{I}}$ is the $d_v \times d_v$ identity matrix, and
\begin{align}
\label{eqn:NSF-sigma-q}
    \mat{\sigma} = \left( \nabla_\vec{x} \, \vec{u} + \left(\nabla_\vec{x} \, \vec{u}\right)^T \right) - \frac{2}{d_v} \left( \, \nabla_\vec{x} \cdot \vec{u} \, \right) \mat{\mathcal{I}} \, \, \, \, \left(\text{stress tensor}\right), \quad
    \vec{h_1} = -\kappa \nabla_{\vec{x}} T \, \, \, \, \left({\mathcal O}\left(\varepsilon\right) \text{heat flux vector}\right). 
\end{align}
In \eqref{eqn:NSF}, the transport coefficients, $\mu$ and $\kappa$, are functions of $\vec{U}$:
\begin{align}
\label{eqn:BGK-transport-coeffs}
  \mu = \frac{p}{\tau} \, \, \, \,  \left(\text{viscosity}\right), \quad 
  \kappa = \left(\frac{d_v + 2}{2}\right) \frac{p}{\tau} \, \, \, \,   \left(\text{heat conductivity}\right), \quad
     \text{Pr} = \left(\frac{d_v + 2}{2}\right) 
       \frac{\mu}{\kappa} = 1 \, \, \, \,   \left(\text{Prandtl \#}\right).
\end{align}
Note that the Prandtl number, which is proportional to the ratio of the viscosity to the heat conductivity, is always one in the BGK model.

\subsection{ES-BGK Collision Operator}
The BGK model has the correct compressible Euler limit when $\varepsilon\rightarrow0^+$ and satisfies
properties \eqref{eqn:collision-conservation}, \eqref{eqn:collision-entropy}, and \eqref{eqn:maxwboltz}. However,  the ${\mathcal O}\left( \varepsilon^2 \right)$ Navier-Stokes-Fourier system \eqref{eqn:NSF}--\eqref{eqn:NSF-sigma-q} has incorrect transport coefficients \eqref{eqn:BGK-transport-coeffs} compared to the Boltzmann collision operator; and in particular, the Prandtl number is always one, while for most gases it should be less than one.

To overcome this problem, Holway \citep{article:Holway1965} proposed the so-called ellipsoidal statistical model (ES-BGK), which through the Chapman-Enskog expansion yields the same Navier-Stokes-Fourier approximation \eqref{eqn:NSF}--\eqref{eqn:NSF-sigma-q}, but with slightly modified transport coefficients, including a tuneable Prandtl number that can be made less than one. The ES-BGK collision operator is similar to the original BGK collision operator. However, the Maxwell-Boltzmann distribution in the relaxation term of \eqref{eqn:collision-bgk} is replaced by an anisotropic Gaussian.
For ES-BGK, conservation property \eqref{eqn:collision-conservation} can be readily shown, as can the equilibrium distribution property \eqref{eqn:maxwboltz}. Entropy property \eqref{eqn:collision-entropy} is more difficult to prove, but was shown by Andries et al. \cite{Andries2000Gaussian-BGKNumber,article:Andries2001} to be satisfied for 
\begin{equation}
\frac{d_v-1}{d_v} \le \text{Pr} < \infty.
\end{equation}

To describe the ES-BGK collision operator, we first define the temperature tensor:
\begin{equation}
\label{eqn:temp_tensor}
    \mat{\mathbb T} = \frac{1}{\rho}\int_{\mathbb{R}^{d_v}} \left(\vec{v} - \vec{u}\right) \otimes \left(\vec{v}-\vec{u}\right) \, f \, dv \, , \quad   T = \frac{1}{d_v} \, \text{tr}\left(\mat{\mathbb T}\right),
\end{equation}
where the trace of the temperature tensor is proportional to the scalar temperature $T$.
Next we introduce the {\it modified temperature tensor}:
\begin{equation}
    \label{eq2-19}
   \mat{\mathcal{T}} =  \left(1-\nu \right) \, T \, \mat{\mathcal{I}} + \nu \, \mat{{\mathbb T}} \, ,
   \quad \text{where} \quad  
   \frac{1}{1-d_v} \le \nu < 1 \, .
\end{equation}
Note that in the case $d_v=1$, the range on $\nu$ is given by $-\infty < \nu < 1$. Also note that
$\text{tr}\left(\mat{\mathcal{T}}\right) = \text{tr}\left(\mat{{\mathbb T}}\right)$.
Using the modified temperature tensor, the ES-BGK collision operator can be written as
\begin{equation}
\label{eq2-15}
    \mathcal{Q}(f) = {\tau} \left(\mathcal{G}[f] - f \right) \, ,
\end{equation}
where the Gaussian distribution, $\mathcal{G}[f]$, is defined through the modified temperature tensor:
\begin{equation}
\label{eqn2-14}
    \mathcal{G}[f] = \frac{\rho}{\sqrt{\text{det}\left(2 \, \pi \, \mat{\mathcal{T}} \, \right)}}\exp\left(-\frac{1}{2} {\left( \, \vec{v} - \vec{u} \, \right) \cdot \mat{\mathcal{T}}^{-1} \left( \, \vec{v}-\vec{u} \, \right)}\right).
\end{equation}
The first $d_v+2$ moments, as defined in \eqref{eqn2-3}, of $\mathcal{G}[f]$ again match those of $f$:
\begin{equation}
\label{eq:moments_ESBGK_G}
\left\langle \vec{m}\left(\vec{v} \right) \, \mathcal{G}[f] \right\rangle \equiv \left\langle \vec{m}\left(\vec{v} \right) \, f \right\rangle,
\end{equation}
but the temperature tensors of each distribution are generally not the same:
\begin{align}
    \mat{{\mathbb T}} = \frac{1}{\rho} \int_{\reals^{d_v}} \left( \vec{v}-\vec{u} \right) \otimes \left( \vec{v}-\vec{u} \right) \, f \, d\vec{v}  \qquad \text{and} \qquad
    \mat{\mathcal{T}} = \frac{1}{\rho} \int_{\reals^{d_v}} \left( \vec{v}-\vec{u} \right) \otimes \left( \vec{v}-\vec{u} \right) \, \mathcal{G}[f] \, d\vec{v} \, .
\end{align}
Similarly, the pressure tensors of each distribution are generally not the same:
\begin{align}
\label{eqn:pressure-tensor-esbgk}
     \mat{{\mathbb P}} := \int_{\reals^{d_v}} \left( \vec{v}-\vec{u} \right) \otimes \left( \vec{v}-\vec{u} \right) \, f \, d\vec{v} = \rho \mat{{\mathbb T}}  \qquad \text{and} \qquad
    \mat{\mathcal{P}} := \int_{\reals^{d_v}} \left( \vec{v}-\vec{u} \right) \otimes \left( \vec{v}-\vec{u} \right) \, \mathcal{G}[f] \, d\vec{v} = \rho \mat{\mathcal T}.
\end{align}
The scalar pressure, $p$, and temperature, $T$, defined in \eqref{eqn:scalar-pressure} and \eqref{eqn:scalar-temperature} satisfy:
\begin{equation}
p = \frac{1}{d_v} \text{tr}\left(\mat{\mathcal{P}}\right) = \frac{1}{d_v} \text{tr}\left(\mat{\mathbb{P}}\right)
\qquad \text{and} \qquad
T = \frac{1}{d_v} \text{tr}\left(\mat{\mathcal{T}}\right) = \frac{1}{d_v} \text{tr}\left(\mat{{\mathbb T}}\right) = 
\frac{p}{\rho}.
\end{equation}
Another important related tensor is the energy tensor:
\begin{align}
\label{eqn:energy-tensor-esbgk}
     \mat{{\mathbb E}} := \int_{\reals^{d_v}} \vec{v} \otimes \vec{v} \, f \, d\vec{v} = \rho \left( \, \vec{u} \otimes \vec{u} + \mat{{\mathbb T}} \, \right) \qquad \text{where} \qquad
     {\mathcal E} = \frac{1}{2} \text{tr}\left(\mat{\mathbb{E}}\right) = 
     \frac{1}{2} \rho \left( \, \| \vec{u} \|^2 + d_v \, T \, \right).
\end{align}

Again using the Chapman-Enskog expansion \eqref{eqn:chapman-enskog}, one arrives at the same Navier-Stokes-Fourier equations \eqref{eqn:NSF}--\eqref{eqn:NSF-sigma-q}, but this time with modified transport coefficients given by
\citep{article:Andries2001,Andries2000Gaussian-BGKNumber,Filbet2011AnEquation}:
\begin{align}
\label{eqn:ESBGK-transport-coeffs}
  \mu =  \frac{p}{\left(1-\nu\right)\tau} \, \, \,  \left(\text{viscosity}\right), \, \, \,
  \kappa = \left(\frac{d_v + 2}{2}\right) \frac{p}{\tau} \, \, \, \,   \left(\text{heat cond.}\right), \, \, \,
     \text{Pr} = \left(\frac{d_v + 2}{2}\right) 
       \frac{\mu}{\kappa} = \frac{1}{1-\nu} \, \, \, \left(\text{Prandtl \#}\right).
\end{align}
Note that the Prandtl number is now a function of the parameter $\nu$ and satisfies:
\begin{equation}
\frac{d_v-1}{d_v} \le \text{Pr} = \frac{1}{1-\nu} < \infty \qquad \text{for} \qquad \frac{1}{1-d_v} \le \nu < 1.
\end{equation}

\begin{remark}
According to Chapman \citep{Chapman1962TheGases}, the Prandtl number and viscosity computed for the full Boltzmann equation for a hard sphere model for a monoatomic gas with $d_v=2$ is given by
\begin{equation}
    \text{Pr} \approx \frac{1}{2} \qquad \text{and} \qquad \mu_B(T) = \frac{\sqrt{2}}{3\pi}\frac{T}{C_{25}} \qquad \text{where} \qquad C_{25} \approx 0.436.
\end{equation}
The value for $\nu$ and the expression $\tau$ need to be chosen such that the viscosity, $\mu$, and heat conductivity, $\kappa$, computed from the asymptotic limit of the ES-BGK model are the same as the ones corresponding to the full Boltzmann equation given above. 
Therefore, setting values for $\nu$ and $\tau$ to match these values yields \cite{Filbet2011AnEquation}:
\begin{equation}
\label{eqn:tau-for-2d-es-bgk}
    \nu = -1 \qquad \text{and} \qquad 
    \tau = \frac{p}{2\mu_B(T)} = \left(\frac{3 C_{25} \, \pi}{2\sqrt{2}} \right) \rho \approx \left(0.4625 \, \pi \right) \rho.
\end{equation}
\end{remark}

\section{1D Micro-Macro Decomposition}
\label{sec:micmac1D}
The micro-macro decomposition method proposed in this work is based on the original approach of Bennoune, Lemou, and Mieussens \citep{Bennoune2008UniformlyAsymptotics} for solving the one-dimensional (i.e., $d_v=1$ and $d_x=1$) BGK model. In this section, we review the approach of
\citep{Bennoune2008UniformlyAsymptotics}, but with two important changes that will be propagated to the multi-dimensional extension (i.e., $d_v=2$ and $d_x=2$) for the ES-BGK model in Section \ref{sec:micmac2D}: (1) we do not stagger the microscopic portion of the distribution function relative to the macroscopic fluid variables, and (2) the numerical update for the macroscopic fluid variables is carried out using a kinetic flux vector splitting (KFVS) approach. While neither of these changes is critical in 1D, they will allow us to directly extend to 2D and handle non-trivial boundary conditions.

\subsection{Mathematical Formulation}
The one-dimensional ($d_v=1$ and $d_x=1$) Boltzmann-BGK equation can be written as:
\begin{equation}
\label{eqn:boltzmann-bgk-1d}
f_{,t} + {v}  f_{,x} = \frac{\tau}{\varepsilon} \left( \maxw\left[f\right] - f \right)  \quad
\text{where} \quad
\maxw\left[f\right] = \frac{\rho}{\sqrt{2\pi T}} \exp\left[
-\frac{\left(v - u \right)^2}{2T} \right],
\end{equation}
and
\begin{equation}
\rho = \left\langle f \right\rangle = \left\langle \maxw\left[f\right] \right\rangle, \quad
u = \frac{1}{\rho} \left\langle v \, f \right\rangle = \frac{1}{\rho} \left\langle v \, \maxw\left[f\right] \right\rangle, \quad T = \frac{1}{\rho} \left\langle \left(v - u \right)^2 \, f \right\rangle = \frac{1}{\rho} \left\langle \left(v - u \right)^2 \, \maxw\left[f\right] \right\rangle.
\end{equation}
In \eqref{eqn:boltzmann-bgk-1d}, we have made use of the following shorthand notation for partial differentiation:
\begin{equation}
f_{,t} := \frac{\partial f}{\partial t} \quad \text{and} \quad
f_{,x} := \frac{\partial f}{\partial x}.
\end{equation}
One of the key steps in the micro-macro approach \citep{Bennoune2008UniformlyAsymptotics} is to decompose the full distribution function, $f$, into a macroscopic Maxwell-Boltzmann distribution, $\mathcal{M}[f]$, plus whatever is left over, which is called the microscopic portion and denoted by $g$: 
\begin{equation}
\label{eq5-1}
f  = \maxw\left[f\right]  + \varepsilon \, g \qquad \text{where} \qquad \varepsilon := \text{Knudsen number}.
\end{equation}
We note that although this looks like a perturbation expansion, it is, in fact, not a perturbation expansion, because we do not assume that $\varepsilon$ is small, nor is anything being truncated. Also note that because $f$ and $\maxw[f]$ have the same first three moments, the first three moments of the microscopic portion, $g$, must be zero:
\begin{equation}
\left\langle \vec{m}\left(v\right) g \right\rangle = 0 \qquad \text{where} \qquad
\vec{m}\left(v\right) = \left( 1, \, v, \, \frac{1}{2} v^2 \right).
\end{equation}
Inserting the micro–macro decomposition \eqref{eq5-1} into the Boltzmann-BGK equation yields:
\begin{equation}
\label{eq5-2}
\maxw_{,t} + {v}  \maxw_{,x} + \varepsilon g_{,t}  +\varepsilon \, v g_{,x} =  - \tau \, g.
\end{equation}
Computing the first three moments of \eqref{eq5-2} as defined by \eqref{eqn2-3} with $d_v=1$, yields the following fluid system:
\begin{equation}
\left(\text{1D Macro}\right): \quad
\label{eq5-11}
\frac{\partial}{\partial t}
\begin{bmatrix}
\rho \\
\rho u \\
{\mathcal E}
\end{bmatrix} + 
\frac{\partial}{\partial x}
\begin{bmatrix}
\rho u \\
\rho u^2 + \rho T \\
u \left( {\mathcal E} + \rho T \right)
\end{bmatrix} = \frac{\partial}{\partial x}
\begin{bmatrix}
0 \\ 0 \\ - h
\end{bmatrix}, \quad
{\mathcal E} = \frac{1}{2} \rho T + \frac{1}{2} \rho u^2, \quad 
h = \frac{1}{2} \varepsilon \left\langle v^3 g \right\rangle,
\end{equation}
which is not closed since the heat flux, $h$, must still be supplied via the microscopic portion $g$.

To close fluid system \eqref{eq5-11}, we need to obtain an evolution equation for the microscopic portion, $g$, from equation \eqref{eq5-2}. One key observation is that the Boltzmann-BGK collision operator \eqref{eqn:collision-bgk} can be rewritten using micro-macro decomposition \eqref{eq5-1}:
\begin{equation}
  \mathcal{Q}(f) = {\tau} \bigl( \mathcal{M}\left[f\right] - f \bigr) \quad \Longrightarrow \quad
  \mathcal{Q}(f) = - {\tau} \varepsilon g,
\end{equation}
which shows that the BGK collision operator is linear in $g$, since $\tau$ depends only on macroscopic quantities. Furthermore, note that the first three moments of $g$ vanish, which means that the nullspace of the BGK collision operator is the set of all functions in the span of the following basis:
\begin{equation}
\label{eq5-3}
    \mathcal{B} = \left\{ \, \frac{1}{\rho}\mathcal{M}, \quad \left( \frac{v-u}{\sqrt{T}} \right) \frac{1}{\rho}\mathcal{M}, \quad \left(\frac{\left(v-u\right)^2}{2T}-\frac{1}{2}\right)\frac{1}{\rho}\mathcal{M} \, \right\},
\end{equation}
which is orthogonal with respect to the weighted inner product:
\begin{equation}
 \left( \psi_1, \, \psi_2 \right)_{\maxw^{-1}} := \int_{-\infty}^{\infty} \psi_1 \, \psi_2 \, \maxw^{-1} \, dv.
\end{equation} 
An orthogonal projection onto this basis for some function, $F(v):\reals \mapsto \reals$, such that
$F \maxw^{-1/2} \in L^2\left(\reals\right)$, is given by:
\begin{equation}
\label{eq5-4}
 \Pi_{\maxw}\left[ F \right]  = \left\{ A_1 + \left(\frac{v - u}{\sqrt{T}}\right) A_2 + \sqrt{2} \left( \frac{\left({v} - {u}\right)^2}{2T} - \frac{1}{2}\right) 
 A_3
 \right\} \maxw[f],
\end{equation}
where
\begin{align}
\label{eqn5-5}
\Bigl( A_1, \, A_2, \, A_3 \Bigr) :=
 \frac{1}{\rho} \int_{-\infty}^{\infty} \left( 1, \, \left( \frac{v-u}{\sqrt{T}} \right), \, 
 \sqrt{2} \left( \frac{(v-u)^2}{2T} - \frac{1}{2} \right) \right) \, F\left( v \right) \, dv.
\end{align}
This projection operator has the property that (see Bennoune et al. \cite{Bennoune2008UniformlyAsymptotics} for details):
\begin{equation}
\label{eq5-8}
  \Pi_{\maxw} \left[ \maxw \right] = \maxw, \quad
  \Pi_{\maxw} \left[ \maxw_{,t} \right] = \maxw_{,t}, \quad
  \Pi_{\maxw} \left[ g_{,t} \right] = 0, \quad
  \Pi_{\maxw} \left[ g \right] = 0,
\end{equation}
and when applied to \eqref{eq5-2} yields:
\begin{equation}
\label{eqn:projected_micmac_1D}
\maxw_{,t} + \Pi_{\maxw} \left[ {v}  \maxw_{,x} \right]  + \varepsilon \, 
\Pi_{\maxw} \left[ v g_{,x} \right] =  0.
\end{equation}
Subtracting \eqref{eqn:projected_micmac_1D} from \eqref{eq5-2}, dividing by $\varepsilon$, and rearranging terms yields an evolution equation for the microscopic portion $g$:
\begin{gather}
\label{eq5-9}
   g_{,t}  +
\left( {\mathcal I} - 
 \Pi_{\maxw} \right) \left[{v}  g_{,x} \right] = -\frac{1}{\varepsilon} \biggl( \tau \, g 
+ \left( {\mathcal I} - \Pi_{\maxw} \right) \left[{v}  \maxw_{,x} \right]  \biggr),
\end{gather}
where ${\mathcal I}$ is the identity operator.

Equation \eqref{eq5-9} can be further simplified, since it turns out that it is possible to analytically compute the orthogonal projection of the Maxwellian transport term (again, see Bennoune et al. \cite{Bennoune2008UniformlyAsymptotics} for details):
\begin{equation}
\label{eq5-15}
 \left( {\mathcal I} - \Pi_{\maxw} \right) \left[ v \maxw_{,x} \right]  = 
 \left\{ \left(\frac{(v - u)^3}{2 T} - \frac{3}{2} \left(v - u\right) \right) \, \frac{T_{,x}}{T}
 \right\} \maxw.
\end{equation}
Equation \eqref{eq5-9} can now be rewritten as follows:
\begin{equation}
\label{eq5-16}
   \left(\text{1D Micro}\right): \quad g_{,t}  +\left( {\mathcal I} - 
 \Pi_{\maxw} \right) \left[{v}  g_{,x} \right] = -\frac{\tau}{\varepsilon} \biggl(   g 
-  \widehat{g} \biggr), \qquad 
\widehat{g} := -\frac{1}{\tau} 
 \left(\frac{(v - u)^3}{2 T} - \frac{3}{2} \left(v - u\right) \right) \, \frac{T_{,x}}{{T}}\maxw.
\end{equation}
The final 1D micro-macro evolution equations are then given by \eqref{eq5-11} and \eqref{eq5-16}.

\begin{remark}
\label{note:note-NSF-limit}
For $\varepsilon \ll 1$, $g$ will be driven towards $\hat{g}$, and the full distribution function, $f$, will have the following form:
\begin{equation}
 f = \maxw[f] + \varepsilon \, \hat{g} + {\mathcal O}\left(\varepsilon^2\right).
\end{equation}
Therefore, the Maxwell-Boltzmann distribution, $\maxw[f]$, is the compressible Euler limit, while 
$\varepsilon \, \hat{g}$ is the Navier-Stokes-Fourier modification, and all other terms
represent higher-order corrections beyond Navier-Stokes-Fourier.
\end{remark}

\subsection{Setup for 1D Micro-Macro Finite Volume Method}
Now that the micro-macro formulation for $d_v=1$ has been established via \eqref{eq5-11} and \eqref{eq5-16}, we turn our attention to the numerical discretization. Although our version is similar to the approach of Bennoune, Lemou, and Mieussens \citep{Bennoune2008UniformlyAsymptotics}, we introduce a few small changes, including a version that does not use a staggered grid for $g$, which simplifies the scheme, and ultimately allows us to extend the method to higher dimensions (see Section \ref{sec:micmac2D}).

We introduce a uniform, cell-centered, finite volume mesh in the domain $\left(x_{\text{min}}, x_{\text{max}} \right) \times \left(v_{\text{min}}, v_{\text{max}} \right)$:
\begin{equation}
   x_i := x_{\text{min}} + \left( i - \frac{1}{2} \right){\Delta x} \quad \text{for} \quad i\in\left[1,N_x\right], \quad
   v_k := v_{\text{min}} + \left( k - \frac{1}{2} \right){\Delta v} \quad \text{for} \quad k\in\left[1,N_v\right],
\end{equation}
where
\begin{equation}
{\Delta x} := \frac{x_{\text{max}}-x_{\text{min}}}{N_x} \qquad \text{and} \qquad
{\Delta v} := \frac{v_{\text{max}}-v_{\text{min}}}{N_v}.
\end{equation}
At each time level, $t^n = n \Delta t$, the macroscopic and microscopic portions of the distribution functions are given by
\begin{equation}
\vec{Q^n_i} := \left( \, \rho^n_i, \, \rho u^n_i, \, {\mathcal E}^n_i \, \right) \quad \text{for} \quad i\in\left[1,N_x\right] 
\quad \text{and} \quad
     G^n_{ik} \quad \text{for} \quad \left(i,k \right) \in\left[1,N_x\right] \times \left[1,N_v\right],
\end{equation}
where
\begin{equation}
\vec{Q^n_i} \approx \frac{1}{2} \int_{-1}^{1} \vec{q}\left(t^n, \, x_i + \frac{\Delta x \, \xi}{2} \right) \, d\xi,
\quad
G^n_{ik} \approx \frac{1}{4} \iint_{-1}^{1} g\left(t^n, \, x_i + \frac{\Delta x \, \xi}{2}, \,
   v_k + \frac{\Delta v \, \mu}{2} \right) \, d\xi \, d\mu.
\end{equation}
The solution is initialized at time $t=0$ and is run to $t=T_{\text{final}}$. The time step is chosen as a global constant and is defined as follows:
\begin{equation}
  \text{Initialize}: \quad \Delta t = \frac{\text{CFL} \cdot \Delta x}{\mathcal{V}}, \quad
  N_{\text{steps}} = \left\lceil \frac{T_{\text{final}}}{\Delta t} \right\rceil \quad \Longrightarrow \quad
  \text{Reset}: \quad \Delta t \leftarrow \frac{T_{\text{final}}}{N_{\text{steps}}}, \quad 
  \text{CFL} \leftarrow  \frac{\mathcal{V} \, \Delta t}{\Delta x},
\end{equation}
where $N_{\text{steps}}$ is the total number of timesteps needed to advance from $t=0$ to $t=T_{\text{final}}$ and
\begin{equation}
\left(\text{CFL number}\right): \quad 0 < \text{CFL} < 1, \qquad
\left(\text{max speed}\right): \quad \mathcal{V} := \max\left\{ \left| v_{\text{min}} \right|, \left|v_{\text{max}}\right| \right\}.
\end{equation}

\subsection{1D Micro Finite Volume Method}
\label{subsec:micro-fvm-1d}
In this section, we demonstrate how to update the current value of the micro portion, $G^n$, to its value at the next time step, $G^{n+1}$:
\begin{equation}
   \text{Input:} \quad \vec{Q^n_i} = \left( \, \rho^n_i, \, \rho u^n_i, \, {\mathcal E}^n_i \, \right) \quad \text{and} \quad
     G^n_{ik} \quad \Longrightarrow \quad \text{Output:} \quad G^{n+1}_{ik}.
\end{equation}
First, we construct an approximate $\hat{g}$ from \eqref{eq5-16} on the numerical grid:
\begin{align}
\label{eqn5-19}
\widehat{G}_{ik}^n &:= -\frac{1}{\tau^n_i} 
 \left(\frac{(v_k - u^n_i)^3}{2 T^n_i} - \frac{3}{2} \left(v_k - u^n_i\right) \right) \, \left( \frac{T^n_{i+1/2}-T^n_{i-1/2}}{\Delta x \, {T^n_i}} \right) \, 
  \maxw^n_{ik}, \quad T^n_{i+1/2} := \frac{1}{2} \left( T^n_{i+1} + T^n_{i} \right), \\ 
  \maxw^n_{ik} &:= \frac{\rho^n_i}{\sqrt{2\pi T^n_i}} \exp\left[
-\frac{\left(v_k - u^n_i \right)^2}{2 T^n_i} 
\right].
\end{align}
We approximate the term $v g_{,x}$ in \eqref{eq5-16} via an upwind gradient:
\begin{align}
  \label{eqn5-20}
  Z^n_{ik} &:= v_k^{-} \left( \frac{G^n_{i+1 k} - G^n_{ik}}{\Delta x} \right)
 + v_k^{+} \left( \frac{G^n_{i k} - G^n_{i-1 k}}{\Delta x}\right), \quad
 v_k^{-} := \min\left(0,v_k\right), \quad
 v_k^{+} := \max\left(0,v_k\right),
\end{align}
and then compute its orthogonal projection onto the nullspace of the collision operator (see \eqref{eq5-4}--\eqref{eqn5-5}):
\begin{align}
\label{eqn5-24}
 \widehat{Z}^n_{ik} &:= 
 \left\{ \left[A_1\right]^n_{i} + \left(\frac{v_k - u^n_i}{\sqrt{T^n_i}}\right) \left[A_2\right]^n_{i} 
    + \sqrt{2} \left( \frac{\left({v_k} - {u^n_i}\right)^2}{2 T^n_i} - \frac{1}{2}\right) 
 \left[A_3\right]^n_{i}
 \right\} \maxw^n_{ik},
\end{align}
where the coefficients are computed via numerical quadrature using the midpoint rule:
\begin{equation}
\int_{-\infty}^{\infty} g(v) \, dv \approx \int_{v_{\text{min}}}^{v_{\text{max}}} g(v) \, dv =
			\sum_{k=1}^{N_v} \int_{v_k - \frac{\Delta v}{2}}^{v_k + \frac{\Delta v}{2}} g(v) \, dv
			\approx 
			\Delta v \sum_{k=1}^{N_v} g\left(v_k\right)
\end{equation}
so that
\begin{align}
 \label{eqn5-21}
\biggl\{ \left[A_1\right]^n_{i}, \quad 
\left[A_2\right]^n_{i}, \quad
\left[A_3\right]^n_{i} \biggr\} := \frac{\Delta v}{\rho^n_i} \sum_{k=1}^{N_v} \left\{ Z^n_{ik}, \quad
\left( \frac{v_k-u^n_i}{\sqrt{T^n_i}} \right) \, Z^n_{ik}, \quad 
\sqrt{2} \left( \frac{(v_k-u^n_i)^2}{2 T^n_i} - \frac{1}{2} \right) Z^n_{ik} \right\}.
\end{align}
Lastly, the numerical micro solution, $G^n$, is updated to its new value, $G^{n+1}$, by applying to the micro equation \eqref{eq5-16} a mixed forward Euler (for the transport term) and backward Euler (for the collision term) time discretization:
\begin{align}
\label{eqn5-25}
 G^{n+1}_{ik} &= \left( \frac{\varepsilon}{\varepsilon + \Delta t \, \tau^n_i} \right) \biggl\{ G^n_{ik} -
\Delta t \left( Z^n_{ik} - 
 \widehat{Z}^n_{ik} \right) \biggr\} 
  +\left( \frac{\Delta t \, \tau^n_i}{\varepsilon + \Delta t \, \tau^n_i} \right) \widehat{G}^{n}_{ik}. 
\end{align}
Note that this update is {\it asymptotic-preserving} \cite{Jin2012AsymptoticReview}, since we can take the limit $\varepsilon \rightarrow 0^+$ of this expression for fixed mesh parameters, $\Delta t$, $\Delta x$, and $\Delta v$, and arrive at the following result:
 \begin{align}
  G^{n+1}_{ik} \rightarrow \widehat{G}^{n}_{ik} \quad \text{as} \quad \varepsilon \rightarrow 0^+.
\end{align}

\subsubsection{Extrapolation Boundary Conditions in the Micro Step}
We briefly mention here how simple boundary conditions are implemented in the 1D micro update. For extrapolation boundary conditions, we do the following:
\begin{align}
\label{eqn:1d_extrap_micro_1}
 \text{Left BC:} & \quad T^n_{1/2} = T^n_{1}, \quad Z^n_{\left(1, \, k\right)} = v_k^{-} \left( \frac{G^n_{\left(2, \, k\right)} - G^n_{\left(1, \, k\right)}}{\Delta x} \right), \\
\label{eqn:1d_extrap_micro_2}
 \text{Right BC:} & \quad T^n_{N_x + 1/2} = T^n_{N_x}, \quad Z^n_{\left(N_x, \, k\right)} = 
 v_k^{+} \left( \frac{G^n_{\left(N_x, \, k\right)} - G^n_{\left(N_x-1, \, k\right)}}{\Delta x}\right).
\end{align}

\subsubsection{Periodic Boundary Conditions in the Micro Step}
For periodic boundary conditions, we do the following:
\begin{align}
\label{eqn:1d_periodic_micro_1}
 \text{Left BC:} & \quad T^n_{1/2} = \frac{T^n_{1} + T^n_{N_x}}{2}, \quad Z^n_{\left(1, \, k\right)} = v_k^{-} \left( \frac{G^n_{\left(2, \, k\right)} - G^n_{\left(1, \, k\right)}}{\Delta x} \right)
 + v_k^{+} \left( \frac{G^n_{\left(1, \, k\right)} - G^n_{\left(N_x, \, k\right)}}{\Delta x}\right), \\
\label{eqn:1d_periodic_micro_2}
 \text{Right BC:} & \quad  T^n_{N_x + 1/2} = \frac{T^n_{1} + T^n_{N_x}}{2}, \quad Z^n_{\left(N_x, \, k\right)} = v_k^{-} \left( \frac{G^n_{\left(1, \, k\right)} - G^n_{\left(N_x, \, k \right)}}{\Delta x} \right)
 + v_k^{+} \left( \frac{G^n_{\left(N_x, \, k\right)} - G^n_{\left(N_x-1, \, k\right)}}{\Delta x}\right).
\end{align}


\subsection{1D Macro Finite Volume Method}
In this section, we demonstrate how to update the current value of the macro variables at $t=t^n$ to the new values at $t=t^{n+1} = t^n + \Delta t$:
\begin{equation}
   \text{Input:} \quad \vec{Q^n_i} = \left( \, \rho^n_i, \, \rho u^n_i, \, {\mathcal E}^n_i \, \right) \quad \text{and} \quad
     H^{n+1}_{i} \quad \Longrightarrow \quad \text{Output:} \quad \vec{Q^{n+1}_i} = \left( \, \rho^{n+1}_i, \, \rho u^{n+1}_i, \, {\mathcal E}^{n+1}_i \, \right),
\end{equation}
where the micro portion has already been updated to $t=t^{n+1}$ via \eqref{eqn5-25}, and the heat flux is computed via
\begin{align}
\label{eqn5-26}
H^{n+1}_i &:=   \frac{\varepsilon}{2} \, \Delta v \sum_{k=1}^{N_v} \left(v_k\right)^3 \, G^{n+1}_{ik} \, .
\end{align}
The macro equation given by \eqref{eq5-11} is discretized using the following finite volume method:
\begin{align}
\label{eqn5-27}
\vec{Q^{n+1}_i} &= \vec{Q^n_i} - \frac{\Delta t}{\Delta x} \left( \, \vec{\mathcal{F}^{ \, n}_{i+1/2}} - \vec{\mathcal{F}^{ \, n}_{i-1/2}} \, \right)
- \frac{\Delta t}{\Delta x} \left( \begin{bmatrix} 0 \\ 0 \\ H^{n+1}_{i+1/2} \end{bmatrix} - 
\begin{bmatrix} 0 \\ 0 \\ H^{n+1}_{i-1/2} \end{bmatrix} \right), \quad
H^{n+1}_{i+1/2} := \frac{H^{n+1}_{i+1} + H^{n+1}_{i}}{2},
\end{align}
where the numerical fluxes, $\mathcal{F}$, are computed via kinetic flux vector splitting (KFVS) \citep{Mandal1994KineticEquations}. We note that other types of numerical fluxes could also be used; we select KFVS for two key reasons: (1) it is less numerically diffusive than many other fluxes (e.g., Lax-Friedrichs, Rusanov, and HLL); and (2) it allows for direct implementation of complicated boundary conditions (i.e., diffusely reflecting walls). In 1D, KFVS produces the following numerical flux at cell interface $x_{i-1/2}$:
\begin{equation}
\label{eqn:flux1d_1}
\begin{split}
  \vec{\mathcal{F}^{ \, n}_{i-1/2}} &=  \int_0^\infty v \, \vec{m}\left(v\right) \, \maxw^n_{i-1} \, dv + 
     \int_{-\infty}^0 v \, \vec{m}\left(v \right) \, \maxw^n_i \, dv \\
  &= \alpha^{n}_{i-1}  \, 
  \begin{pmatrix}
  \rho^n_{i-1} \\ \rho^n_{i-1} u^n_{i-1} \\ \frac{1}{2} \rho^n_{i-1} \left( 2 T^n_{i-1} + \left(u^n_{i-1}\right)^2 \right)
  \end{pmatrix} + 
  \beta^{n+}_{i-1}  \, \vec{f}\left( \vec{Q^n_{i-1}} \right)
  -\alpha^n_i \, 
  \begin{pmatrix}
  \rho^n_i \\ \rho^n_i u^n_i \\ \frac{1}{2} \rho^n_i \left( 2 T^n_i + \left(u^n_i\right)^2 \right)
  \end{pmatrix} + 
  \beta^{n-}_i  \, \vec{f}\left( \vec{Q^n_i} \right), 
\end{split}
\end{equation}
where
\begin{align}
\label{eqn:flux1d_2}
    \vec{f}\left( \vec{Q^n_i} \right) = \begin{pmatrix} \rho^n_i u^n_i \\ \rho^n_i \left( T^n_i + \left( u^n_i\right) ^2\right) \\ 
    \frac{1}{2}\rho^n_i u^n_i \left( 3 T^n_i + \left(u^n_i\right)^2\right) \end{pmatrix}, \quad 
    \alpha^n_{i} =  \sqrt{\frac{T^n_i}{2\pi}}\exp\left(-\frac{\left(u^n_i\right)^2}{2T^n_i}\right), \quad 
    \beta_{i}^{n\pm} = \frac{1}{2} \left(1 \pm \mathrm{erf}\left(\frac{u^n_i}{\sqrt{2T^n_i}}\right) \right).
\end{align}

\subsubsection{Extrapolation Boundary Conditions in the Macro Step}
We briefly mention how simple boundary conditions are implemented in the 1D macro update. For extrapolation boundary conditions, we do the following:
\begin{align}
\label{eqn:1d_extrap_macro_1}
 \text{Left BC:} & \quad H^n_{1/2} = H^n_{1}, \quad \vec{\mathcal{F}^{ \, n}_{1/2}} =  \int_{-\infty}^\infty v \, \vec{m}\left(v\right) \, \maxw^n_{1} \, dv = \vec{f}\left( \vec{Q^n_1} \right), \\
 \label{eqn:1d_extrap_macro_2}
 \text{Right BC:} & \quad H^n_{N_x + 1/2} = H^n_{N_x}, \quad \vec{\mathcal{F}^{ \, n}_{N_x + 1/2}} =  \int_{-\infty}^\infty v \, \vec{m}\left(v\right) \, \maxw^n_{N_x} \, dv = \vec{f}\left( \vec{Q^n_{N_x}} \right).
\end{align}

\subsubsection{Periodic Boundary Conditions in the Macro Step}
For periodic boundary conditions, we do the following:
\begin{align}
\label{eqn:1d_periodic_macro}
 H^n_{1/2} = H^n_{N_x + 1/2} = \frac{H^n_{1} + H^n_{N_x}}{2}, \quad 
 \vec{\mathcal{F}^{ \, n}_{1/2}} = \vec{\mathcal{F}^{ \, n}_{N_x+1/2}} = 
 \int_0^\infty v \, \vec{m}\left(v\right) \, \maxw^n_{N_x} \, dv + 
     \int_{-\infty}^0 v \, \vec{m}\left(v \right) \, \maxw^n_1 \, dv.
\end{align}


\section{1D Numerical Examples}
\label{sec:numerical1D}
In this section, we verify the correctness and accuracy of the 1D micro-macro method described above
on three test cases. The first example demonstrates the convergence of the proposed method using a manufactured solution. The second example shows the ability of the scheme to compute Riemann solutions. The final example illustrates how the scheme can compute boundary-driven flows for Knudsen numbers in all flow regimes.

\begin{table}[!th]
\begin{Large}
\begin{tabular}{|c|c|c|c|c|c|} \hline
{\normalsize {$N_x$}} & {\normalsize $N_{v}$} & {\normalsize Macro Error} & {\normalsize $\log_2\left(\text{Error Ratio}\right)$} & {\normalsize Micro Error} & {\normalsize $\log_2\left(\text{Error Ratio}\right)$} \\ \hline \hline
{\normalsize 10}  & {\normalsize 10}  & {\normalsize $4.841399 \times 10^{-2}$} & {\normalsize --} & {\normalsize $9.339605 \times 10^{-2}$} & {\normalsize --} \\ \hline{\normalsize 20}  & {\normalsize 20}  & {\normalsize $2.550892 \times 10^{-2}$} & {\normalsize $0.924422$} & {\normalsize $2.801581 \times 10^{-2}$} & {\normalsize $1.737120$} \\ \hline{\normalsize 40}  & {\normalsize 40}  & {\normalsize $1.334602 \times 10^{-2}$} & {\normalsize $0.934592$} & {\normalsize $1.505672 \times 10^{-2}$} & {\normalsize $0.895833$} \\ \hline{\normalsize 80}  & {\normalsize 80}  & {\normalsize $6.844527 \times 10^{-3}$} & {\normalsize $0.963387$} & {\normalsize $7.832841 \times 10^{-3}$} & {\normalsize $0.942800$} \\ \hline{\normalsize 160}  & {\normalsize 160}  & {\normalsize $3.468271 \times 10^{-3}$} & {\normalsize $0.980734$} & {\normalsize $4.002026 \times 10^{-3}$} & {\normalsize $0.968805$} \\ \hline{\normalsize 320}  & {\normalsize 320}  & {\normalsize $1.745691 \times 10^{-3}$} & {\normalsize $0.990418$} & {\normalsize $2.022980 \times 10^{-3}$} & {\normalsize $0.984248$} \\ \hline{\normalsize 640}  & {\normalsize 640}  & {\normalsize $8.752563 \times 10^{-4}$} & {\normalsize $0.996021$} & {\normalsize $1.015872 \times 10^{-3}$} & {\normalsize $0.993764$} \\ \hline\end{tabular}
\caption{(\S\ref{subsec:method_of_manufacutured_solutions}: 1D1V method of manufactured solutions example) Error table showing first order convergence in both the macroscopic and microscopic variables with $\varepsilon = 0.1$, $(x,v) \in \left[0,1\right] \times \left[-6.5, 6.5 \right]$ with periodic boundary conditions in $x$, $t=0.9351$, $\text{CFL}=0.95$, and $\tau$ is taken as \eqref{eqn:collision_parameter_1d}. \label{tab5-2}}
\end{Large}
\end{table}

\subsection{Method of Manufactured Solutions Example}
\label{subsec:method_of_manufacutured_solutions}
One challenge with the Boltzmann-BGK system is that there is no simple way to generate
exact solutions of \eqref{eq5-8} other than in some simple cases. 
For this reason, we consider here the method of manufactured solutions with the following invented ``solution''
on $0\le x \le 1$ with periodic boundary conditions (see \eqref{eqn:1d_periodic_micro_1}, \eqref{eqn:1d_periodic_micro_2}, and \eqref{eqn:1d_periodic_macro}):
\begin{align}
    f(t,x,v) = \left(e^{-(v-1)^2}+2e^{-(v+1)^2}\right)\left(2 + \sin\left(2\pi\left(x-t\right)\right)\right),
\end{align}
from which we calculate the quantities:
\begin{gather}
\left(\rho, \, u, \, T \right)(t,x) = \Biggl( 3\sqrt{\pi} \Bigl(2 + \sin\bigl(2\pi\left(x-t\right)\bigr)\Bigr), \,
	-\frac{1}{3}, \, \frac{25}{18} \Biggr), \\
	g(t,x,v) =  \frac{1}{\varepsilon} \left( e^{-(v-1)^2} + 2 e^{-(v+1)^2}
	   - 1.8 e^{- 0.04 \left( {3 v + 1} \right)^2} \right) \Bigl(2 + \sin\bigl(2\pi\left(x-t\right)\bigr)\Bigr).
\end{gather}
Importantly, this distribution has a spatially and temporally varying heat flux. Since this ``solution'' does not actually solve the Boltzmann-BGK equation, we need to compute the residual,
\begin{align}
S(t,x,v) = f_{,t} + v f_{,x} - \dfrac{\tau}{\varepsilon}(\mathcal{M} - f),
\end{align}
and include it as a source term in the micro-macro method. In particular, the 1D micro-macro scheme described above can be run as normal with two small exceptions. First, in the micro update, we modify $\widehat{G}$ as defined in \eqref{eqn5-19} as follows: 
\begin{align}
\widehat{G}_{ik}^n &= -\frac{1}{\tau^n_i} 
 \left(\frac{(v_k - u^n_i)^3}{2 T^n_i} - \frac{3}{2} \left(v_k - u^n_i\right) \right) \, \left( \frac{T^n_{i+1/2}-T^n_{i-1/2}}{\Delta x \, {T^n_i}} \right) \, 
  \maxw^n_{ik}  + \frac{1}{\tau^n_i} \left( {\mathcal I} - \Pi_{\maxw} \right) \left[S\right]\left(t^n, x_i, v_k\right),
\end{align}
Second, in the macro update we modify \eqref{eqn5-27} by adding a source term:
\begin{align}
\vec{Q^{n+1}_i} &= \vec{Q^n_i} - \frac{\Delta t}{\Delta x} \left( \, \vec{\mathcal{F}^{ \, n}_{i+1/2}} - \vec{\mathcal{F}^{ \, n}_{i-1/2}} \, \right)
- \frac{\Delta t}{\Delta x} \left( \begin{bmatrix} 0 \\ 0 \\ H^{n+1}_{i+1/2} \end{bmatrix} - 
\begin{bmatrix} 0 \\ 0 \\ H^{n+1}_{i-1/2} \end{bmatrix} \right)
+ \Delta t \, \pi^{\frac{3}{2}} \, \cos\left(2\pi\left(x_i-t^n\right)\right) \, \begin{bmatrix} -8 \\ 11 \\ -7 \end{bmatrix}.
\end{align}

A convergence study for this problem is shown in Table \ref{tab5-2}. In these calculations we take $(x,v) \in \left[0,1\right] \times \left[-6.5, 6.5 \right]$ with periodic boundary conditions in $x$, the solution is computed out to $t=0.9351$, the time step is chosen so that $\text{CFL}=0.95$, the Knudsen number is $\varepsilon = 0.1$. The collision parameter is set according to the hard-sphere collision value described in Bird \cite{Bird1994MolecularEdition}:
\begin{equation}
\label{eqn:collision_parameter_1d}
\tau^n_i = \frac{16}{5} \sqrt{\frac{T^n_i}{2 \pi}}.
\end{equation} 
The numerical errors shown are the relative $L_2$ errors of the macroscopic, $\vec{Q} = \left(\rho, \, \rho u, \,{\mathcal E} \right)$, and microscopic, $g$, variables: 
\begin{equation}
\begin{matrix}
\text{macro} \\
\text{error}
\end{matrix} := \sqrt{ {{\displaystyle \sum_{i=1}^{N_x}} \left\| \vec{Q_i} - \vec{Q^{\text{exact}}_i} \right\|^2}
\Bigg/ {{\displaystyle \sum_{i=1}^{N_x}} \left\| \vec{Q^{\text{exact}}_i} \right\|^2}}, \quad
\begin{matrix}
\text{micro} \\
\text{error}
\end{matrix} := \sqrt{ {{\displaystyle \sum_{i=1}^{N_x}} \sum_{k=1}^{N_{v}}} \left| {G_{ik}} - {g^{\text{exact}}_{ik}} \right|^2
\Bigg/ {\displaystyle \sum_{i=1}^{N_x} \sum_{k=1}^{N_{v}}} \left| {g^{\text{exact}}_{ik}} \right|^2}.
\end{equation}
The results in Table \ref{tab5-2} clearly show first-order convergence in all variables. This table can be generated using the Python/C++ companion code \cite{code:RossmanithSar2025} by going to the code directory and typing:
\begin{tcolorbox}
\begin{verbatim}
cd 1d/example1-manufactured-solution/; python make_paper_figures.py
\end{verbatim}
\end{tcolorbox}

\begin{figure}[!th]
\begin{center}
    \begin{tabular}{cc}
    (a)\includegraphics[width=0.44\textwidth]{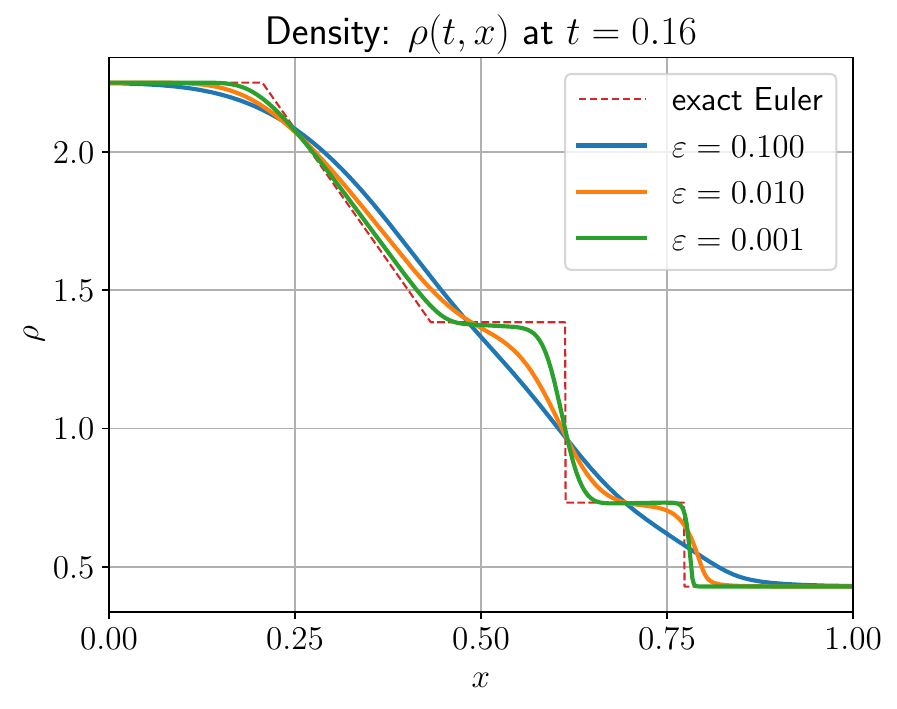} & 
    (b)\includegraphics[width=0.44\textwidth]{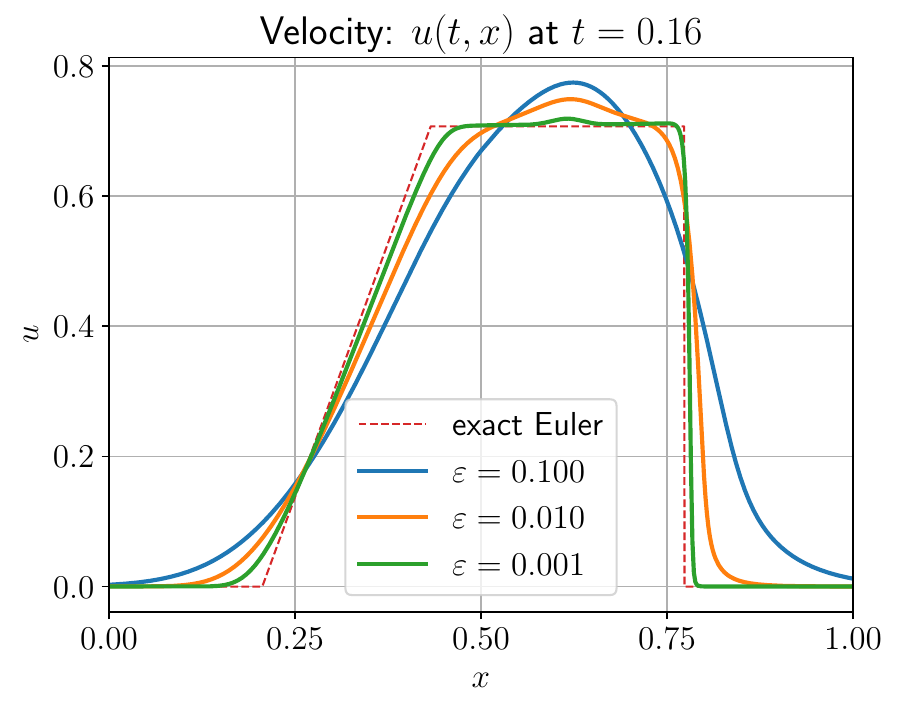} \\
    (c)\includegraphics[width=0.44\textwidth]{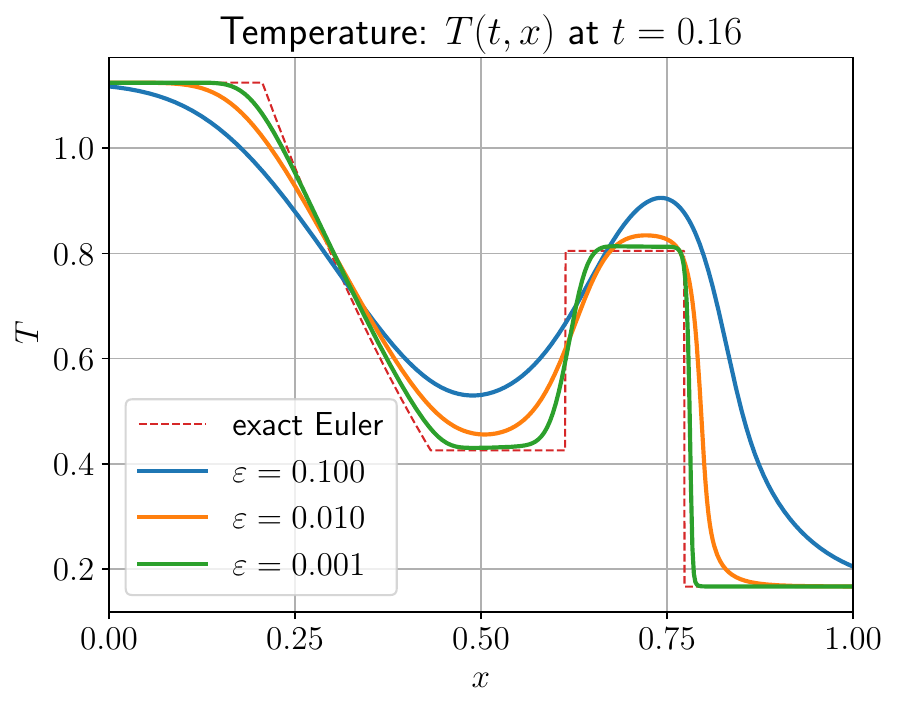} &
    (d)\includegraphics[width=0.44\textwidth]{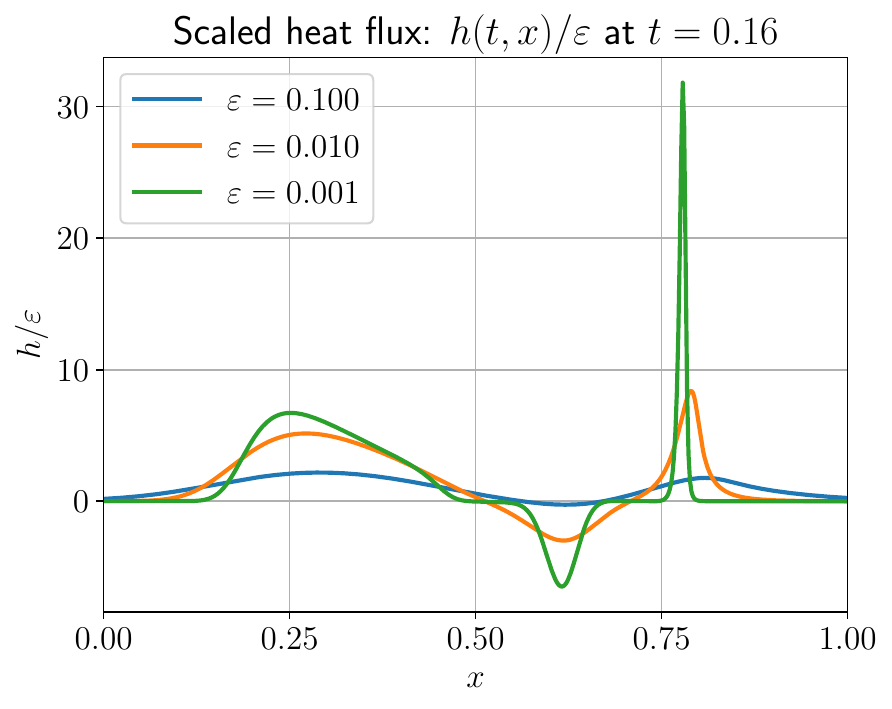}
    \end{tabular}
    \caption{(\S\ref{subsec:shock_tube}: shock tube example) Shock tube problem using the 1D micro-macro solver with three different Knudsen numbers: $\varepsilon = 0.1, 0.01, 0.001$, as well as the exact Euler solution ($\varepsilon \rightarrow 0^+$). 
  The full domain is $(x,v) \in \left[-0.25,1.25\right] \times \left[-4.5, 4.5 \right]$ with constant extrapolation boundary conditions in $x$, the grid resolution is $N_x \times N_v=768 \times 128$,
   the CFL number is $\text{CFL}=0.991$, and $\tau$ is \eqref{eqn:collision_parameter_1d}.
    Shown in the panels are the following quantities for $x \in \left[0,1\right]$: (a) density: $\rho$, (b) macroscopic velocity: $u$,
    (c) scalar temperature: $T$, and (d) heat flux divided by $\varepsilon$: $h/\varepsilon$.}
    \label{fig5-1}
\end{center}
\end{figure}

\subsection{Shock-Tube Problem}
\label{subsec:shock_tube}
Next, we apply the 1D micro-macro scheme to the Sod shock tube problem \cite{article:Sod1978}, which can be found in many places, including Bennoune et al. \cite{Bennoune2008UniformlyAsymptotics}. The initial conditions are
\begin{equation}
    (\rho, u, T)(t=0,x) = 
    \begin{cases}
        \bigl( 1.000, \, \, 0, \, \, 1.0 \bigr) & x < 0.5 \\
        \bigl( 0.125, \, \, 0, \, \, 0.8 \bigr) & x > 0.5 
    \end{cases} \qquad \text{and} \qquad
    g(t=0,x,v) = 0.
\end{equation}
In the calculations we take $(x,v) \in \left[-0.25,1.25\right] \times \left[-4.5, 4.5 \right]$ with constant extrapolation boundary conditions in $x$ (see \eqref{eqn:1d_extrap_micro_1}--\eqref{eqn:1d_extrap_micro_2} and \eqref{eqn:1d_extrap_macro_1}--\eqref{eqn:1d_extrap_macro_2}), the grid resolution is $N_x \times N_v=768 \times 128$, the solution is computed out to $t=0.16$ with 372 time steps ($\text{CFL} \approx 0.991$ and
$\Delta t \approx 4.301075\times10^{-4}$), and the collision parameter, $\tau$, is again \eqref{eqn:collision_parameter_1d}. The plots for the macroscopic quantities in the range $0 \le x \le 1$ are shown in Figure \ref{fig5-1} with three different Knudsen numbers: $\varepsilon = 0.1, 0.01, 0.001$. For the density ($\rho$), macroscopic velocity ($u$), and temperature ($T$), we also include the exact Riemann solution for the 1D compressible Euler equations (i.e., the $\varepsilon \rightarrow 0^+$ limit) (e.g., see \cite{article:Ketcheson2023}). These results are consistent with those in the literature (e.g., Bennoune et al. \cite{Bennoune2008UniformlyAsymptotics} and Pieraccini and Puppo \cite{Pieraccini2007Implicit-explicitEquations}). This figure can be generated using the Python/C++ companion code \cite{code:RossmanithSar2025} by going to the code directory and typing:
\begin{tcolorbox}
\begin{verbatim}
cd 1d/example2-shocktube/; python make_paper_figures.py
\end{verbatim}
\end{tcolorbox}

\begin{figure}[!th]
\begin{center}
    \begin{tabular}{cc}
    (a)\includegraphics[width=0.44\textwidth]{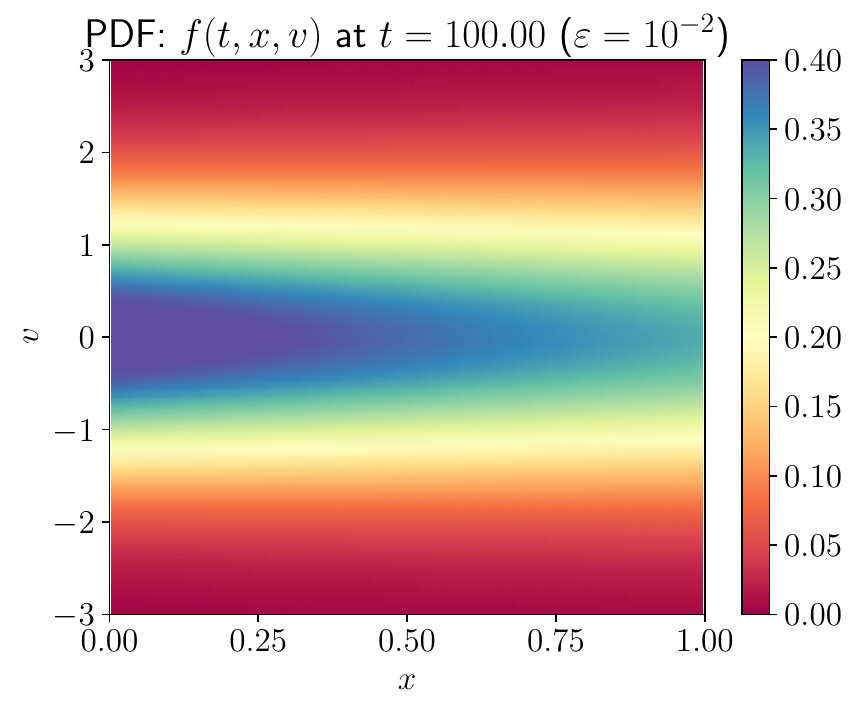} & 
    (b)\includegraphics[width=0.44\textwidth]{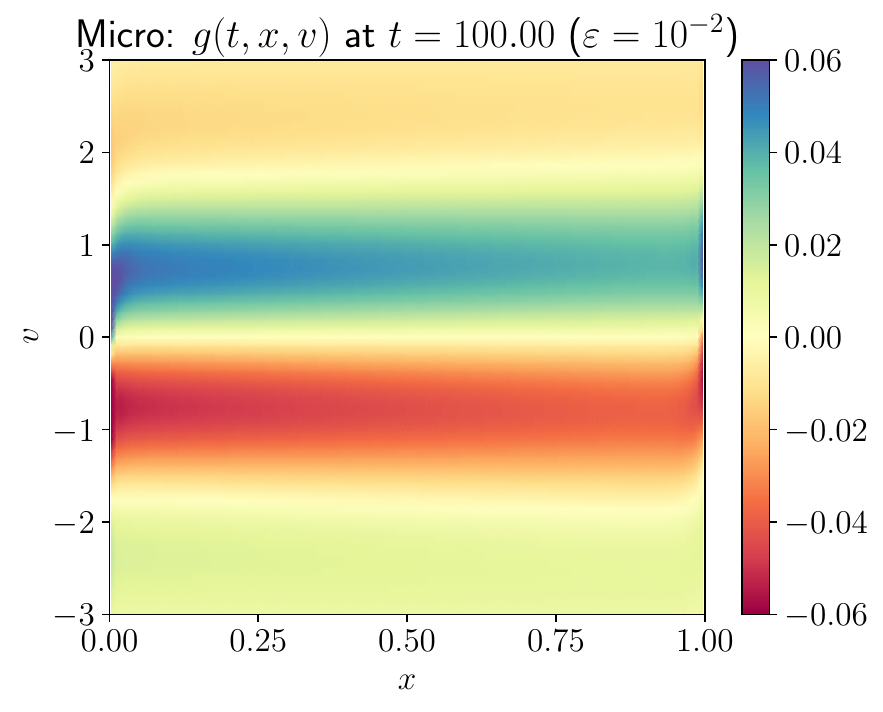} \\
    (c)\includegraphics[width=0.44\textwidth]{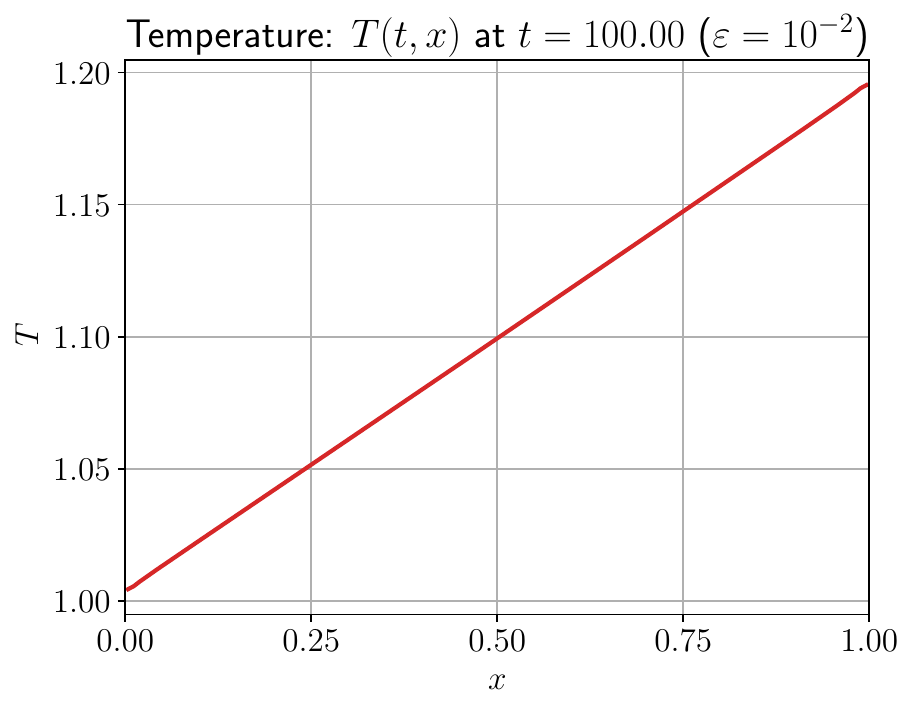} &
    (d)\includegraphics[width=0.44\textwidth]{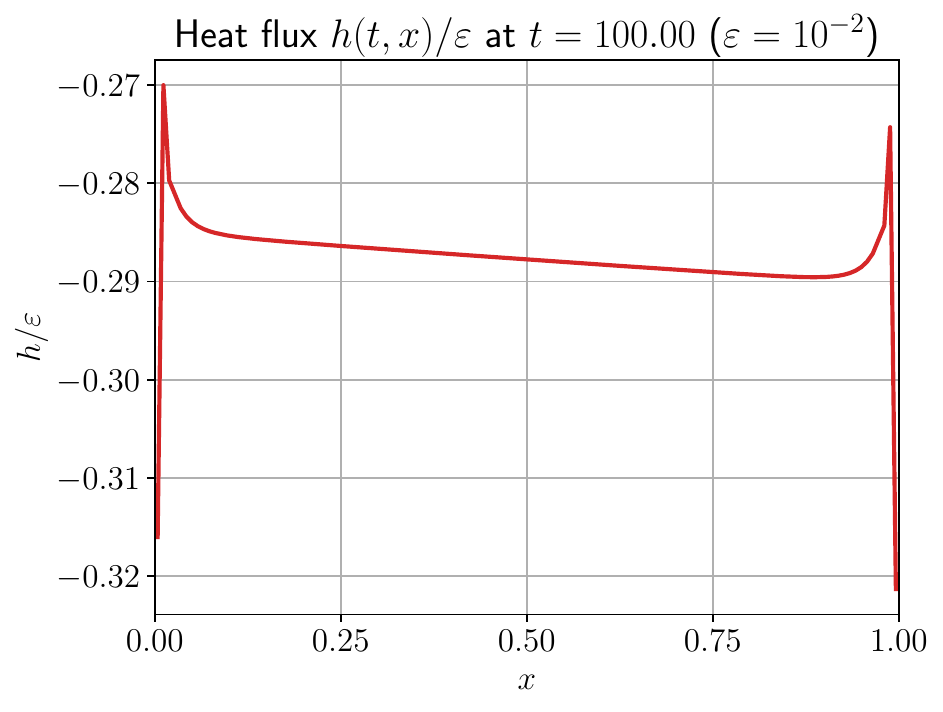}
    \end{tabular}
    \caption{(\S\ref{subsec:heat_transfer}: heat transfer problem) Micro-macro solution at time $t=100$ with Knudsen number $\varepsilon=10^{-2}$ on $(x,v) \in \left[0,1\right] \times \left[-6,6\right]$ with $N_x = N_v = 129$. Shown in the panels are the (a) full PDF (the plot only shows $v\in[-3,3]$): $f= \maxw+\varepsilon g$, (b) microscopic portion (the plot only shows $v\in[-3,3]$): $g$,
    (c) scalar temperature: $T$, and (d) heat flux divided by $\varepsilon$: $h/\varepsilon$.}
    \label{fig:1dtest3a}
\end{center}
\end{figure}

\begin{figure}[!th]
\begin{center}
    \begin{tabular}{cc}
    (a)\includegraphics[width=0.44\textwidth]{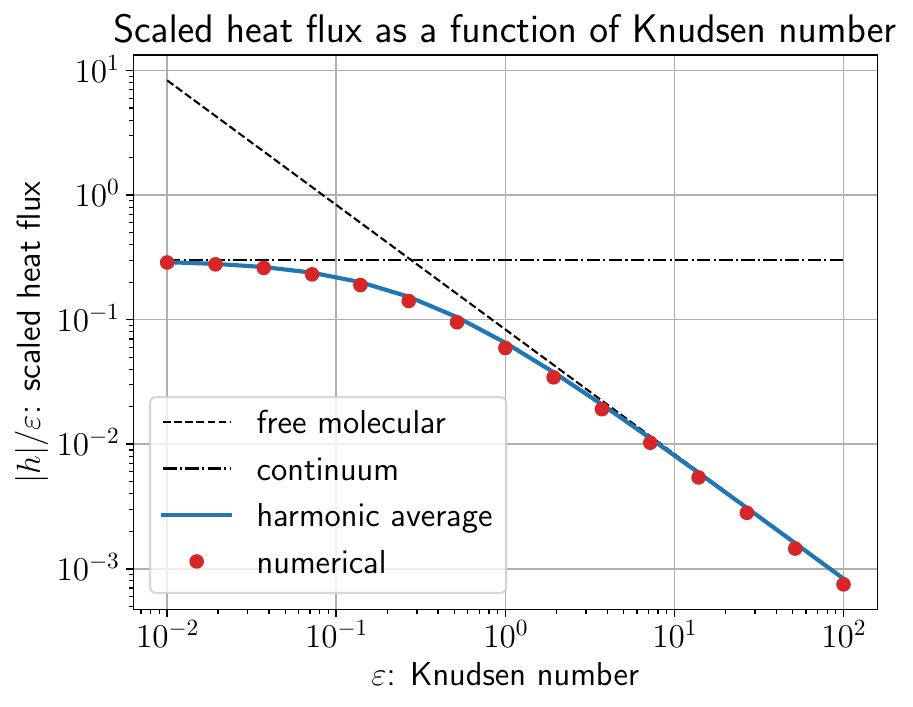} & 
    (b)\includegraphics[width=0.44\textwidth]{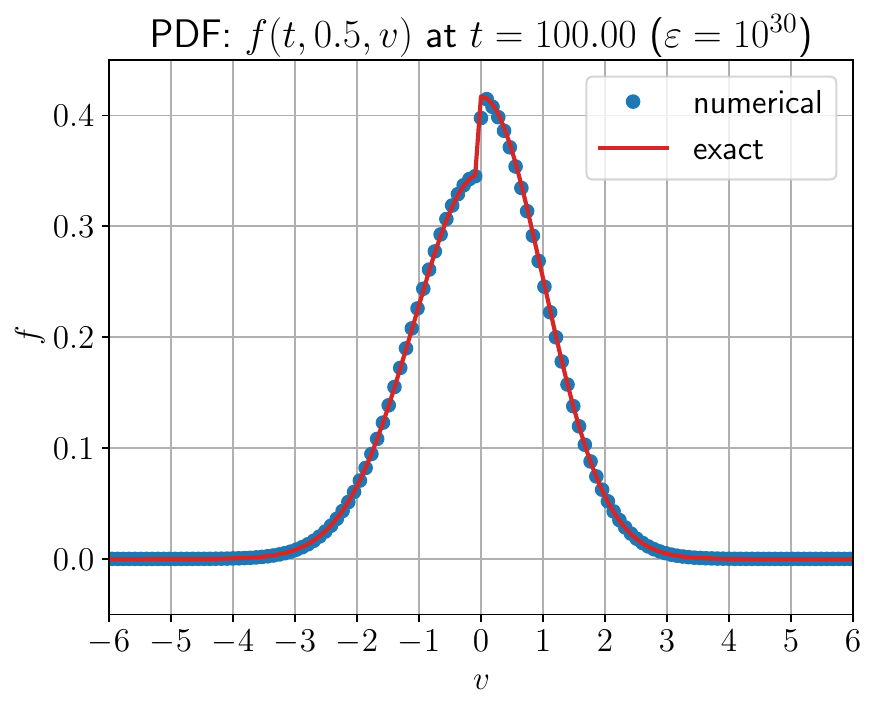}
    \end{tabular}
    \caption{(\S\ref{subsec:heat_transfer}: heat transfer problem) Micro-macro solution at time $t=100$  on $(x,v) \in \left[0,1\right] \times \left[-6,6\right]$ with $N_x = N_v = 129$. Shown in Panel (a) is the absolute value of the scaled heat flux $(|h|/\varepsilon)$ at $t=100$ and $x=0.5$ as a function of the Knudsen number, $\varepsilon$, alongside the values predicted by
    free molecular flow ($\varepsilon \rightarrow \infty$), continuum flow (i.e., Navier-Stokes-Fourier), and
    the harmonic average of free molecular and continuum flow. Shown in Panel (b) is a comparison of the exact
    solution for the full PDF with $\varepsilon \rightarrow \infty$ at $x=0.5$ and the one computed with the micro-macro
    scheme with $\varepsilon = 10^{30}$.}
    \label{fig:1dtest3b}
\end{center}
\end{figure}

\subsection{Heat Transfer Problem}
\label{subsec:heat_transfer}
The two previous examples had simple boundary conditions (periodic and constant extrapolation). To test the micro-macro approach on a problem with more complex boundary conditions, we consider the heat transfer problem involving a gas on the domain $0 \le x \le 1$ between two walls at fixed temperatures (diffusively reflecting). This problem has been considered by several papers in the literature (e.g., Gamba et al. \cite{article:Gamba2010}, Abdelmalik et al. \cite{Abdelmalik2017ErrorEquation}, and van der Woude et al. \cite{article:vanderWoude2024}). 
At the kinetic level, the boundary conditions are 
\begin{alignat}{2}
\label{eqn:wall_bc_1}
f(t,x=0,v) &= \maxw_C(t,v) = \frac{\rho_C(t)}{\sqrt{2\pi T_C}} \exp\left[-\frac{v^2}{2T_C}\right] \quad &v>0&, \\
\label{eqn:wall_bc_2}
f(t,x=1,v) &= \maxw_H(t,v) = \frac{\rho_H(t)}{\sqrt{2\pi T_H}} \exp\left[-\frac{v^2}{2T_H}\right] \quad &v<0&.
\end{alignat}
Following \cite{Abdelmalik2017ErrorEquation}, we take $T_C=1.0$ and $T_H=1.2$; the wall densities, $\rho_C(t)$ and $\rho_H(t)$, are set to ensure that there is no mass flow through the walls:
\begin{equation}
\int_{0}^{\infty} v \maxw_C^{n} \, dv + \int_{-\infty}^{0} v \maxw_1^{n} \, dv = 0 \qquad \text{and} \qquad
\int_{0}^{\infty} v \maxw_{N_x}^{n} \, dv + \int_{-\infty}^{0} v \maxw_H^{n} \, dv = 0.
\end{equation}

Several boundary values need to be specified to implement the above boundary conditions in the micro-macro approach. In the micro update as described in Section \ref{subsec:micro-fvm-1d}, we need boundary values in \eqref{eqn5-19} and \eqref{eqn5-20}. The boundary temperature values  at the wall boundaries in \eqref{eqn5-19} are computed via kinetic flux vector splitting:
\begin{align}
T^n_{1/2} &= \frac{\int_{0}^{\infty} v^2 \maxw_C^{n} \, dv + \int_{-\infty}^{0} v^2 \maxw_1^{n} \, dv}{\int_{0}^{\infty} \maxw_C^{n} \, dv + \int_{-\infty}^{0} \maxw_1^{n} \, dv} 
= \left[\frac{\left( \sqrt{\frac{\pi T_C}{2}} - u \right) \rho \alpha + \left(T - u
\sqrt{\frac{\pi T_C}{2}} + u^2 \right) \rho  \beta^{-}}{\sqrt{\frac{2 \pi}{T_C}} \rho \left(  \alpha -  u \beta^{-} \right)}\right]_{ \, i=1}^{ \, n}, \\
T^n_{N_x + 1/2} &= \frac{\int_{0}^{\infty} v^2 \maxw_{N_x}^{n} \, dv + \int_{-\infty}^{0} v^2 \maxw_H^{n} \, dv}{\int_{0}^{\infty} \maxw_{N_x}^{n} \, dv + \int_{-\infty}^{0} \maxw_H^{n} \, dv} 
= \left[ \frac{\left(\sqrt{\frac{\pi T_H}{2}} + u \right) \rho \alpha + \left(T + u \sqrt{\frac{\pi T_H}{2}} + u^2 \right) \rho \beta^+}{\sqrt{\frac{2 \pi}{T_H}} \rho \left( \alpha + \beta^{+} \right)}
\right]_{ \, i=N_x}^{ \, n},
\end{align}
where $\alpha$ and $\beta^{\pm}$ are defined in \eqref{eqn:flux1d_2}.
In the micro update equation \eqref{eqn5-20}, all inflow fluxes are set to zero:
\begin{align}
  Z^n_{1k} = v_k^{-} \left( \frac{G^n_{2 k} - G^n_{1 k}}{\Delta x} \right)
 + v_k^{+} \left( \frac{G^n_{1 k} - \inred{0}}{\Delta x}\right), \quad
  Z^n_{N_x k} = v_k^{-} \left( \frac{\inred{0} - G^n_{N_x k}}{\Delta x} \right)
 + v_k^{+} \left( \frac{G^n_{N_x k} - G^n_{N_x-1 k}}{\Delta x}\right).
\end{align}

In the macro update equation \eqref{eqn5-27}, we need to specify fluxes at the wall; we accomplish this via kinetic flux vector splitting between fluid (interior) and wall (exterior) Maxwell-Boltzmann distributions:
\begin{alignat}{2}
\vec{\mathcal{F}^{ \, n}_{1/2}} &= \int_0^\infty v \, \vec{m}\left(v\right) \, \maxw^n_{C} \, dv
+ \int_{-\infty}^0 v \, \vec{m}\left(v \right) \, \maxw^n_1 \, dv, \quad  
&\rho_C^n =& \, \sqrt{\frac{2\pi}{T_C}} \left( \rho^n_{1} \alpha^n_{1} - \rho^n_{1} u^n_{1} \beta^{n-}_{1} \right), \\
\vec{\mathcal{F}^{ \, n}_{N_x+1/2}} &= \int_0^\infty v \, \vec{m}\left(v\right) \, \maxw^n_{N_x} \, dv
+ \int_{-\infty}^0 v \, \vec{m}\left(v \right) \, \maxw^n_{H} \, dv, \quad 
&\rho_H^n =& \, \sqrt{\frac{2\pi}{T_H}} \left( \rho^n_{N_x} \alpha^n_{N_x} + \rho^n_{N_x} u^n_{N_x} \beta^{n+}_{N_x} \right),
\end{alignat}
where $\alpha$ and $\beta^{\pm}$ are defined in \eqref{eqn:flux1d_2}.
Finally, we also need to specify the heat flux at the boundaries; for that, we take:
\begin{equation}
H^{n+1}_{0} = 0 \quad \text{and} \quad H^{n+1}_{N_x + 1} = 0 \quad \Longrightarrow \quad
H^{n+1}_{1/2} = \frac{1}{2} H^{n+1}_{1} \quad \text{and} \quad H^{n+1}_{N_x + 1/2} = \frac{1}{2} H^{n+1}_{N_x}.
\end{equation}

We run a series of numerical simulations all on $(x,v) \in \left[0,1\right] \times \left[-6,6\right]$ with $N_x = N_v = 129$ and with the initial condition:
\begin{equation}
\left(\rho, \, u, \, T \right)(t=0,x) = \left( 1, \, 0, \, 1 \right), \quad 
g(t=0,x,v) = 0.
\end{equation}
In all these runs, the final time is $t=100$ and the CFL number is set to $0.95$.
We use the diffusely reflecting wall boundary conditions with temperatures $T_C=1.0$ on the left and $T_H=1.2$ right wall, which are implemented in micro-macro as described above, and use the following collision parameter so that we arrive at a constant viscosity and heat conductivity:
\begin{equation}
\label{eqn:collision_parameter_1d_p}
\tau = p \quad \Longrightarrow \quad \mu \equiv 1, \quad \kappa \equiv 1.5.
\end{equation}

The results of one such run with $\varepsilon=10^{-2}$ are shown in Figure \ref{fig:1dtest3a}. The panels of
Figure \ref{fig:1dtest3a} show the (a) full PDF (the plot only shows $v\in[-3,3]$): $f= \maxw+\varepsilon g$, (b) microscopic portion (the plot only shows $v\in[-3,3]$): $g$, (c) scalar temperature: $T$, and (d) heat flux divided by $\varepsilon$: $h/\varepsilon$. These results are consistent with similar results shown in the literature (e.g., Abdelmalik et al. \cite{Abdelmalik2017ErrorEquation}).
This figure can be generated using the Python/C++ companion code \cite{code:RossmanithSar2025} by going to the code directory and typing:
\begin{tcolorbox}
\begin{verbatim}
cd 1d/example3-heat-transfer; python make_paper_figures.py
\end{verbatim}
\end{tcolorbox}

The results of a series of runs with 15 different Knudsen numbers, ranging from $10^{-2}$ to 
$10^2$, are shown in Figure \ref{fig:1dtest3b}. Shown in Panel (a) of Figure \ref{fig:1dtest3b} is the absolute value of the scaled heat flux at $x=0.5$ at time $t=100$:
\begin{equation}
|h_1| :=  \frac{1}{2} \left| \int_{-\infty}^{\infty} v^3 g\left(t=100,x=0.5,v\right) \, dv \right|.
\end{equation}
Following van der Woude et al. \cite{article:vanderWoude2024}, and as a point of comparison, we include in Figure \ref{fig:1dtest3b}(a) three additional scaled heat fluxes:
\begin{align}
 \text{Navier-Stokes-Fourier:}& \quad \left| h_1^{\text{NSF}} \right| \sim \kappa T_{,x} = 
 \kappa \left( T_H - T_C \right) = 0.3, \\
 \text{Free molecular flow:}& \quad 
    \left| h_1^{\text{FreeM}} \right| \sim \frac{\sqrt{T_H} - \sqrt{T_C}}{\varepsilon}
       \sqrt{\frac{2 T_C T_H}{\pi}} 
    \approx \frac{0.083422728489775170621}{\varepsilon}, \\
  \text{Harmonic average:}& \quad \left| h_1^{\text{Ave}} \right| \sim 
  \frac{\left| h_1^{\text{NSF}} \right| \left| h_1^{\text{FreeM}} \right|}{\left| h_1^{\text{NSF}} \right| + \left| h_1^{\text{FreeM}} \right|} \approx 
  \frac{0.0834227284897751706}{0.2780757616325839021 + \varepsilon}. 
\end{align}
Figure \ref{fig:1dtest3b}(a) shows that the scaled heat flux at $x=0.5$ is close to the Navier-Stokes-Fourier scaled heat flux for small $\varepsilon$, close to the free-molecular scaled heat flux for large $\varepsilon$, and approximately equal to the harmonic average for all $\varepsilon$. This result is consistent with van der Woude et al. \cite{article:vanderWoude2024}. In Figure \ref{fig:1dtest3b}(b), we show the distribution at $t=100$ and $x=0.5$ for $v\in[-6,6]$ with $\varepsilon = 10^{30}$ alongside the exact free molecular flow distribution function. The results between the numerical computation via micro-macro and the exact solution are very close, and in particular, the jump across $v=0$ is well approximated.
This figure can be generated using the Python/C++ companion code \cite{code:RossmanithSar2025} by going to the code directory and typing:
\begin{tcolorbox}
\begin{verbatim}
cd 1d/example4-HT-Knudsen; python make_paper_figures.py
\end{verbatim}
\end{tcolorbox}


\section{2D Micro-Macro Decomposition Scheme}
\label{sec:micmac2D}
In this section, we extend the 1D micro-macro scheme from \S\ref{sec:micmac1D} to 2D (i.e., $d_x = d_v = 2$) and to the ES-BGK collision operator. To handle the ES-BGK collision operator, we solve the extended compressible Euler equations (i.e., 6-moment model when $d_v=2$) in the macro step so that we have access to the full pressure tensor.

\subsection{Mathematical Formulation}
The two-dimensional ($d_v=2$ and $d_x=2$) Boltzmann-ES-BGK equation can be written as:
\begin{equation}
\label{eq5-28}
    f_{,t} + \vec{v} \cdot \vec{\nabla}_{\vec{x}} f = \frac{\tau}{\varepsilon} \left( \gauss\left[f\right] - f \right),
    \quad
    \gauss\left[f\right] := \frac{\rho}{\sqrt{\text{det}\left(2 \, \pi \, \mat{\mathcal{T}} \, \right)}}\exp\left(-\frac{1}{2} {\left( \, \vec{v} - \vec{u} \, \right) \cdot \mat{\mathcal{T}}^{-1} \left( \, \vec{v}-\vec{u} \, \right)}\right),
\end{equation}
where $\gauss\left[f\right]$ is a Gaussian distribution with the same mass, momentum, and scalar energy as $f$, but with a different temperature tensor given by \eqref{eq2-19}. The temperature tensor of $f$ is $\mat{{\mathbb T}}$ (see \eqref{eqn:temp_tensor}). The trace of $\mat{{\mathbb T}}$ and $\mat{\mathcal T}$ are equal and proportional to the scalar temperature: $T$ (see definitions in Equation \eqref{eqn:temp_tensor}). The pressure tensors for the Gaussian, $\gauss\left[f\right]$, and $f$ are also different:
\begin{equation}
    \label{eqn:mod-press-esbgk}
   \mat{\mathcal{P}} =  \left(1-\nu \right) \, p \, \mat{\mathcal{I}} + \nu \, \mat{{\mathbb P}} \, ,
   \quad \text{with} \quad  
   -1 \le \nu < 1 \,,
\end{equation}
where $\mat{\mathcal{P}}$ and $\mat{{\mathbb P}}$ are defined in \eqref{eqn:pressure-tensor-esbgk}, and
the scalar pressure is defined in \eqref{eqn:scalar-pressure}.

\begin{remark}
From \eqref{eq5-28} it appears at first glance that the equilibrium of the ES-BGK collision operator is $f\equiv \gauss\left[f\right]$; however, by construction $f$ and $\gauss\left[f\right]$ are not equivalent, because they have different temperature tensors: $\mat{\mathcal{T}} \not\equiv \mat{{\mathbb T}}$ \, for $-1 \le \nu < 1$. The only exception to this non-equivalence is when $\mat{{\mathbb T}} \equiv T \mat{\mathcal{I}}$, in which case $\mat{\mathcal{T}} \equiv T \mat{\mathcal{I}}$, which follows from \eqref{eqn:temp_tensor}; this exception is precisely the case when $f\equiv \maxw\left[f\right]$, which means the equilibrium of \eqref{eq5-28} is again the Maxwell-Boltzmann distribution (see \eqref{eq5-31} below).
\end{remark}

Similar to the 1D case, in the 2D micro-macro approach, the distribution function $f$ is written as the sum of the macroscopic Maxwell-Boltzmann distribution and the microscopic portion denoted by $g$:
\begin{equation}
\label{eq5-31}
f = \maxw\left[f\right]
   + \varepsilon \, g, \qquad
   \maxw\left[f\right] := \frac{\rho}{2\pi T} \exp\left[-\frac{\left\| \vec{v} - \vec{u} \right\|^2}{2T}\right].
\end{equation}
Inserting the micro-macro decomposition \eqref{eq5-31} into the Boltzmann-ES-BGK equation \eqref{eq5-28} yields:
\begin{gather}
\label{eq5-32}
\maxw_{,t} + \vec{v} \cdot \vec{\nabla}_{\vec{x}} \maxw + \varepsilon g_{,t}  +
\varepsilon \, \vec{v} \cdot \vec{\nabla}_{\vec{x}} g = \frac{\tau}{\varepsilon} \left( \gauss - \maxw \right) - \tau \, g .
\end{gather}
Computing the first six moments of \eqref{eq5-32}, namely $\left(1, v_1, v_2, v_2^2, v_1 v_2, v_2^2 \right)$, yields the following fluid system:
\begin{align}
\label{eq5-40}
\left(\text{2D Macro}\right): \quad
  \begin{bmatrix} 
	\rho \\ 
	\rho u_1 \\ 
	\rho u_2 \\ 
	\rho u_1^2 + \pr_{11} \\
	\rho u_1 u_2 + \pr_{12} \\
	\rho u_2^2 + \pr_{22}
	\end{bmatrix}_{,t} + \vec{f}\left(\vec{q}\right)_{,x} + \vec{g}\left(\vec{q}\right)_{,y} = \frac{\tau \left(1-\nu\right)}{2\varepsilon} 
\begin{bmatrix}
0 \\ 0 \\ 0 \\ \pr_{22} - \pr_{11} \\ -2 \pr_{12} \\ \pr_{11} - \pr_{22}
\end{bmatrix},
\end{align}
where
\begin{equation}
\label{eq5-41}
\begin{gathered}
		\vec{f}\left(\vec{q}\right) = \begin{bmatrix} 
	\rho u_1 \\ 
	\rho u_1^2 + \pr_{11} \\
	\rho u_1 u_2 + \pr_{12} \\
	\rho u_1^3 + 3 u_1 \pr_{11} + {\mathbb H}_{111} \\
	\rho u_1^2 u_2 + u_2 \pr_{11} + 2 u_1 \pr_{12}  + {\mathbb H}_{112} \\
	\rho u_1 u_2^2 + 2 u_2 \pr_{12} + u_1 \pr_{22}  + {\mathbb H}_{122}
	\end{bmatrix} \qquad \text{and} \qquad
	\vec{g}\left(\vec{q}\right) = \begin{bmatrix} 
	\rho u_2 \\ 
    \rho u_1 u_2 + \pr_{12} \\
	\rho u_2^2 + \pr_{22} \\
	\rho u_1^2 u_2 + u_2 \pr_{11} + 2 u_1 \pr_{12} + {\mathbb H}_{112} \\
    \rho u_1 u_2^2 + 2 u_2 \pr_{12} + u_1 \pr_{22} + {\mathbb H}_{122}  \\
    \rho u_2^3 + 3 u_2 \pr_{22} + {\mathbb H}_{222}
	\end{bmatrix},
\end{gathered}
\end{equation}
which is not closed since the heat flux tensor must still be supplied via the microscopic portion $g$:
\begin{equation}
   \begin{pmatrix} {\mathbb H}_{111} \\ {\mathbb H}_{112} \\
     {\mathbb H}_{122} \\ {\mathbb H}_{222} \end{pmatrix}
      = \varepsilon \iint_{-\infty}^{\infty} \begin{pmatrix} \left(v_1-u_1\right)^3 \\
         \left(v_1-u_1\right)^2 \left(v_2-u_2\right) \\ \left(v_1-u_1\right) \left(v_2-u_2\right)^2 \\ 
         \left(v_2-u_2\right)^3 \end{pmatrix} \, g \, dv_1\,dv_2.
\end{equation}
In \eqref{eq5-40}, $\varepsilon$ is the Knudsen number, $\nu$ is the Prandtl number, and $\tau$ is the collision parameter.

To close fluid system \eqref{eq5-40}--\eqref{eq5-41}, we need to obtain an evolution equation for the microscopic portion, $g$, from equation \eqref{eq5-32}. Similar to the $d_v=1$ case, a basis for the nullspace of \eqref{eq2-15} is given by \cite{Bennoune2008UniformlyAsymptotics}:
\begin{equation}
\label{eqn:nullspace_basis_2d}
    \mathcal{B} = \left\{ \, \frac{1}{\rho}\mathcal{M}, \quad \left(\frac{v_1-u_1}{\sqrt{T}} \right) \frac{1}{\rho}\mathcal{M}, 
    \quad \left( \frac{v_2-u_2}{\sqrt{T}} \right) \frac{1}{\rho}\mathcal{M}, \quad \left(\frac{\left(v_1 - u_1 \right)^2 + \left(v_2 - u_2 \right)^2}{2T}-1\right)\frac{1}{\rho}\mathcal{M} \, \right\},
\end{equation}
which is orthogonal with respect to the weighted inner product:
\begin{equation}
 \left( \psi_1, \, \psi_2 \right)_{\maxw^{-1}} := \iint_{-\infty}^{\infty} \psi_1 \, \psi_2 \, \maxw^{-1} \, dv_1 \, dv_2.
\end{equation} 
An orthogonal projection onto this basis for some function, $F\left(\vec{v}\right):\reals^2 \mapsto \reals$, such that $F \maxw^{-1/2} \in L^2\left(\reals^2\right)$, is given by:
\begin{equation}
\label{eq5-33}
 \Pi_{\maxw}\left[ F \right]  = \left\{ A_1 + \left(\frac{v_1 - u_1}{\sqrt{T}}\right) A_2 
 + \left( \frac{v_2 - u_2}{\sqrt{T}} \right) A_3 + \left( \frac{\left(v_1 - u_1 \right)^2 + \left(v_2 - u_2 \right)^2}{2T} - 1\right) A_4
 \right\} \maxw[f],
\end{equation}
where
\begin{align}
\label{eq5-33-plus}
\left(  A_1, \, A_2, \, A_3, \, A_4 \right) &=
\frac{1}{\rho} \iint_{-\infty}^{\infty} \left( \, 1, \, \frac{v_1 - u_1}{\sqrt{T}}, \, \frac{v_2 - u_2}{\sqrt{T}}, \, \frac{\left(v_1 - u_1 \right)^2 + \left(v_2 - u_2 \right)^2}{2T} - 1 \right) F  \, dv_1 \, dv_2.
\end{align}
This projector has the property that (see Bennoune et al. \cite{Bennoune2008UniformlyAsymptotics} for details)
\begin{equation}
\label{eq5-38}
  \Pi_{\maxw} \left( \maxw \right) = \maxw, \quad
  \Pi_{\maxw} \left( \gauss \right) = \maxw, \quad
  \Pi_{\maxw} \left( \maxw_{,t} \right) = \maxw_{,t}, \quad
  \Pi_{\maxw} \left( g_{,t} \right) = 0, \quad
  \Pi_{\maxw} \left( g \right) = 0,
\end{equation}
and when applied to \eqref{eq5-32} yields:
\begin{equation}
\label{eqn:projected_micmac_2D}
\maxw_{,t} + \Pi_{\maxw} \left[ \vec{v} \cdot \vec{\nabla}_{\vec{x}} \maxw \right]  + \varepsilon \, 
\Pi_{\maxw} \left[ \vec{v} \cdot \vec{\nabla}_{\vec{x}} g \right] =  0.
\end{equation}
Subtracting \eqref{eqn:projected_micmac_2D} from \eqref{eq5-32}, dividing by $\varepsilon$, and rearranging terms yields an evolution equation for the microscopic portion $g$:
\begin{gather}
\label{eq5-39}
   g_{,t}  +
\left( {\mathcal I} - 
 \Pi_{\maxw} \right) \left(\vec{v} \cdot \vec{\nabla}_{\vec{x}} g \right) = -\frac{1}{\varepsilon} \biggl( \tau \, g 
+ \left( {\mathcal I} - \Pi_{\maxw} \right) \left(\vec{v} \cdot \vec{\nabla}_{\vec{x}} \maxw \right) - 
\frac{\tau}{\varepsilon} \left( \gauss - \maxw \right) \biggr),
\end{gather}
where ${\mathcal I}$ is the identity operator. Note that for $0 < \varepsilon \ll 1$, $\gauss - \maxw = {\mathcal O}\left(\varepsilon\right)$, which means the last term in the parenthesis on the right-hand side of \eqref{eq5-39} is 
${\mathcal O}\left(1\right)$ for $0 < \varepsilon \ll 1$.

Just as in 1D, equation \eqref{eq5-39} can be further simplified, since it turns out that it is possible to analytically compute the orthogonal projection of the Maxwellian transport term (again, see Bennoune et al. \cite{Bennoune2008UniformlyAsymptotics} for details):
\begin{equation}
\label{eq5-69}
\begin{gathered}
 \left( {\mathcal I} - \Pi_{\maxw} \right) \left( \, \vec{v} \cdot \vec{\nabla}_{\vec{x}} \maxw \, \right)  = \Bigr\{ \, \mat{B} : \mat{\sigma} + \vec{C} \cdot
 {\vec{\nabla} T}  \, \Bigl\} \maxw, \qquad \text{where} \qquad
 \mat{\sigma} =
 \begin{bmatrix}[1.5]
 \frac{\partial u_1}{\partial x} -  \frac{\partial u_2}{\partial y} &   \frac{\partial u_1}{\partial y} 
    + \frac{\partial u_2}{\partial x} \\
  \frac{\partial u_1}{\partial y} + \frac{\partial u_2}{\partial x} & - \frac{\partial u_1}{\partial x} 
     + \frac{\partial u_2}{\partial y}
 \end{bmatrix}, \\
\vec{C} = \left( \frac{\|\vec{v} - \vec{u} \|^2}{2T} - 2\right) \frac{\vec{v} - \vec{u}}{T}, \qquad
\text{and} \qquad
\mat{B} = \frac{1}{2T} \begin{bmatrix}[1.5]
 -\left( v_2 - u_2 \right)^2 & \left(v_1 - u_1 \right) \left( v_2 - u_2 \right) \\
\left(v_1- u_1 \right) \left( v_2 - u_2 \right) & -\left( v_1 - u_1 \right)^2
\end{bmatrix}.
 \end{gathered}
\end{equation}
In \eqref{eq5-69} we used the following shorthand for the Frobenius inner product:
\begin{equation}
	\mat{B} : \mat{\sigma} = \sum_{\ell=1}^{2} \sum_{k=1}^{2} B_{\ell k} \sigma_{\ell k}.
\end{equation}
Substituting \eqref{eq5-69} into \eqref{eq5-39} results in the following micro evolution equation:
\begin{equation}
\label{eqn5-75}
\left(\text{2D Micro}\right): \quad
    g_{,t}  +
\left( {\mathcal I} - 
 \Pi_{\maxw} \right) \left(\vec{v} \cdot \vec{\nabla}_{\vec{x}} g \right) 
 = -\frac{\tau}{\varepsilon} \biggl( g - \hat{g} \biggr), \quad
 \hat{g} = -\frac{1}{\tau} \Bigr\{ \, \mat{B}: \mat{\sigma} + \vec{C} \cdot
 {\vec{\nabla} T} \, \Bigl\} \maxw +
\frac{1}{\varepsilon} \left( \gauss - \maxw \right). 
\end{equation}
For $0 < \varepsilon \ll 1$, $\hat{g} = {\mathcal O}\left(1\right)$, and the conclusion of Remark \ref{note:note-NSF-limit}
again applies in this 2D case.

\subsection{Setup for 2D Micro-Macro Finite Volume Method}
Now that the micro-macro formulation for $d_v=2$ has been established via \eqref{eq5-40} and \eqref{eqn5-75}, we turn our attention to the numerical discretization. We introduce a uniform, cell-centered, finite volume mesh in the domain $\left(x_{\text{min}}, x_{\text{max}} \right) \times \left(y_{\text{min}}, y_{\text{max}} \right)
\times \left(\left[v_{1}\right]_{\text{min}}, \left[v_{1}\right]_{\text{max}} \right) \times \left(\left[v_{2}\right]_{\text{min}}, \left[v_{2}\right]_{\text{max}} \right)$:
\begin{equation}
\begin{gathered}
   x_i := x_{\text{min}} + \left( i - \frac{1}{2} \right){\Delta x} \quad \forall i\in\left[1,N_x\right], \quad
   y_j := y_{\text{min}} + \left( j - \frac{1}{2} \right){\Delta y} \quad \forall j=\in\left[1,N_y\right], \\
   \left[v_1\right]_k := \left[v_1\right]_{\text{min}} + \left( k - \frac{1}{2} \right){\Delta v_1} \quad \forall k\in\left[1,N_{v_1}\right], \quad
   \left[v_2\right]_{\ell} := \left[v_2\right]_{\text{min}} + \left( \ell - \frac{1}{2} \right){\Delta v_2} \quad \forall \ell\in\left[1,N_{v_2}\right],
\end{gathered}
\end{equation}
where
\begin{equation}
{\Delta x} := \frac{x_{\text{max}}-x_{\text{min}}}{N_x}, \quad
{\Delta y} := \frac{y_{\text{max}}-y_{\text{min}}}{N_y}, \quad
{\Delta v_1} := \frac{\left[v_1\right]_{\text{max}}-\left[v_1\right]_{\text{min}}}{N_{v_1}}, \quad
{\Delta v_2} := \frac{\left[v_2\right]_{\text{max}}-\left[v_2\right]_{\text{min}}}{N_{v_2}}.
\end{equation}
At each time level, $t^n = n \Delta t$, the macroscopic and microscopic portions of the distribution functions are given:
\begin{equation}
\vec{Q^n_{ij}} := \left( \, \rho^n_{ij}, \, \rho \vec{u^n_{ij}}, \, \mat{{\mathbb E}^n_{ij}} \, \right) \quad  
 \text{and} \quad
     G^n_{ijk\ell} \quad 
\forall i\in\left[1,N_x\right], \, \, \, \forall j\in\left[1,N_y\right], \, \, \, \forall k\in\left[1,N_{v_1}\right], \, \, \, \forall \ell\in\left[1,N_{v_2}\right],
\end{equation}
where
\begin{align}
\vec{Q^n_{ij}} &\approx \frac{1}{4} \iint_{-1}^{1} \vec{q}\left(t^n, \, x_i + \frac{\Delta x \, \xi}{2}, \,
y_j + \frac{\Delta y \, \eta}{2} \right) \, d\xi \, d\eta, \\
G^n_{ijk\ell} &\approx \frac{1}{16} \iiiint_{-1}^{1} g\left(t^n, \, x_i +  \frac{\Delta x \, \xi}{2}, \,
y_j + \frac{\Delta y \, \eta}{2}, \, \left[v_1\right]_k + \frac{\Delta v_1 \, \mu_1}{2},
\, \left[v_2\right]_{\ell} + \frac{\Delta v_2 \, \mu_2}{2} \right) \, d\xi \, d\eta \, d\mu_1 \, d\mu_2.
\end{align}
The solution is initialized at time $t=0$ and runs forward until $t=T_{\text{final}}$. 
The time step is chosen as a global constant and is defined as follows:
\begin{equation}
\label{eqn:CFL-2D-part1}
  \Delta t = \text{CFL} \cdot \text{min} \left\{ \frac{\Delta x}{\mathcal{V}_1}, \,
      \frac{\Delta y}{\mathcal{V}_2} \right \}, \quad
  N_{\text{steps}} = \left\lceil \frac{T_{\text{final}}}{\Delta t} \right\rceil \quad \Longrightarrow \quad
   \Delta t \leftarrow \frac{T_{\text{final}}}{N_{\text{steps}}}, \quad 
  \text{CFL} \leftarrow  \max\left\{ \frac{\mathcal{V}_1 \, \Delta t}{\Delta x}, \, \frac{\mathcal{V}_2 \, \Delta t}{\Delta y}\right\},
\end{equation}
where $N_{\text{steps}}$ is the total number of timesteps needed to advance from $t=0$ to $t=T_{\text{final}}$ and
\begin{equation}
\label{eqn:CFL-2D-part2}
\left(\text{CFL number}\right): \quad 0 < \text{CFL} < 1, \qquad
\left(\text{max speed}\right): \quad 
\begin{cases}
\mathcal{V}_{1} := \max\left\{ \bigl| \left[v_1\right]_{\text{min}} \bigr|, \bigl|\left[v_1\right]_{\text{max}}\bigr| \right\}, \\
\mathcal{V}_{2} := \max\left\{ \bigl| \left[v_2\right]_{\text{min}} \bigr|, \bigl|\left[v_2\right]_{\text{max}}\bigr| \right\}.
\end{cases}
\end{equation}

\subsection{2D Micro Finite Volume Method}
In this section, we demonstrate how to update the current value of the micro portion, $G^n$, to its value at the next time step, $G^{n+1}$:
\begin{equation}
   \text{Input:} \quad \vec{Q^n_{ij}} = \left( \, \rho^n_{ij}, \, \rho \vec{u^n_{ij}}, \, \mat{{\mathbb E}^n_{ij}} \, \right) \quad \text{and} \quad
     G^n_{ijk\ell} \quad \Longrightarrow \quad \text{Output:} \quad G^{n+1}_{ijk\ell}.
\end{equation}
To guarantee stability up to a CFL number of one as defined in \eqref{eqn:CFL-2D-part1}, we elect
to solve the micro equation \eqref{eqn5-75} using operator splitting:
\begin{alignat}{2}
\label{eqn:micro-2D-split-1}
  \Delta t:& \quad g_{,t}  =  \Pi_{\maxw} \left( v_1 \, g_{,x} \right) - v_1 \, g_{,x} & \qquad & 
  \left( \, G^n \rightarrow G^{\star} \, \right), \\
\label{eqn:micro-2D-split-2}
 \Delta t:& \quad g_{,t}  =  \Pi_{\maxw} \left( v_2 \, g_{,y} \right) - v_2 \, g_{,y}  & \qquad & 
 \left( \, G^{\star} \rightarrow G^{\star \star} \, \right), \\
\label{eqn:micro-2D-split-3}
 \Delta t:& \quad g_{,t}
 = -\frac{\tau}{\varepsilon} \biggl( g - \hat{g} \biggr) & \qquad & \left( \, G^{\star \star} \rightarrow G^{n+1} \, \right). 
\end{alignat}

For the first part of the operator splitting as described by \eqref{eqn:micro-2D-split-1}, we first approximate the term $v \, g_{,x}$ via simple upwinding:
\begin{equation}
\label{eqn:zmicrox}
Z^n_{{ijk\ell}} = \left[v_1\right]_{k}^{-} \left( \frac{G^n_{\left(i+1, j, k, \ell\right)} - G^n_{\left(i, j, k, \ell\right)}}{\Delta x} \right)
+ \left[v_1\right]_{k}^{+} \left( \frac{G^n_{\left(i, j, k, \ell\right)} - G^n_{\left(i-1, j, k, \ell\right)}}{\Delta x} \right).
\end{equation}
We also compute the projection of this term on the nullspace of the collision operator via
\begin{equation}
\label{eqn:zhat-2D-xdir}
\widehat{Z}^n_{ijk\ell} = \left\{ \left[A_1\right]^n_{ij}
 + \left(\frac{\vec{v_{k\ell}} - \vec{u^n_{ij}}}{\sqrt{T^n_{ij}}}\right) \cdot \begin{bmatrix}
  A_2 \\ A_3 \end{bmatrix}^n_{ij}
+ \left( \frac{\Bigl\| \vec{v_{k\ell}} - \vec{u^n_{ij}} \Bigr\|^2}{2 T^n_{ij}} - 1\right) \left[A_4\right]^n_{ij} \right\} \maxw^n_{ijk\ell},
\end{equation}
where the coefficients are computed via numerical quadrature using the midpoint rule:
\begin{equation}
\begin{split}
\iint_{-\infty}^{\infty} g\left(v_1,v_2\right) \, d\vec{v} &\approx
			\int_{\left[v_1\right]_{\text{min}}}^{\left[v_1\right]_{\text{max}}}
			\int_{\left[v_2\right]_{\text{min}}}^{\left[v_2\right]_{\text{max}}}
			 g\left(v_1,v_2\right) \, d\vec{v} =
			\sum_{k=1}^{N_{v_1}} \sum_{\ell=1}^{N_{v_2}} \int_{\left[v_1\right]_k - \frac{\Delta v_1}{2}}^{\left[v_1\right]_k + \frac{\Delta v_1}{2}} \int_{\left[v_2\right]_{\ell} - \frac{\Delta v_2}{2}}^{\left[v_2\right]_{\ell} + \frac{\Delta v_2}{2}} g\left(v_1,v_2\right) \, d\vec{v} \\
			&\approx 
			\Delta v_1 \Delta v_2 \sum_{k=1}^{N_{v_1}} \sum_{\ell=1}^{N_{v_2}} g\Bigl(\left[v_1\right]_k, \left[v_2\right]_{\ell} \Bigr)
\end{split}
\end{equation}
such that
\begin{equation}
\label{eqn:coeffs-2D-xdir}
\left( \, \left[A_1\right]^n_{ij}, \, \, \begin{bmatrix} A_2 \\ A_3 \end{bmatrix}^n_{ij}, \, \, \left[A_4\right]^n_{ij} \, \right) = \frac{\Delta v_1 \, \Delta v_2}{\rho^n_{ij}} \sum_{k=1}^{N_{v_1}}\sum_{\ell=1}^{N_{v_2}} 
\left( \, 1, \, \, \frac{\vec{v_{k\ell}} - \vec{u^n_{ij}}}{\sqrt{T^n_{ij}}},  \, \,
\frac{\Bigl\| \vec{v_{k\ell}} - \vec{u^n_{ij}} \Bigr\|^2}{2 T^n_{ij}} - 1 \, \right) Z^n_{ijk\ell}.
\end{equation}
Using forward Euler in time yields the following approximation to \eqref{eqn:micro-2D-split-1}:
\begin{equation}
\label{eqn:update-micro-2D-xdir}
G^{\star}_{ijk\ell} =  G^n_{ijk\ell} + \Delta t \left( \widehat{Z}^n_{ijk\ell} - Z^n_{ijk\ell} \right).
\end{equation}

For the second part of the operator splitting as described by \eqref{eqn:micro-2D-split-2}, we first approximate the term $v \, g_{,y}$ via simple upwinding:
\begin{equation}
\label{eqn:zmicroy}
Z^{\star}_{ijk\ell} = \left[v_2\right]^{-}_{\ell} 
\left( \frac{G^{\star}_{\left(i, j+1, k, \ell\right)} - G^{\star}_{\left(i,j,k,\ell\right)}}{\Delta y} \right)+ \left[v_2\right]^{+}_{\ell} \left( \frac{G^{\star}_{\left(i, j, k, \ell\right)} - G^{\star}_{\left(i, j-1, k, \ell\right)}}{\Delta y} \right).
\end{equation}
We also compute the projection of this term on the nullspace of the collision operator via
\begin{equation}
\label{eqn:zhat-2D-ydir}
\widehat{Z}^{\star}_{ijk\ell} = \left\{ \left[A_1\right]^{\star}_{ij}
 + \left(\frac{\vec{v_{k\ell}} - \vec{u^n_{ij}}}{\sqrt{T^n_{ij}}}\right) \cdot \begin{bmatrix}
  A_2 \\ A_3 \end{bmatrix}^{\star}_{ij}
+ \left( \frac{\Bigl\| \vec{v_{k\ell}} - \vec{u^n_{ij}} \Bigr\|^2}{2 T^n_{ij}} - 1\right) \left[A_4\right]^{\star}_{ij} \right\} \maxw^n_{ijk\ell},
\end{equation}
where
\begin{equation}
\label{eqn:coeffs-2D-ydir}
\left( \, \left[A_1\right]^{\star}_{ij}, \, \, \begin{bmatrix} A_2 \\ A_3 \end{bmatrix}^{\star}_{ij}, \, \, \left[A_4\right]^{\star}_{ij} \, \right) = \frac{\Delta v_1 \, \Delta v_2}{\rho^n_{ij}} \sum_{k=1}^{N_{v_1}}\sum_{\ell=1}^{N_{v_2}} 
\left( \, 1, \, \, \frac{\vec{v_{k\ell}} - \vec{u^n_{ij}}}{\sqrt{T^n_{ij}}},  \, \,
\frac{\Bigl\| \vec{v_{k\ell}} - \vec{u^n_{ij}} \Bigr\|^2}{2 T^n_{ij}} - 1 \, \right) Z^{\star}_{ijk\ell}.
\end{equation}
Using forward Euler in time yields the following approximation to \eqref{eqn:micro-2D-split-2}:
\begin{equation}
\label{eqn:update-micro-2D-ydir}
G^{\star\star}_{ijk\ell} =  G^{\star}_{ijk\ell} + \Delta t \left( \widehat{Z}^{\star}_{ijk\ell} - Z^{\star}_{ijk\ell} \right).
\end{equation}

To complete the final part of the operator splitting as described by \eqref{eqn:micro-2D-split-3}, we first construct an approximate $\hat{g}$ from \eqref{eqn5-75} on the numerical grid:
\begin{align}
\label{eqn5-97} 
\begin{split}
\widehat{G}_{ijk\ell}^{n} &= -\frac{1}{\tau^n_{ij}} 
 \Bigg(\mat{B^n_{ijk\ell}}:\mat{\sigma^n_{ij}} \, + \,  \vec{C^n_{ijk\ell}} \cdot \left( \frac{T^n_{\left( i+1/2, \, j \right)}-T^n_{\left( i-1/2, \, j \right)}}{\Delta x}, \, \frac{T^n_{\left( i, \, j+1/2 \right)}-T^n_{\left( i, \, j-1/2 \right)}}{\Delta y} \right) \Bigg) \, \maxw^n_{ijk\ell} \\ 
 &+ \frac{1}{\varepsilon}\left(\mathcal{G}^n_{ijk\ell} - \maxw^n_{ijk\ell}\right), \qquad
 T^n_{\left( i+1/2, \, j \right)} = \frac{1}{2} \left( T^n_{\left(i+1, j\right)} + T^n_{ij} \right), \qquad
 T^n_{\left( i, \, j+1/2 \right)} = \frac{1}{2} \left( T^n_{\left(i, j+1\right)} + T^n_{ij} \right),
\end{split}
\end{align}
where
\begin{equation}
\label{eqn5-97more}
\mat{\sigma^n_{ij}} :=  \begin{bmatrix}[2]
 \frac{\Delta^x_{ij} u^n_1}{\Delta x} -  \frac{\Delta^y_{ij} u^n_2}{\Delta y} &   \frac{\Delta^y_{ij} u^n_1}{\Delta y} 
    + \frac{\Delta^x_{ij} u^n_2}{\Delta x} \\
  \frac{\Delta^y_{ij} u^n_1}{\Delta y} + \frac{\Delta^x_{ij} u^n_2}{\Delta x} & - \frac{\Delta^x_{ij} u^n_1}{\Delta x} 
     + \frac{\Delta^y_{ij} u^n_2}{\Delta y}
     %
 \end{bmatrix}
 \quad
 \Delta^x_{ij} u^n := \frac{u^n_{(i+1, \, j)}-u^n_{(i-1, \, j)}}{2}, \quad
 \Delta^y_{ij} u^n := \frac{u^n_{(i, \, j+1)}-u^n_{(i, \, j-1)}}{2},
\end{equation}
and $\mat{B}$ and $\mat{C}$ are defined in \eqref{eq5-69}.
Lastly, the numerical micro solution, $G^{\star\star}$, is updated to its new value, $G^{n+1}$,
 by applying to micro equation \eqref{eqn:micro-2D-split-3} a backward Euler time discretization:
\begin{align}
\label{eqn5-98}
G^{n+1}_{ijk\ell} &= \left( \frac{\varepsilon}{\varepsilon + \Delta t \, \tau^n_{ij}} \right) G^{\star\star}_{ijk\ell} 
+\left( \frac{\Delta t \, \tau^n_{ij}}{\varepsilon + \Delta t \, \tau^n_{ij}} \right) \widehat{G}^{n}_{ijk\ell}.
\end{align}
Note that this update is {\it asymptotic-preserving} \cite{Jin2012AsymptoticReview}, since we can take the limit
 $\varepsilon \rightarrow 0^+$ of this expression for fixed mesh parameters, $\Delta t$, $\Delta x$, $\Delta y$, 
 $\Delta v_1$, and $\Delta v_2$,
 and arrive at the following result:
 \begin{align}
  G^{n+1}_{ijk\ell} \rightarrow \widehat{G}^{n}_{ijk\ell} \quad \text{as} \quad \varepsilon \rightarrow 0^+.
\end{align}


\subsubsection{Extrapolation Boundary Conditions in the Micro Step}
We briefly mention how simple boundary conditions are implemented in the 2D micro update. For extrapolation boundary conditions, in the first step of operator splitting, we do the following:
\begin{align}
\label{eqn:2d_extrap_micro_1}
 \text{Left BC:} & \quad Z^n_{{\left(1,j,k,\ell \right)}} = \left[v_1\right]_{k}^{-} \left( \frac{G^n_{\left(2, j, k, \ell\right)} - G^n_{\left(1, j, k, \ell\right)}}{\Delta x} \right), \\
\label{eqn:2d_extrap_micro_2}
 \text{Right BC:} & \quad Z^n_{\left(N_x,j,k,\ell\right)} = \left[v_1\right]_{k}^{+} \left( \frac{G^n_{\left(N_x, j, k, \ell\right)} - G^n_{\left(N_x-1, j, k, \ell\right)}}{\Delta x} \right).
\end{align}
In the second step of operator splitting, we set
\begin{align}
\label{eqn:2d_extrap_micro_3}
\text{Bottom BC:} & \quad 
Z^{\star}_{\left(i,1,k,\ell\right)} = \left[v_2\right]^{-}_{\ell} 
\left( \frac{G^{\star}_{\left(i, 2, k, \ell\right)} - G^{\star}_{\left(i,1,k,\ell\right)}}{\Delta y} \right), \\
\label{eqn:2d_extrap_micro_4}
\text{Top BC:} & \quad
Z^{\star}_{\left(i,N_y,k,\ell\right)} = \left[v_2\right]^{+}_{\ell} \left( \frac{G^{\star}_{\left(i, N_y, k, \ell\right)} - G^{\star}_{\left(i, N_y-1, k, \ell\right)}}{\Delta y} \right).
\end{align}
In the third and final step of operator splitting, we set
\begin{equation}
\label{eqn:2d_extrap_micro_5}
\begin{gathered}
 T^n_{\left(1/2, j\right)}     = T^n_{\left(1, j\right)}, \quad
 T^n_{\left(N_x+1/2, j\right)} = T^n_{\left(N_x, j\right)}, \quad
 T^n_{\left(i, 1/2\right)}     = T^n_{\left(i, 1\right)}, \quad
 T^n_{\left(i, N_y+1/2\right)} = T^n_{\left(i, N_y\right)}, \\
 \Delta^x_{\left(1,j\right)} u^n = \frac{1}{2} \left( u^n_{\left(2,j\right)}-u^n_{\left(1,j\right)} \right), \quad
 \Delta^x_{\left(N_x,j\right)} u^n = \frac{1}{2} \left( u^n_{\left(N_x,j\right)}-u^n_{\left(N_x-1,j\right)} \right), \\
 \Delta^y_{\left(i,1\right)} u^n = \frac{1}{2} \left( u^n_{\left(i,2\right)}-u^n_{\left(i,1\right)} \right), \quad
 \Delta^y_{\left(i,N_y\right)} u^n = \frac{1}{2} \left( u^n_{\left(i,N_y\right)}-u^n_{\left(i,N_y-1\right)} \right).
 \end{gathered}
\end{equation}

\subsubsection{Periodic Boundary Conditions in the Micro Step}
For periodic boundary conditions, in the first step of operator splitting, we do the following:
\begin{align}
\label{eqn:2d_periodic_micro_1}
 \text{Left BC:} & \quad Z^n_{{\left(1,j,k,\ell \right)}} = \left[v_1\right]_{k}^{-} \left( \frac{G^n_{\left(2, j, k, \ell\right)} - G^n_{\left(1, j, k, \ell\right)}}{\Delta x} \right)
+ \left[v_1\right]_{k}^{+} \left( \frac{G^n_{\left(1, j, k, \ell\right)} - G^n_{\left(N_x, j, k, \ell\right)}}{\Delta x} \right), \\
\label{eqn:2d_periodic_micro_2}
 \text{Right BC:} & \quad Z^n_{\left(N_x,j,k,\ell\right)} = \left[v_1\right]_{k}^{-} \left( \frac{G^n_{\left(1, j, k, \ell\right)} - G^n_{\left(N_x, j, k, \ell\right)}}{\Delta x} \right)
+ \left[v_1\right]_{k}^{+} \left( \frac{G^n_{\left(N_x, j, k, \ell\right)} - G^n_{\left(N_x-1, j, k, \ell\right)}}{\Delta x} \right).
\end{align}
In the second step of operator splitting, we set
\begin{align}
\label{eqn:2d_periodic_micro_3}
\text{Bottom BC:} & \quad 
Z^{\star}_{\left(i,1,k,\ell\right)} = \left[v_2\right]^{-}_{\ell} 
\left( \frac{G^{\star}_{\left(i, 2, k, \ell\right)} - G^{\star}_{\left(i,1,k,\ell\right)}}{\Delta y} \right)+ \left[v_2\right]^{+}_{\ell} \left( \frac{G^{\star}_{\left(i, 1, k, \ell\right)} - G^{\star}_{\left(i, N_y, k, \ell\right)}}{\Delta y} \right), \\
\label{eqn:2d_periodic_micro_4}
\text{Top BC:} & \quad
Z^{\star}_{\left(i,N_y,k,\ell\right)} = \left[v_2\right]^{-}_{\ell} 
\left( \frac{G^{\star}_{\left(i, 1, k, \ell\right)} - G^{\star}_{\left(i,N_y,k,\ell\right)}}{\Delta y} \right)+ \left[v_2\right]^{+}_{\ell} \left( \frac{G^{\star}_{\left(i, N_y, k, \ell\right)} - G^{\star}_{\left(i, N_y-1, k, \ell\right)}}{\Delta y} \right).
\end{align}
In the third and final step of operator splitting, we set
\begin{equation}
\label{eqn:2d_periodic_micro_5}
\begin{gathered}
 T^n_{\left(1/2, j\right)}     = T^n_{\left(N_x+1/2, j\right)} = \frac{1}{2} \left( T^n_{\left(1, j\right)} + T^n_{\left(N_x, j\right)} \right), \quad
 T^n_{\left(i, 1/2\right)}     = T^n_{\left(i, N_y+1/2\right)} = \frac{1}{2} \left( T^n_{\left(i, 1\right)}+ T^n_{\left(i, N_y\right)} \right), \\
 \Delta^x_{\left(1,j\right)} u^n = \frac{1}{2} \left( u^n_{\left(2,j\right)}-u^n_{\left(N_x,j\right)} \right), \quad
 \Delta^x_{\left(N_x,j\right)} u^n = \frac{1}{2} \left( u^n_{\left(1,j\right)}-u^n_{\left(N_x-1,j\right)} \right), \\
 \Delta^y_{\left(i,1\right)} u^n = \frac{1}{2} \left( u^n_{\left(i,2\right)}-u^n_{\left(i,N_y\right)} \right), \quad
 \Delta^y_{\left(i,N_y\right)} u^n = \frac{1}{2} \left( u^n_{\left(i,1\right)}-u^n_{\left(i,N_y-1\right)} \right).
 \end{gathered}
\end{equation}


\subsection{2D Macro Finite Volume Method}
In this section, we demonstrate how to update the current value of the macro variables at $t=t^n$
to the new values at $t=t^{n+1} = t^n + \Delta t$:
\begin{equation}
   \text{Input:} \quad \vec{Q^n_{ij}} = \left( \, \rho^n_{ij}, \, \rho \vec{u^n_{ij}}, \, \mat{{\mathbb E}^n_{ij}} \, \right) \quad  
 \text{and} \quad
     \mat{{\mathbb H}^{n+1}_{ij}} \quad \Longrightarrow \quad \text{Output:} \quad \vec{Q^{n+1}_{ij}} = \left( \, \rho^{n+1}_{ij}, \, \rho \vec{u^{n+1}_{ij}}, \, \mat{{\mathbb E}^{n+1}_{ij}} \, \right),
\end{equation}
where the micro portion has already been updated to $t=t^{n+1}$ via \eqref{eqn:update-micro-2D-xdir}, \eqref{eqn:update-micro-2D-ydir}, and \eqref{eqn5-98}, and the heat flux tensor is computed via
\begin{equation}
\label{eqn5-99}
   \begin{pmatrix}[1.5] {\mathbb H}_{111} \\ {\mathbb H}_{112} \\
     {\mathbb H}_{122} \\ {\mathbb H}_{222} \end{pmatrix}^{n+1}_{ij} =
     \varepsilon \, \Delta v_1 \Delta v_2 \sum_{k=1}^{N_{v_1}} \sum_{\ell=1}^{N_{v_2}}
     \begin{pmatrix}[1.5]
      \left(v_{1k} - u^n_{1i} \right)^3 \\ \left( v_{1k} - u^n_{1i} \right)^2 \left(v_{2\ell} - u^n_{2j} \right) \\  
      \left( v_{1k} - u^n_{1i} \right) \left(v_{2\ell} - u^n_{2j} \right)^2 \\ \left(v_{2\ell} - u^n_{2j} \right)^3 
     \end{pmatrix} \, G^{n+1}_{ijk\ell}.
\end{equation}
To maintain the asymptotic-preserving property in the macro update, we elect to solve the macro equation \eqref{eq5-40}--\eqref{eq5-41} using Strang \cite{article:St68} operator splitting such that the collision term will be handled implicitly, while the transport terms will be handled explicitly:
\begin{alignat}{2}
\label{eqn:macro2d-split1}
\frac{\Delta t}{2}:& \quad \begin{bmatrix}
 \pr_{11} \\ \pr_{12} \\ \pr_{22}
\end{bmatrix}_{,t} = \frac{\tau \left(1-\nu\right)}{2\varepsilon}  \begin{bmatrix}
 \pr_{22} - \pr_{11} \\ -2 \pr_{12} \\ \pr_{11} - \pr_{22}
\end{bmatrix} & \qquad & \left( \, \vec{Q^n_{ij}} \rightarrow \vec{Q^{\star}_{ij}} \, \right), \\
\label{eqn:macro2d-split2}
\Delta t:& \quad 
 \vec{q}_{,t} + \vec{f}\left(\vec{q}\right)_{,x}  = \vec{0} 
  & \qquad & \left( \, \vec{Q^{\star}_{ij}} \rightarrow \vec{Q^{\star\star}_{ij}} \, \right), \\
\label{eqn:macro2d-split3}
\Delta t:& \quad 
 \vec{q}_{,t} + \vec{g}\left(\vec{q}\right)_{,y} = \vec{0} 
  & \qquad & \left( \, \vec{Q^{\star\star}_{ij}} \rightarrow \vec{Q^{\star\star\star}_{ij}} \, \right), \\
\label{eqn:macro2d-split4}
\frac{\Delta t}{2}:& \quad \begin{bmatrix}
 \pr_{11} \\ \pr_{12} \\ \pr_{22}
\end{bmatrix}_{,t} = \frac{\tau \left(1-\nu\right)}{2\varepsilon}  \begin{bmatrix}
 \pr_{22} - \pr_{11} \\ -2 \pr_{12} \\ \pr_{11} - \pr_{22}
\end{bmatrix} & \qquad & \left( \, \vec{Q^{\star\star\star}_{ij}} \rightarrow \vec{Q^{n+1}_{ij}} \, \right).
\end{alignat}

The first collision step is discretized using the L-stable TR-BDF2 scheme \cite{article:TRBDF2}, which can be compactly written as follows:
\begin{align}
  \begin{bmatrix}
  \pr_{11} \\
  \pr_{12} \\
  \pr_{22}
  \end{bmatrix}^{\star}_{ij}
  = \frac{1}{2} 
  \begin{bmatrix}
  1 + W^{n}_{ij} & 0 & 1 - W^{n}_{ij} \\
  0 & 2 W^{n}_{ij} & 0 \\
  1 - W^{n}_{ij} & 0 & 1 + W^{n}_{ij}
  \end{bmatrix}
  \begin{bmatrix}
  \pr_{11} \\
  \pr_{12} \\
  \pr_{22}
  \end{bmatrix}^{n}_{ij},
  \,
  W^n_{ij} = \frac{ 48 \varepsilon^2 - 10 \, \tau^n_{ij} \left(1-\nu \right) \varepsilon  \Delta t}{48 \varepsilon^2 + 14 \, \tau^n_{ij} \left(1-\nu \right) \varepsilon  \Delta t + \left( \tau^n_{ij}  \left( 1 - \nu \right) \Delta t \right)^2},
\end{align}
where then
\begin{equation} 
\left( \, \, \rho^{\star}, \quad \rho \vec{u^{\star}}, \quad \mat{\mathbb{E}^{\star}} \right)_{ij} = 
   \left( \, \, \rho^{n}, \quad \rho \vec{u^{n}}, \quad \rho \vec{u^{n}} \otimes \vec{u^{n}} + \mat{\mathbb{P}^{\star}}  \, \, \right)_{ij}.
\end{equation}

In the second and third steps of the Strang operator splitting, the moment equations without the collision operator are solved using the following finite volume method:
 \begin{gather}
\label{eqn5-103}
\begin{split}
\vec{Q^{\star\star}_{ij}} &= \vec{Q^{\star}_{ij}} 
- \frac{\Delta t}{\Delta x} \left( \, \vec{\left[\mathcal{F}_1\right]^{\star}_{{\left(i+{1}/{2}, \, j\right)}}} - \vec{\left[\mathcal{F}_1 \right]^{\star}_{{\left(i-{1}/{2}, \, j\right)}}} \, \right)
- {\Delta t} \begin{pmatrix} 0,  \, 0,  \, 0, \,
\frac{\Delta^x_{ij}\left(\mathbb{H}^{n+1}_{111}\right)}{\Delta x}, \,
\frac{\Delta^x_{ij}\left(\mathbb{H}^{n+1}_{112}\right)}{\Delta x}, \,
\frac{\Delta^x_{ij}\left(\mathbb{H}^{n+1}_{122}\right)}{\Delta x} \end{pmatrix},
\end{split} \\
\begin{split}
\vec{Q^{\star\star\star}_{ij}} &= \vec{Q^{\star\star}_{ij}} 
- \frac{\Delta t}{\Delta y} \left( \, \vec{\left[\mathcal{F}_2\right]^{\star\star}_{{\left(i, \, j+{1}/{2}\right)}}} - \vec{\left[\mathcal{F}_2\right]^{\star\star}_{{\left(i, \, j-{1}/{2}\right)}}} \, \right)
- {\Delta t} \begin{pmatrix} 0, \, 0, \, 0, \,
\frac{\Delta^y_{ij}\left(\mathbb{H}^{n+1}_{112}\right)}{\Delta y}, \,
\frac{\Delta^y_{ij}\left(\mathbb{H}^{n+1}_{122}\right)}{\Delta y}, \,
\frac{\Delta^y_{ij}\left(\mathbb{H}^{n+1}_{222}\right)}{\Delta y} \end{pmatrix},
\end{split}
\end{gather}
where the heat flux values are computed from \eqref{eqn5-99}, 
\begin{equation}
\begin{gathered}
\Delta^x_{ij}\left(\mathbb{H}\right) = {\mathbb{H}_{\left( i+1/2, \, j \right)} - \mathbb{H}_{\left( i-1/2, \, j \right)}}, \quad
\Delta^y_{ij}\left(\mathbb{H}\right) = {\mathbb{H}_{\left( i, \, j+1/2 \right)} - \mathbb{H}_{\left( i, \, j-1/2 \right)}}, \\
\mathbb{H}_{\left( i+1/2, \, j \right)} = \frac{1}{2} \left( \mathbb{H}_{\left(i+1, \, j\right)} + \mathbb{H}_{\left(i, \, j\right)} \right), \quad
\mathbb{H}_{\left( i, \, j+1/2 \right)} = \frac{1}{2} \left( \mathbb{H}_{\left(i, \, j+1\right)} + \mathbb{H}_{\left(i, \, j\right)} \right),
\end{gathered}
\end{equation}
and the numerical fluxes, $\mathcal{F}_1$ and $\mathcal{F}_2$,
 are computed via kinetic flux vector splitting (KFVS) \citep{Mandal1994KineticEquations}, which we briefly describe below.

To compute the fluxes via KFVS, we assume that at each cell interface, the distribution is Gaussian with moments
matching the first six moments in each cell: $\vec{m}\left(v_1, v_2 \right) := \left(1, v_1, v_2, v_2^2, v_1 v_2, v_2^2 \right)$. The left-going (right-going) portion of the flux is based on the fluid variables on the right (left) 
side of the interface, and all the resulting integrals can be computed analytically:
\begin{align}
\label{eqn:flux-2D-1}
\begin{split}
    \vec{\left[{\mathcal{F}_1}\right]^{\star}_{i-{1}/{2} j}} &=  \int_{-\infty}^{\infty}\int_{-\infty}^0 v_1 \, \vec{m}\left(v_1, v_2 \right) \, \gauss^{\star}_{ij} \, dv_1 \, dv_2 + \int_{-\infty}^{\infty}\int_0^{\infty}v_1 \, \vec{m}\left(v_1, v_2 \right) \, \gauss^{\star}_{i-1 j} \, dv_1 \, dv_2 \\
    &= \frac{1}{2} \left( \left[ a_1 \right]^{\star}_{ij} \left[ J_1 \right]^{\star}_{ij}
        + \left( 1 + \left[ b_1 \right]^{\star}_{ij} \right) \left[ K_1 \right]^{\star}_{ij}
        - \left[ a_1 \right]^{\star}_{i-1 j} \left[ J_1 \right]^{\star}_{i-1 j}
        + \left( 1 - \left[ b_1 \right]^{\star}_{i-1 j} \right) \left[ K_1 \right]^{\star}_{i-1 j}
         \right), 
    \end{split} \\
\label{eqn:flux-2D-2}
\begin{split}
\vec{\left[{\mathcal{F}_2}\right]^{\star\star}_{i j-{1}/{2}}} &=  \int_{-\infty}^{\infty}\int_{-\infty}^0 v_2 \, \vec{m}\left(v_1, v_2 \right) \, \gauss^{\star\star}_{ij} \, dv_2 \, dv_1 + \int_{-\infty}^{\infty}\int_0^{\infty}v_2 \, \vec{m}\left(v_1, v_2 \right) \, \gauss^{\star\star}_{i j-1} \, dv_2 \, dv_1 \\
    &= \frac{1}{2} \left( \left[ a_2 \right]^{\star\star}_{ij} \left[ J_2 \right]^{\star\star}_{ij}
        + \left( 1 + \left[ b_2 \right]^{\star\star}_{ij} \right) \left[ K_2 \right]^{\star\star}_{ij}
        - \left[ a_2 \right]^{\star\star}_{i j-1} \left[ J_2 \right]^{\star\star}_{i j-1}
        + \left( 1 - \left[ b_2 \right]^{\star\star}_{i j-1} \right) \left[ K_2 \right]^{\star\star}_{i j-1}
         \right), 
\end{split}
\end{align}
where
\begin{gather}
\begin{split}
\begin{matrix}
a_1 := \sqrt{\frac{2\pr_{11}}{\pi \rho}}\exp\left(-\frac{\rho u_1^2}{2\pr_{11}}\right) \\
b_1 := \mathrm{erf}\left(u_1\sqrt{\frac{\rho}{2\pr_{11}}}\right)
\end{matrix}, \quad
J_1 := \begin{bmatrix}
    \rho \\
    \rho u_1\\
    \rho u_2\\
    \rho u_1^2 + 2 \pr_{11}\\
    \rho u_1 u_2 + 2 \pr_{12}\\
    \rho u_2^2 + \pr_{22} + \frac{\pr^2_{12}}{\pr_{11}}\\
\end{bmatrix}, \quad 
K_1 = \begin{bmatrix}
    \rho u_1\\
    \rho u_1^2 + \pr_{11}\\
    \rho u_1 u_2 + \pr_{12}\\
    \rho u_1^3 + 3 u_1 \pr_{11}\\
    \rho u_1^2 u_2 + u_2 \pr_{11} + 2 u_1 \pr_{12}\\
    \rho u_1 u_2^2 + u_1 \pr_{22} + 2 u_2 \pr_{12}
\end{bmatrix}, \\
\begin{matrix}
a_2 = \sqrt{\frac{2 \pr_{22}}{\pi \rho}}\exp\left(-\frac{\rho u_2^2}{2\pr_{22}}\right) \\
b_2 = \mathrm{erf}\left(u_2\sqrt{\frac{\rho}{2\pr_{22}}}\right)
\end{matrix}, \quad
J_2 = \begin{bmatrix}
    \rho \\
    \rho u_1\\
    \rho u_2\\
    \rho u_1^2 + \pr_{11} + \frac{\pr^2_{12}}{\pr_{22}}\\
    \rho u_1 u_2 + 2 \pr_{12}\\
    \rho u_2^2 + 2 \pr_{22}\\
\end{bmatrix}, \quad 
K_2 = \begin{bmatrix}
    \rho u_2\\
    \rho u_1 u_2 + \pr_{12}\\
    \rho u_2^2 + \pr_{22}\\
    \rho u_1^2 u_2 + u_2 \pr_{11} + 2 u_1 \pr_{12}\\
    \rho u_1 u_2^2 + u_1 \pr_{22} + 2 u_2 \pr_{12}\\
    \rho u_2^3 + 3 u_2 \pr_{22}
\end{bmatrix}.
\end{split}
\end{gather}

The fourth and final step in the Strang splitting is a second collision step, which can be compactly written as follows:
\begin{gather}
  \begin{bmatrix}
  \pr_{11} \\
  \pr_{12} \\
  \pr_{22}
  \end{bmatrix}^{n+1}_{ij}
  = \frac{1}{2} 
  \begin{bmatrix}
  1 + W^{\star\star\star}_{ij} & 0 & 1 - W^{\star\star\star}_{ij} \\
  0 & 2 W^{\star\star\star}_{ij} & 0 \\
  1 - W^{\star\star\star}_{ij} & 0 & 1 + W^{\star\star\star}_{ij}
  \end{bmatrix}
  \begin{bmatrix}
  \pr_{11} \\
  \pr_{12} \\
  \pr_{22}
  \end{bmatrix}^{\star\star\star}_{ij},
  \\
  W^{\star\star\star}_{ij} = \frac{ 48 \varepsilon^2 - 10 \, \tau^{\star\star\star}_{ij} \left(1-\nu \right) \varepsilon  \Delta t}{48 \varepsilon^2 + 14 \, \tau^{\star\star\star}_{ij} \left(1-\nu \right) \varepsilon  \Delta t + \left( \tau^{\star\star\star}_{ij}  \left( 1 - \nu \right) \Delta t\right)^2},
\end{gather}
where then
\begin{equation} 
\left( \, \, \rho^{n+1}, \quad \rho \vec{u^{n+1}}, \quad \mat{\mathbb{E}^{n+1}} \right)_{ij} = 
   \left(\rho^{\star\star\star}, \quad \rho \vec{u^{\star\star\star}}, \quad \rho \vec{u^{\star\star\star}} \otimes \vec{u^{\star\star\star}} + \mat{\mathbb{P}^{n+1}} \, \, \right)_{ij}.
\end{equation}

\subsubsection{Extrapolation Boundary Conditions in the Macro Step}
We briefly mention here how simple boundary conditions are implemented in the 2D macro update. For extrapolation boundary conditions, we do the following:
\begin{align}
\label{eqn:2d_extrap_macro_1}
 \text{Left BC:} & \quad \mathbb{H}_{\left( 1/2, \, j \right)} =  \mathbb{H}_{\left(1, \, j\right)}, \quad \vec{\left[\mathcal{F}_1\right]^{\star}_{{\left({1}/{2}, \, j\right)}}} = \vec{f}\left( \vec{Q^{\star}_{\left(1, \, j\right)}} \right), \\
 \label{eqn:2d_extrap_macro_2}
 \text{Right BC:} & \quad \mathbb{H}_{\left( N_x+1/2, \, j \right)} =  \mathbb{H}_{\left(N_x, \, j\right)}, \quad \vec{\left[\mathcal{F}_1\right]^{\star}_{{\left(N_x+{1}/{2}, \, j\right)}}} = \vec{f}\left( \vec{Q^{\star}_{\left(N_x, \, j\right)}} \right), \\
 \label{eqn:2d_extrap_macro_3}
 \text{Bottom BC:} & \quad \mathbb{H}_{\left( i, \, 1/2 \right)} =  \mathbb{H}_{\left(i, \, 1/2\right)}, \quad \vec{\left[\mathcal{F}_2\right]^{\star\star}_{{\left(i, \, {1}/{2}\right)}}} = \vec{g}\left( \vec{Q^{\star\star}_{\left(i, \, 1\right)}} \right), \\
 \label{eqn:2d_extrap_macro_4}
 \text{Top BC:} & \quad \mathbb{H}_{\left( i, \, N_y+1/2 \right)} =  \mathbb{H}_{\left(i, \, N_y\right)}, \quad \vec{\left[\mathcal{F}_2\right]^{\star\star}_{{\left(i, \, N_y+{1}/{2} \right)}}} = \vec{g}\left( \vec{Q^{\star\star}_{\left(i, \, N_y\right)}} \right).
\end{align}

\subsubsection{Periodic Boundary Conditions in the Macro Step}
For periodic boundary conditions, we do the following:
\begin{align}
\label{eqn:2d_periodic_macro_1}
\begin{split}
\text{Left/Right BCs:} & \quad
 \mathbb{H}_{\left( 1/2, \, j \right)} = \mathbb{H}_{\left( N_x + 1/2, \, j \right)} =  
    \frac{1}{2} \left( \mathbb{H}_{\left(1, \, j \right)} + \mathbb{H}_{\left(N_x, \, j \right)} \right), \\ 
   \vec{\left[\mathcal{F}_1\right]^{\star}_{{\left({1}/{2}, \, j\right)}}} &= \vec{\left[\mathcal{F}_1\right]^{\star}_{{\left(N_x+{1}/{2}, \, j\right)}}} = 
 \int_{-\infty}^{\infty}\int_{-\infty}^0 v_1 \, \vec{m} \, \gauss^{\star}_{\left(1, \, j\right)} \, dv_1 \, dv_2 + \int_{-\infty}^{\infty}\int_0^{\infty}v_1 \, \vec{m} \, \gauss^{\star}_{\left(N_x, \, j\right)} \, dv_1 \, dv_2,
     \end{split} \\
\label{eqn:2d_periodic_macro_2}
\begin{split}
\text{Bottom/Top BCs:} & \quad
 \mathbb{H}_{\left( i, \, 1/2 \right)} = \mathbb{H}_{\left( i, \, N_y + 1/2 \right)} =  
    \frac{1}{2} \left( \mathbb{H}_{\left(i, \, 1 \right)} + \mathbb{H}_{\left(i, \, N_y \right)} \right), \\ 
   \vec{\left[\mathcal{F}_2\right]^{\star\star}_{{\left(i, \, {1}/{2} \right)}}} &= \vec{\left[\mathcal{F}_2\right]^{\star\star}_{{\left(i, \, N_y+{1}/{2} \right)}}} = 
 \int_{-\infty}^{\infty}\int_{-\infty}^0 v_2 \, \vec{m} \, \gauss^{\star\star}_{\left(i, \, 1\right)} \, dv_2 \, dv_1 + \int_{-\infty}^{\infty}\int_0^{\infty}v_2 \, \vec{m} \, \gauss^{\star\star}_{\left(i, \, N_y \right)} \, dv_2 \, dv_1.
     \end{split}
\end{align}

\subsection{Message Passing Interface (MPI) Implementation}
\label{subsec:MPI}
We briefly describe here how the 2D code is parallelized using MPI (Message Passing Interface) \cite{mpi50}.  We assume a total mesh of size 
\begin{equation}
N_x \times N_y \times N_{v_1} \times N_{v_2} \quad \text{over} \quad
\left(x_{\text{min}}, x_{\text{max}} \right) \times \left(y_{\text{min}}, y_{\text{max}} \right) \times \left(\left[v_{1}\right]_{\text{min}}, \left[v_{1}\right]_{\text{max}} \right) \times \left(\left[v_{2}\right]_{\text{min}}, \left[v_{2}\right]_{\text{max}} \right),
\end{equation}
 and access to $N_{\text{Proc}}$ processors, where $N_{\text{Proc}}$ is a perfect square. The proposed parallelization divides the spatial domain into $N_{\text{Proc}}$ equally-sized sub-domains, so that each processor has a mesh of size 
\begin{equation}
\frac{N_x}{\sqrt{N_{\text{Proc}}}} \times \frac{N_y}{\sqrt{N_{\text{Proc}}}} \times N_{v_1} \times N_{v_2} \quad \text{over} \quad
\left(x^{(p)}_{\text{min}}, x^{(p)}_{\text{max}} \right) \times \left(y^{(p)}_{\text{min}}, y^{(p)}_{\text{max}} \right) \times \left(\left[v_{1}\right]_{\text{min}}, \left[v_{1}\right]_{\text{max}} \right) \times \left(\left[v_{2}\right]_{\text{min}}, \left[v_{2}\right]_{\text{max}} \right),
\end{equation}
 where
 \begin{alignat}{2}
 x_{\text{min}}^{(p)} &= x_{\text{min}} + {p_{\text{m}}} \cdot \left(\frac{x_{\text{max}} - x_{\text{min}}}{\sqrt{N_{\text{proc}}}}\right), \quad
 &x_{\text{max}}^{(p)} &= x_{\text{min}}^{(p)} + \left(\frac{x_{\text{max}} - x_{\text{min}}}{\sqrt{N_{\text{proc}}}}\right),
 \\
 y_{\text{min}}^{(p)} &= y_{\text{min}} + \left( p - p_{\text{m}} \right) \cdot \left(\frac{y_{\text{max}} - y_{\text{min}}}{{N_{\text{proc}}}}\right), \quad
 &y_{\text{max}}^{(p)} &= y_{\text{min}}^{(p)} + \left(\frac{y_{\text{max}} - y_{\text{min}}}{\sqrt{N_{\text{proc}}}}\right),
 \end{alignat}
 $p_{\text{m}} = p\left(\text{mod} \, \sqrt{N_{\text{proc}}}\right)$, and $p=0,1,2,\ldots,\left(N_{\text{proc}}-1\right)$ is the processor index. We illustrate a parallel setup as described in this section with $16$ processors in Figure \ref{fig:mpi}.

For each processor, $p$, we need to know what other processors share an internal boundary. As a first pass, we can use the following formulas to find what processors are to the east, west, north, and south of the current processor $p$:
\begin{align}
  \text{east}^{(p)} = p + 1, \quad
  \text{west}^{(p)} = p - 1, \quad
  \text{north}^{(p)} = p + \sqrt{N_{\text{proc}}}, \quad
  \text{south}^{(p)} = p - \sqrt{N_{\text{proc}}}.
\end{align}
However, if processor $p$ borders the actual physical boundary on any of its four sides, we overwrite the first pass by the following:
\begin{align}
  \text{if} \quad p < \sqrt{N_{\text{proc}}}:& \quad
    \text{south}^{(p)} = \begin{cases}
    {N_{\text{proc}}}-\sqrt{N_{\text{proc}}} + p,  &\text{(periodic BCs)}\\
    -1,  &\text{(other BCs)}
    \end{cases} \\
  \text{if} \quad p > N_{\text{proc}} - \sqrt{N_{\text{proc}}}-1:& \quad 
  \text{north}^{(p)} = \begin{cases}
  p\left(\text{mod} \, \sqrt{N_{\text{proc}}}\right), &\text{(periodic BCs)}\\
  -1, &\text{(other BCs)}
  \end{cases} \\
  \text{if} \quad p\left(\text{mod} \, \sqrt{N_{\text{proc}}}\right)=0:& \quad 
  \text{west}^{(p)} = \begin{cases}
  p+\sqrt{N_{\text{proc}}}-1,  &\text{(periodic BCs)}\\
  -1, & \text{(other BCs)}
  \end{cases} \\
  \text{if} \quad \left(p+1\right)\left(\text{mod} \, \sqrt{N_{\text{proc}}}\right)=0:& \quad   
  \text{east}^{(p)} = \begin{cases}
  p + 1 - \sqrt{N_{\text{proc}}}, &\text{(periodic BCs)}\\
  -1. &\text{(other BCs)}
  \end{cases}
\end{align}
Note that periodic boundary conditions are handled via domain wraparound and all other boundary conditions by setting the directional pointer to $-1$.

Now that the mesh is divided among the $N_{\text{proc}}$ processors, we can execute the code independently on each processor with the caveat that we will need boundary condition data from neighboring processors. To explain how this is done, we use the following notation for variables on processor $p$:
\begin{align}
G^n_{\left(i, j, k, \ell\right)} := \text{solution on proc $p$, at time $t^n$, and on grid cell $(i,j,k,\ell)$ as indexed on proc $p$},
\end{align}
and the following notation for variables on neighboring processors:
\begin{align}
\inred{$\text{west}^{(p)}$: $G^n_{\left(i, j, k, \ell\right)}$} :=& \, \, \text{solution on proc $\text{west}^{(p)}$ as indexed on proc $\text{west}^{(p)}$}, \\
\inred{$\text{east}^{(p)}$: $G^n_{\left(i, j, k, \ell\right)}$} :=& \, \, \text{solution on proc $\text{east}^{(p)}$ as indexed on proc $\text{east}^{(p)}$}, \\
\inred{$\text{south}^{(p)}$: $G^n_{\left(i, j, k, \ell\right)}$} :=& \, \, \text{solution on proc $\text{south}^{(p)}$ as indexed on proc $\text{south}^{(p)}$}, \\
\inred{$\text{north}^{(p)}$: $G^n_{\left(i, j, k, \ell\right)}$} :=& \, \, \text{solution on proc $\text{north}^{(p)}$ as indexed on proc $\text{north}^{(p)}$}.
\end{align}

\subsubsection{MPI Micro Update}
In the first part of the micro update, the following data must be sent to the neighboring processors (except when the directional pointer is $-1$):
\begin{align}
{\tt MPI\_Send}: \quad G^n_{\left(1, j, k, \ell\right)} \, \text{to proc $\text{west}^{(p)}$} \quad \text{and} \quad G^n_{\left(N_x, j, k, \ell\right)} \, \text{to proc $\text{east}^{(p)}$}, \quad \forall j,k,\ell,
\end{align}
and the following information must either be received from other processors, or if the directional pointer is $-1$, an appropriate boundary value must be supplied (e.g., extrapolation or wall conditions):
\begin{align}
{\tt MPI\_Recv}: \quad \inred{$\text{west}^{(p)}$: $G^n_{\left(N_x, j, k, \ell\right)}$} \quad \text{and} \quad 
\inred{$\text{east}^{(p)}$: $G^n_{\left(1, j, k, \ell\right)}$}, \quad \forall j,k,\ell.
\end{align}
This received data is then used in the numerical updates of the cells that are bordering other processors:
\begin{align}
{Z^n_{\left(1,j,k,\ell\right)}} &= {\left[v_1\right]_{k}^{-} \left( \frac{G^n_{\left(2, j, k, \ell\right)} - G^n_{\left(1, j, k, \ell\right)}}{\Delta x} \right)}
+ {\left[v_1\right]_{k}^{+} \left( \frac{G^n_{\left(1, j, k, \ell\right)} - {\inred{$\text{west}^{(p)}$: $G^n_{\left(N_x, j, k, \ell\right)}$}}}{\Delta x} \right)}, \\
{Z^n_{{\left(N_x,j,k,\ell\right)}}} &= {\left[v_1\right]_{k}^{-} \left( \frac{\inred{$\text{east}^{(p)}$: $G^n_{\left(1, j, k, \ell\right)}$} - \, G^n_{\left(N_x, j, k, \ell\right)}}{\Delta x} \right)}
+ {\left[v_1\right]_{k}^{+} \left( \frac{G^n_{\left(N_x, j, k, \ell\right)} - G^n_{\left(N_x-1, j, k, \ell\right)}}{\Delta x} \right)}.
\end{align}

In the second part of the micro update, the following data must be sent to the neighboring processors (except when the directional pointer is $-1$):
\begin{align}
{\tt MPI\_Send}: \quad G^n_{\left(i, 1, k, \ell\right)} \, \text{to proc $\text{south}^{(p)}$} \quad \text{and} \quad G^n_{\left(i, N_y, j, k, \ell\right)} \, \text{to proc $\text{north}^{(p)}$}, \quad \forall i,k,\ell,
\end{align}
and the following information must either be received from other processors, or if the directional pointer is $-1$, an appropriate boundary value must be supplied (e.g., extrapolation or wall conditions):
\begin{align}
{\tt MPI\_Recv}: \quad \inred{$\text{south}^{(p)}$: $G^n_{\left(i, N_y, k, \ell\right)}$} \quad \text{and} \quad 
\inred{$\text{north}^{(p)}$: $G^n_{\left(i, 1, k, \ell\right)}$}, \quad \forall i,k,\ell.
\end{align}
This received data is then used in the numerical updates of the cells that are bordering other processors:
\begin{align}
Z^{\star}_{\left(i,1,k,\ell\right)} &= \left[v_2\right]^{-}_{\ell} 
\left( \frac{G^{\star}_{\left(i, 2, k, \ell\right)} - G^{\star}_{\left(i,1,k,\ell\right)}}{\Delta y} \right)+ \left[v_2\right]^{+}_{\ell} \left( \frac{G^{\star}_{\left(i, 1, k, \ell\right)} - \inred{$\text{south}^{(p)}$: $G^{\star}_{\left(i, N_y, k, \ell\right)}$}}{\Delta y} \right), \\
Z^{\star}_{\left(i,N_y,k,\ell\right)} &= \left[v_2\right]^{-}_{\ell} 
\left( \frac{\inred{$\text{north}^{(p)}$: $G^{\star}_{\left(i, 1, k, \ell\right)}$} - G^{\star}_{\left(i,N_y,k,\ell\right)}}{\Delta y} \right)+ \left[v_2\right]^{+}_{\ell} \left( \frac{G^{\star}_{\left(i, N_y, k, \ell\right)} - G^{\star}_{\left(i, N_y-1, k, \ell\right)}}{\Delta y} \right).
\end{align}

Finally, in the third part of the micro update, the following data must be sent to the neighboring processors (except when the directional pointer is $-1$):
\begin{align}
{\tt MPI\_Send}:& \quad T^n_{\left(1, j\right)} \quad \text{and} \quad u^n_{\left(1, j\right)} \quad \text{to proc $\text{west}^{(p)}$}, \quad T^n_{\left(N_x, j\right)} \quad \text{and} \quad u^n_{\left(N_x, j\right)} \quad \text{to proc $\text{east}^{(p)}$}, \quad \forall j, \\
{\tt MPI\_Send}:& \quad T^n_{\left(i, 1\right)} \quad \text{and} \quad u^n_{\left(i, 1\right)} \quad \text{to proc $\text{south}^{(p)}$}, \quad T^n_{\left(i, N_y \right)} \quad \text{and} \quad u^n_{\left(i, N_y\right)} \quad \text{to proc $\text{north}^{(p)}$}, \quad \forall i,
\end{align}
and the following information must either be received from other processors, or if the directional pointer is $-1$, an appropriate boundary value must be supplied (e.g., extrapolation or wall conditions):
\begin{align}
{\tt MPI\_Recv}:& \quad \inred{$\text{west}^{(p)}$: $T^n_{\left(N_x, j\right)}$}, \quad
\inred{$\text{west}^{(p)}$: $u^n_{\left(N_x, j\right)}$}, \quad
\inred{$\text{east}^{(p)}$: $T^n_{\left(1, j\right)}$}, \quad
\inred{$\text{east}^{(p)}$: $u^n_{\left(1, j\right)}$}, \quad \forall j, \\
{\tt MPI\_Recv}:& \quad \inred{$\text{south}^{(p)}$: $T^n_{\left(i,N_y\right)}$}, \quad
\inred{$\text{south}^{(p)}$: $u^n_{\left(i,N_y\right)}$}, \quad
\inred{$\text{north}^{(p)}$: $T^n_{\left(i,1\right)}$}, \quad
\inred{$\text{north}^{(p)}$: $u^n_{\left(i,1\right)}$}, \quad \forall i.
\end{align}
This received data is then used to define edge temperature and velocity values:
\begin{gather}
 T^n_{\left(1/2, j\right)} = \frac{1}{2} \left( T^n_{\left(1, j\right)} + \inred{$\text{west}^{(p)}$: $T^n_{\left(N_x, j\right)}$} \right), \quad
 T^n_{\left(N_x+1/2, j\right)} = \frac{1}{2} \left( \inred{$\text{east}^{(p)}$: $T^n_{\left(1, j\right)}$} + T^n_{\left(N_x, j\right)} \right), \\
 T^n_{\left(i, 1/2\right)} = \frac{1}{2} \left( T^n_{\left(i, 1\right)}+ \inred{$\text{south}^{(p)}$: $T^n_{\left(i, N_y\right)}$} \right), \quad
T^n_{\left(i, N_y+1/2\right)} = \frac{1}{2} \left( \inred{$\text{north}^{(p)}$: $T^n_{\left(i, 1\right)}$} + T^n_{\left(i, N_y\right)} \right), \\
 \Delta^x_{\left(1,j\right)} u^n = \frac{1}{2} \left( u^n_{\left(2,j\right)}-\inred{$\text{west}^{(p)}$: $u^n_{\left(N_x,j\right)}$} \right), \quad
 \Delta^x_{\left(N_x,j\right)} u^n = \frac{1}{2} \left( \inred{$\text{east}^{(p)}$: $u^n_{\left(1,j\right)}$} - u^n_{\left(N_x-1,j\right)} \right), \\
 \Delta^y_{\left(i,1\right)} u^n = \frac{1}{2} \left( u^n_{\left(i,2\right)}-\inred{$\text{south}^{(p)}$: $u^n_{\left(i,N_y\right)}$} \right), \quad
 \Delta^y_{\left(i,N_y\right)} u^n = \frac{1}{2} \left( \inred{$\text{north}^{(p)}$: $u^n_{\left(i,1\right)}$} - u^n_{\left(i,N_y-1\right)} \right),
\end{gather}
which are then used in \eqref{eqn5-97} and \eqref{eqn5-97more}.

\subsubsection{MPI Macro Update}
In the second part of the macro update (see \eqref{eqn:macro2d-split2}), the following data must be sent to the neighboring processors (except when the directional pointer is $-1$):
\begin{align}
{\tt MPI\_Send}:& \quad \mathbb{H}_{\left(1, \, j \right)} \quad \text{and} \quad \vec{\left[\mathcal{F}_1\right]^{-}_{{\left({1}/{2}, \, j\right)}}} = \int_{-\infty}^{\infty}\int_{-\infty}^{0} v_1 \, \vec{m} \, \gauss^{\star}_{\left(1, \, j\right)} \, dv_1 \, dv_2
	 \quad \text{to proc $\text{west}^{(p)}$}, \quad \forall j, \\
{\tt MPI\_Send}:& \quad  \mathbb{H}_{\left(N_x, \, j \right)} \quad \text{and} \quad \vec{\left[\mathcal{F}_1\right]^{+}_{{\left(N_x+{1}/{2}, \, j\right)}}} = \int_{-\infty}^{\infty}\int_{0}^{\infty} v_1 \, \vec{m} \, \gauss^{\star}_{\left(N_x, \, j\right)} \, dv_1 \, dv_2 \quad \text{to proc $\text{east}^{(p)}$}, \quad \forall j,
\end{align}
and the following information must either be received from other processors, or if the directional pointer is $-1$, an appropriate boundary value must be supplied (e.g., extrapolation or wall conditions):
\begin{align}
{\tt MPI\_Recv}: \quad \inred{$\text{west}^{(p)}$: $\mathbb{H}_{\left(N_x, j\right)}$}, \quad 
\inred{$\text{west}^{(p)}$: $\vec{\left[\mathcal{F}_1\right]^{+}_{{\left(N_x+{1}/{2}, \, j\right)}}}$},  
\quad \inred{$\text{east}^{(p)}$: $\mathbb{H}_{\left(1, j\right)}$}, \quad 
\inred{$\text{east}^{(p)}$: $\vec{\left[\mathcal{F}_1\right]^{-}_{{\left({1}/{2}, \, j\right)}}}$}, \quad \forall j.
\end{align}
This received data is then used in the numerical updates of the cells that are bordering other processors:
\begin{align}
 \mathbb{H}_{\left( 1/2, \, j \right)} &=  
    \frac{1}{2} \left( \mathbb{H}_{\left(1, \, j \right)} + \inred{$\text{west}^{(p)}$: $\mathbb{H}_{\left(N_x, \, j \right)}$} \right), 
    \quad
  \mathbb{H}_{\left( N_x + 1/2, \, j \right)} =  
    \frac{1}{2} \left( \inred{$\text{east}^{(p)}$: $\mathbb{H}_{\left(1, \, j \right)}$} + \mathbb{H}_{\left(N_x, \, j \right)} \right), \\
   \vec{\left[\mathcal{F}_1\right]^{\star}_{{\left({1}/{2}, \, j\right)}}} &= 
 \int_{-\infty}^{\infty}\int_{-\infty}^0 v_1 \, \vec{m} \, \gauss^{\star}_{\left(1, \, j\right)} \, dv_1 \, dv_2 + \int_{-\infty}^{\infty}\int_0^{\infty}v_1 \, \vec{m} \, \inred{$\text{west}^{(p)}$: $\gauss^{\star}_{\left(N_x, \, j\right)}$} \, dv_1 \, dv_2, \\
 \vec{\left[\mathcal{F}_1\right]^{\star}_{{\left(N_x+{1}/{2}, \, j\right)}}} &= 
 \int_{-\infty}^{\infty}\int_{-\infty}^0 v_1 \, \vec{m} \, \inred{$\text{east}^{(p)}$: $\gauss^{\star}_{\left(1, \, j\right)}$} \, dv_1 \, dv_2 + \int_{-\infty}^{\infty}\int_0^{\infty}v_1 \, \vec{m} \, \gauss^{\star}_{\left(N_x, \, j\right)} \, dv_1 \, dv_2.
 \end{align}

In the third part of the macro update (see \eqref{eqn:macro2d-split3}), the following data must be sent to the neighboring processors (except when the directional pointer is $-1$):
\begin{align}
{\tt MPI\_Send}:& \quad \mathbb{H}_{\left(i, 1 \right)} \quad \text{and} \quad \vec{\left[\mathcal{F}_2\right]^{-}_{{ \left(i, {1}/{2} \right)}}} = \int_{-\infty}^{\infty}\int_{-\infty}^{0} v_2 \, \vec{m} \, \gauss^{\star\star}_{\left(i, 1 \right)} \, dv_2 \, dv_1
	 \quad \text{to proc $\text{south}^{(p)}$}, \quad \forall i, \\
{\tt MPI\_Send}:& \quad  \mathbb{H}_{\left(i, N_y \right)} \quad \text{and} \quad \vec{\left[\mathcal{F}_2\right]^{+}_{{\left(i, N_y+{1}/{2} \right)}}} = \int_{-\infty}^{\infty}\int_{0}^{\infty} v_2 \, \vec{m} \, \gauss^{\star\star}_{\left(i, N_y \right)} \, dv_2 \, dv_1 \quad \text{to proc $\text{north}^{(p)}$}, \quad \forall i,
\end{align}
and the following information must either be received from other processors, or if the directional pointer is $-1$, an appropriate boundary value must be supplied (e.g., extrapolation or wall conditions):
\begin{align}
{\tt MPI\_Recv}: \, \, \inred{$\text{south}^{(p)}$: $\mathbb{H}_{\left(i,N_y\right)}$}, \quad 
\inred{$\text{south}^{(p)}$: $\vec{\left[\mathcal{F}_2\right]^{+}_{{\left(i, N_y+{1}/{2} \right)}}}$},  
\quad \inred{$\text{north}^{(p)}$: $\mathbb{H}_{\left(i, 1\right)}$}, \quad 
\inred{$\text{north}^{(p)}$: $\vec{\left[\mathcal{F}_2\right]^{-}_{{ \left(i, {1}/{2} \right)}}}$}, \quad \forall i.
\end{align}
The received data is then used in the numerical updates of the cells that border other processors:
\begin{align}
 \mathbb{H}_{\left( i, \, 1/2 \right)} &= 
    \frac{1}{2} \left( \mathbb{H}_{\left(i, \, 1 \right)} + \inred{$\text{south}^{(p)}$: $\mathbb{H}_{\left(i, \, N_y \right)}$} \right), 
    \quad
    \mathbb{H}_{\left( i, \, N_y + 1/2 \right)} =  
    \frac{1}{2} \left( \inred{$\text{north}^{(p)}$: $\mathbb{H}_{\left(i, \, 1 \right)}$} + \mathbb{H}_{\left(i, \, N_y \right)} \right), \\ 
   \vec{\left[\mathcal{F}_2\right]^{\star\star}_{{\left(i, \, {1}/{2} \right)}}} &= 
 \int_{-\infty}^{\infty}\int_{-\infty}^0 v_2 \, \vec{m} \, \gauss^{\star\star}_{\left(i, \, 1\right)} \, dv_2 \, dv_1 + \int_{-\infty}^{\infty}\int_0^{\infty}v_2 \, \vec{m} \, \inred{$\text{south}^{(p)}$: $\gauss^{\star\star}_{\left(i, \, N_y \right)}$} \, dv_2 \, dv_1, \\
 \vec{\left[\mathcal{F}_2\right]^{\star\star}_{{\left(i, \, N_y+{1}/{2} \right)}}} &= 
 \int_{-\infty}^{\infty}\int_{-\infty}^0 v_2 \, \vec{m} \, \inred{$\text{north}^{(p)}$: $\gauss^{\star\star}_{\left(i, \, 1\right)}$} \, dv_2 \, dv_1 + \int_{-\infty}^{\infty}\int_0^{\infty}v_2 \, \vec{m} \, \gauss^{\star\star}_{\left(i, \, N_y \right)} \, dv_2 \, dv_1.
\end{align}

\begin{figure}[!h]
\begin{center}
\begin{tikzpicture}[thick,scale=1.5]
\draw (0,0) grid (4,4);
\draw[<->, red, thick] (0.875,0.5) -- (1.125,0.5);
\draw[<->, red, thick] (1.875,0.5) -- (2.125,0.5);
\draw[<->, red, thick] (2.875,0.5) -- (3.125,0.5);
\draw[<->, red, thick] (0.875,1.5) -- (1.125,1.5);
\draw[<->, red, thick] (1.875,1.5) -- (2.125,1.5);
\draw[<->, red, thick] (2.875,1.5) -- (3.125,1.5);
\draw[<->, red, thick] (0.875,2.5) -- (1.125,2.5);
\draw[<->, red, thick] (1.875,2.5) -- (2.125,2.5);
\draw[<->, red, thick] (2.875,2.5) -- (3.125,2.5);
\draw[<->, red, thick] (0.875,3.5) -- (1.125,3.5);
\draw[<->, red, thick] (1.875,3.5) -- (2.125,3.5);
\draw[<->, red, thick] (2.875,3.5) -- (3.125,3.5);
\draw[<->, red, thick] (0.5,0.875) -- (0.5,1.125);
\draw[<->, red, thick] (0.5,1.875) -- (0.5,2.125);
\draw[<->, red, thick] (0.5,2.875) -- (0.5,3.125);
\draw[<->, red, thick] (1.5,0.875) -- (1.5,1.125);
\draw[<->, red, thick] (1.5,1.875) -- (1.5,2.125);
\draw[<->, red, thick] (1.5,2.875) -- (1.5,3.125);
\draw[<->, red, thick] (2.5,0.875) -- (2.5,1.125);
\draw[<->, red, thick] (2.5,1.875) -- (2.5,2.125);
\draw[<->, red, thick] (2.5,2.875) -- (2.5,3.125);
\draw[<->, red, thick] (3.5,0.875) -- (3.5,1.125);
\draw[<->, red, thick] (3.5,1.875) -- (3.5,2.125);
\draw[<->, red, thick] (3.5,2.875) -- (3.5,3.125);
\filldraw [blue] (0.0,0.5) circle (1.5pt);
\filldraw [blue] (0.0,1.5) circle (1.5pt);
\filldraw [blue] (0.0,2.5) circle (1.5pt);
\filldraw [blue] (0.0,3.5) circle (1.5pt);
\filldraw [blue] (4.0,0.5) circle (1.5pt);
\filldraw [blue] (4.0,1.5) circle (1.5pt);
\filldraw [blue] (4.0,2.5) circle (1.5pt);
\filldraw [blue] (4.0,3.5) circle (1.5pt);
\filldraw [blue] (0.5,0.0) circle (1.5pt);
\filldraw [blue] (1.5,0.0) circle (1.5pt);
\filldraw [blue] (2.5,0.0) circle (1.5pt);
\filldraw [blue] (3.5,0.0) circle (1.5pt);
\filldraw [blue] (0.5,4.0) circle (1.5pt);
\filldraw [blue] (1.5,4.0) circle (1.5pt);
\filldraw [blue] (2.5,4.0) circle (1.5pt);
\filldraw [blue] (3.5,4.0) circle (1.5pt);
\node at (0.5,0.5) {0};
\node at (1.5,0.5) {1};
\node at (2.5,0.5) {2};
\node at (3.5,0.5) {3};
\node at (0.5,1.5) {4};
\node at (1.5,1.5) {5};
\node at (2.5,1.5) {6};
\node at (3.5,1.5) {7};
\node at (0.5,2.5) {8};
\node at (1.5,2.5) {9};
\node at (2.5,2.5) {10};
\node at (3.5,2.5) {11};
\node at (0.5,3.5) {12};
\node at (1.5,3.5) {13};
\node at (2.5,3.5) {14};
\node at (3.5,3.5) {15};
  \draw[->, very thick] (-0.35,-0.35) -- (0.65,-0.35) node[below] {$x$};
  \draw[->, very thick] (-0.35,-0.35) -- (-0.35,0.65) node[left] {$y$};
  \filldraw [black] (-0.35,-0.35) circle (1pt);
\end{tikzpicture}
\caption{(\S\ref{subsec:MPI} MPI implementation): Shown here is a parallel setup with $N_{\text{Proc}}=16$ processors, numbered 0 to 15. Internal boundary fluxes (shown as red arrows) are computed by locally computing the on-processor portion of the flux, sending it to the neighboring processor via \text{\tt MPI\_Send} directives, receiving the remaining portion of the flux via \text{\tt MPI\_Recv} directives, and assembling the total interface flux. Domain boundary fluxes (shown as blue dots) are handled as normal for extrapolation or wall boundary conditions, but again require MPI send/receive directives for periodic boundary conditions.\label{fig:mpi}}
\end{center}
\end{figure}


\section{2D Numerical Examples}
\label{sec:numerical2D}
In this section, we verify the correctness and accuracy of the 2D micro-macro method described above
on two test cases. The first example demonstrates the convergence of the proposed method using a manufactured solution. The second example illustrates how the scheme can compute boundary-driven flows in multi-dimensions. For the second example, we also carry out a parallel scaling analysis.
The second example in particular shows the power of the micro-macro decomposition approach, as we can 
efficiently and accurately compute complex multi-dimensional gas flows in the rarefied regime.

\begin{table}[!th]
\begin{Large}
\begin{tabular}{|c|c|c|c|c|} \hline
{\normalsize {$N_x \times N_y \times N_{v_1} \times N_{v_2}$}} & {\normalsize Macro Error} & {\normalsize $\log_2\left(\text{Error Ratio}\right)$} & {\normalsize Micro Error} & {\normalsize $\log_2\left(\text{Error Ratio}\right)$} \\ \hline \hline
\hline{\normalsize $40 \times 40 \times 40 \times 40$}    & {\normalsize $1.090341 \times 10^{-1}$} & {\normalsize --} & {\normalsize $2.220934 \times 10^{-1}$} & {\normalsize --} \\ 
\hline{\normalsize $80 \times 80 \times 80 \times 80$}   & {\normalsize $6.874683 \times 10^{-2}$} & {\normalsize $0.665414$} & {\normalsize $1.414246 \times 10^{-1}$} & {\normalsize $0.651134$} \\ 
\hline{\normalsize $160 \times 160 \times 160 \times 160$}   & {\normalsize $3.858803 \times 10^{-2}$} & {\normalsize $0.833140$} & {\normalsize $7.802459 \times 10^{-2}$} & {\normalsize $0.858032$} \\ 
\hline{\normalsize $320 \times 320 \times 320 \times 320$}   & {\normalsize $2.043854 \times 10^{-2}$} & {\normalsize $0.916861$} & {\normalsize $4.142719 \times 10^{-2}$} & {\normalsize $0.913351$}
\\ \hline\end{tabular}
\caption{(\S\ref{subsec:mom_2D}: 2D2V method of manufactured solutions example) Error table showing first order convergence in both the macroscopic and microscopic variables with $\varepsilon = 0.08$, $\left(\vec{x},\vec{v}\right) \in \left[0,1\right]^2 \times \left[-6, 6 \right]^2$ with double periodic boundary conditions in $\vec{x}$, $t=0.25$, $\text{CFL}=0.35$, and $\tau$ is taken as \eqref{eqn:tau-for-2d-es-bgk}. \label{tab5-4}}
\end{Large}
\end{table}

\subsection{Method of Manufactured Solutions Example}
\label{subsec:mom_2D}
We consider here the method of manufactured solutions with the following invented ``solution''
on $(x,y) \in [0,1]^2$ with double periodic boundary conditions (i.e., \eqref{eqn:2d_periodic_micro_1}--\eqref{eqn:2d_periodic_micro_5} for the micro step and \eqref{eqn:2d_periodic_macro_1}--\eqref{eqn:2d_periodic_macro_2} for the macro step):
\begin{align}
    f\left( t, \vec{x}, \vec{v} \right) = 
    \left( e^{-\left(v_1 - 1\right)^2 - \left(v_2 - 1\right)^2} + 
   2e^{-\left(v_1 + 1\right)^2 - \left(v_2 + 1\right)^2} \right) \Bigl(2 - 
   \sin\bigl(2 \pi \left(t-x\right) \bigr)  \cos\bigl(2 \pi \left(t-y\right)\bigr) \Bigr),
\end{align}
from which we calculate the quantities:
\begin{gather}
\left(\rho, \, \rho u_1, \, \rho u_2, \, \pr_{11}, \, \pr_{12}, \, \pr_{22} \right) = 
\Bigl( 18, \, -6, \, -6, \, 25, \, 16, \, 25 \Bigr) \left(\frac{\pi}{3}-\frac{\pi}{6}\cos\bigl(2 \pi \left(t-y\right)\bigr) \sin\bigl(2 \pi \left(t-x\right)\bigr) \right), \\
g = \frac{f - {\mathcal M}}{\varepsilon}, \quad \text{where} \quad
\mathcal{M} = \frac{27}{25} \, e^{-\left(\frac{3v_1}{5}+\frac{1}{5}\right)^2
- \left(\frac{3v_2}{5}+\frac{1}{5}\right)^2} \Bigl(2 - 
   \sin\bigl(2 \pi \left(t-x\right) \bigr)  \cos\bigl(2 \pi \left(t-y\right)\bigr) \Bigr).
%
\end{gather}
Importantly, this distribution has a spatially and temporally varying heat flux. Since this ``solution'' does not actually solve the Boltzmann-ES-BGK equation, we need to compute the residual,
\begin{align}
S\left(t,\vec{x},\vec{v}\right) = f_{,t} + v_1 f_{,x} + v_2 f_{,y} 
 - \dfrac{\tau}{\varepsilon}\left(\mathcal{G} - f\right),
\end{align}
and include it as a source term in the micro-macro method.
In particular, the 2D micro-macro scheme described above can be run as normal with two small exceptions. First, in the micro update, we modify $\widehat{G}$ as defined in \eqref{eqn5-97} as follows: 
\begin{align}
\widehat{G}_{ijk\ell}^n \leftarrow \widehat{G}_{ijk\ell}^n 
   + \frac{1}{\tau^n_{ij}} \left( {\mathcal I} - \Pi_{\maxw} \right) \left[S\right]\left(t^n, \vec{x_{ij}}, \vec{v_{k\ell}} \right).
\end{align}
Second, in the macro update, we add one extra step in between \eqref{eqn:macro2d-split3} 
and \eqref{eqn:macro2d-split4}:
\begin{align}
\label{eqn:macro2d-split3b}
\Delta t:& \quad 
 \vec{Q^{\star\star\star}_{ij}} \leftarrow \vec{Q^{\star\star\star}_{ij}} + \Delta t \int_{\reals^2}
 \left(1, v_1, v_2, v_1^2, v_1 v_2, v_2^2 \right)^T \, S\left(t^n, \vec{x_{ij}}, \, \vec{v} \right) \, d\vec{v}.
\end{align}

A convergence study for this problem is shown in Table \ref{tab5-4}. In these calculations we take $\left(\vec{x},\vec{v}\right) \in \left[0,1\right]^2 \times \left[-6, 6 \right]^2$ with double periodic boundary conditions in $\vec{x}$, the solution is computed out to $t=0.25$ with $\text{CFL}=0.35$, and the Knudsen number is $\varepsilon = 0.08$. For the collision parameters, we take \eqref{eqn:tau-for-2d-es-bgk}.
The numerical errors shown are the relative $L_2$ errors of the macroscopic, $\vec{Q} = \left(\rho, \, \rho u_1, \, \rho u_2, \, \mathbb{E}_{11},   \, \mathbb{E}_{12},  \, \mathbb{E}_{22} \right)$, and microscopic, $g$, variables: 
\begin{align}
\begin{matrix}
\text{macro} \\
\text{error}
\end{matrix} :=& \sqrt{ {{\displaystyle \sum_{i=1}^{N_x}\sum_{j=1}^{N_y}} \left\| \vec{Q_{ij}} - \vec{Q^{\text{exact}}_{ij}} \right\|^2}
\Bigg/ {{\displaystyle \sum_{i=1}^{N_x}\sum_{j=1}^{N_y}} \left\| \vec{Q^{\text{exact}}_{ij}} \right\|^2}}, \\
\begin{matrix}
\text{micro} \\
\text{error}
\end{matrix} :=& \sqrt{ {{\displaystyle \sum_{i=1}^{N_x}} \sum_{j=1}^{N_y} \sum_{k=1}^{N_{v_1}} 
\sum_{\ell=1}^{N_{v_2}}} \left| {G_{ijk\ell}} - {g^{\text{exact}}_{ijk\ell}} \right|^2
\Bigg/ {\displaystyle \sum_{i=1}^{N_x}  \sum_{j=1}^{N_y} \sum_{k=1}^{N_{v_1}} \sum_{\ell=1}^{N_{v_2}}}\left| {g^{\text{exact}}_{ijk\ell}} \right|^2}.
\end{align}
The results in Table \ref{tab5-4} clearly show first-order convergence in all variables. This table can be generated using the Python/C++/MPI companion code \cite{code:RossmanithSar2025} by going to the code directory and typing:
\begin{tcolorbox}
\begin{verbatim}
cd MPI-2d/example1-manufactured-solution/; sbatch myjob
\end{verbatim}
\end{tcolorbox}


\subsection{Cylindrical Sod Shock Tube Problem}
\label{subsec:cylindrical-sod}
Next, we apply the 2D micro-macro scheme to a cylindrically symmetric version of the 1D Sod shock tube problem \cite{article:Sod1978}, which was considered by Boscheri and Dimarco \cite{article:Boscheri2020} (a similar problem with slightly different initial conditions was also considered by Boscheri and Dumbser \cite{article:Boscheri2013}). The initial conditions are
\begin{align}
    \left( \rho, \rho u_1, \rho u_2, {\mathbb E}_{11},  {\mathbb E}_{12}, {\mathbb E}_{22} \right)(t=0,x,y) &= 
    \begin{cases}
        \, \vec{Q_{\text{in}}} \, = \bigl( 1.000, \, \, 0, \, \, 0, \, \, 5.0, \, \, 0.0, \, \, 5.0\bigr) & r < 0.5, \\
        \vec{Q_{\text{out}}} = \bigl( 0.125, \, \, 0, \, \, 0, \, \, 0.5, \, \, 0.0, \, \, 0.5 \bigr) & r > 0.5, 
    \end{cases} \\
    g(t=0,x,y,v_1,v_2) &= 0,
\end{align}
where $r=\sqrt{x^2+y^2}$.

These piecewise constant initial conditions are implemented by looping over each grid cell and computing the 
minimum and maximum distance from the four cell corners to the origin:
\begin{equation}
  r_{\text{min}} = \sqrt{\min\left\{ \left(x_i\pm \Delta x/2\right)^2 + \left(y_j\pm \Delta y/2\right)^2 \right\}}
  \quad \text{and} \quad
  r_{\text{max}} = \sqrt{\max\left\{ \left(x_i\pm \Delta x/2\right)^2 + \left(y_j\pm \Delta y/2\right)^2 \right\}},
\end{equation}
and setting:
\begin{equation}
  \vec{Q^0_{ij}} = \begin{cases}
        \vec{Q_{\text{in}}} & \text{if} \quad r_{\text{max}} < 0.5, \\
        \vec{Q_{\text{out}}}& \text{if} \quad r_{\text{min}} > 0.5, \\
        w \, \vec{Q_{\text{in}}} + \left( 1 - w \right) \, \vec{Q_{\text{out}}} & \text{otherwise},
    \end{cases}
\end{equation}
where $w$ is the percentage of the cell inside the disk of radius 0.5. In practice, we estimate $w$ by laying down $N^2$ uniformly spaced points inside the cell, computing the radius of each point, and then counting the number of points that are a distance of less than 0.5 away from the origin:
\begin{equation}
  w \approx \frac{\sum\limits_{a=1}^{N} \sum\limits_{b=1}^N \left( \sqrt{\left(x_i - \frac{\Delta x}{2} + 
   \left(a-\frac{1}{2}\right) \frac{\Delta x}{N} \right)^2 + \left(y_j - \frac{\Delta y}{2} + 
   \left(b-\frac{1}{2}\right) \frac{\Delta y}{N} \right)^2} < 0.5\right)}{N^2}.
\end{equation}
In the calculations below we take $N=20$.

\begin{figure}[!th]
\begin{tabular}{cc}
    (a)\includegraphics[width=0.45\textwidth]{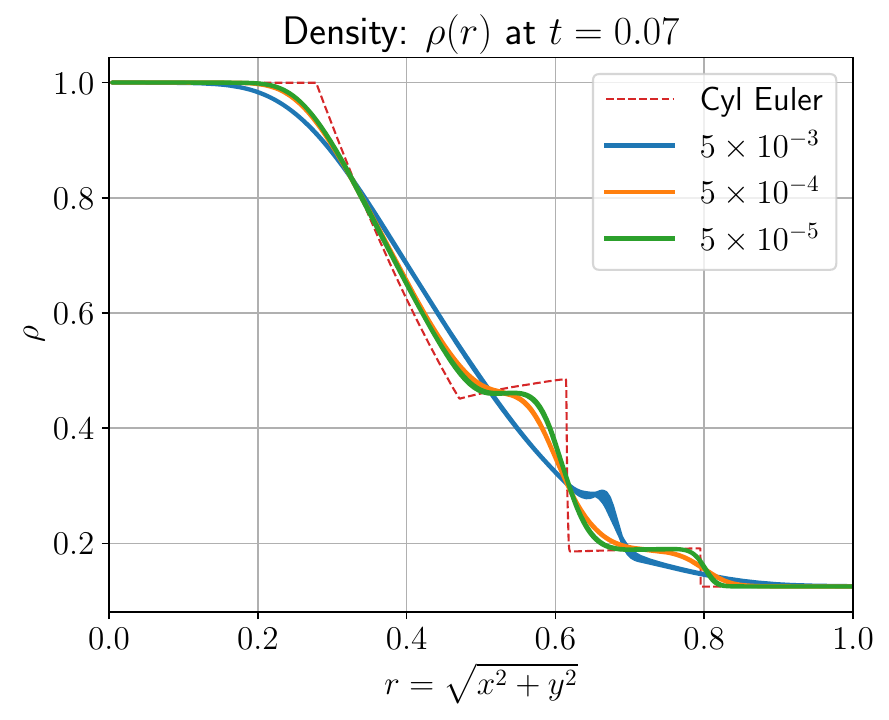} &
    (b)\includegraphics[width=0.46\textwidth]{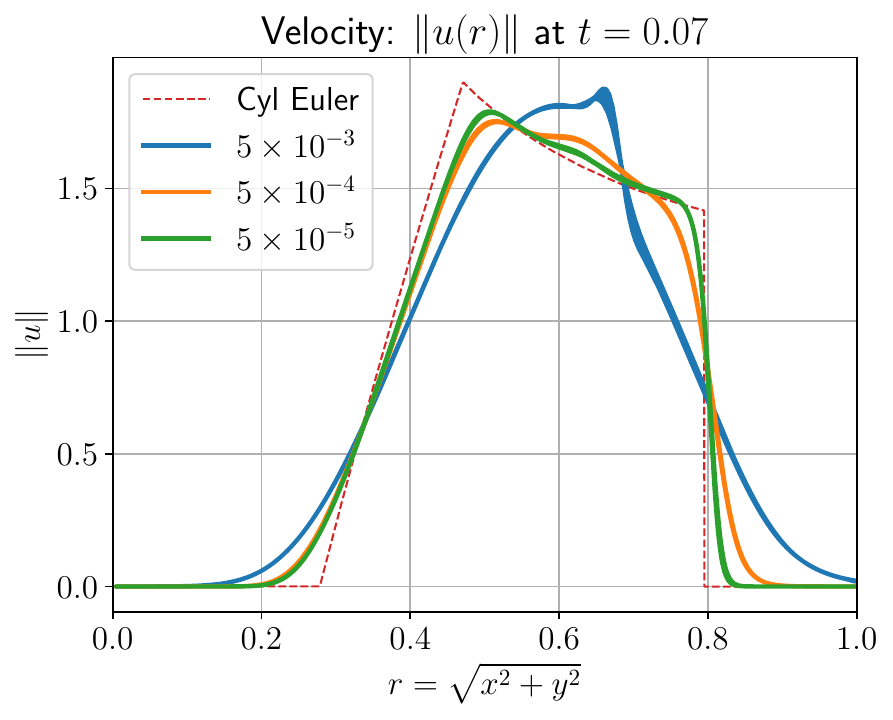} \\
    (c)\includegraphics[width=0.46\textwidth]{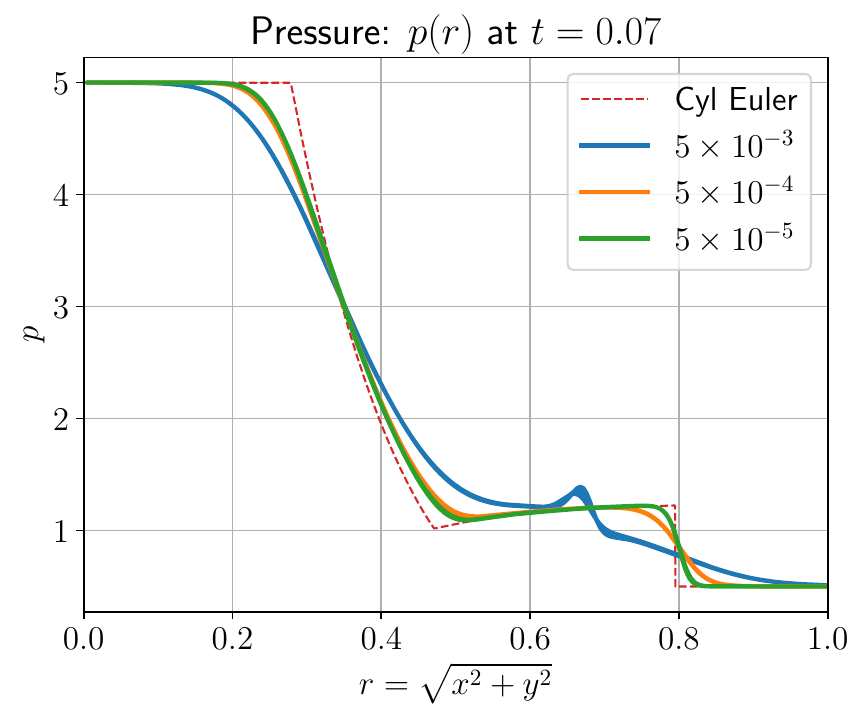} &
    (d)\includegraphics[width=0.46\textwidth]{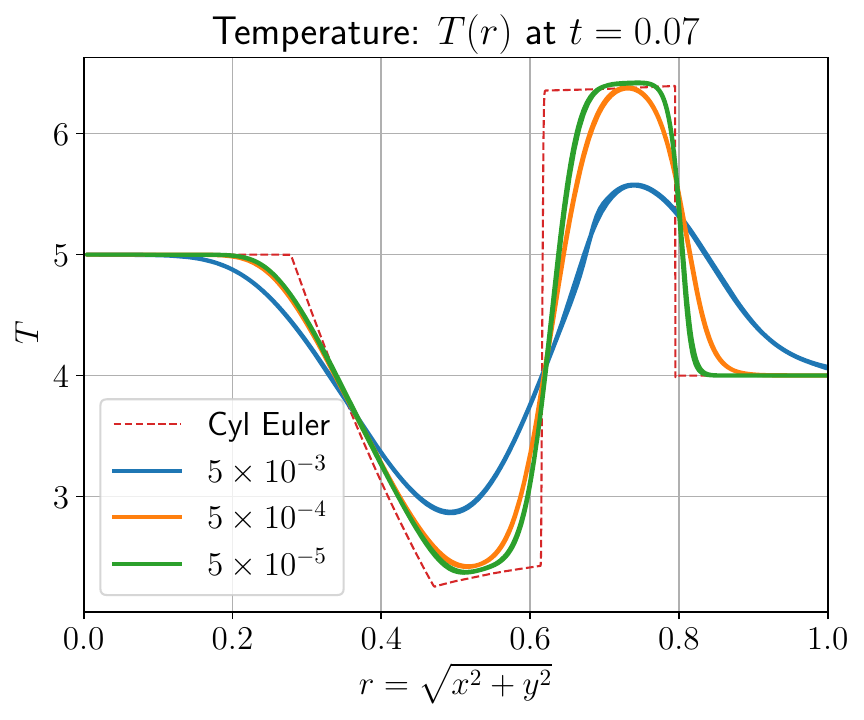}
\end{tabular}
    \caption{(\S\ref{subsec:cylindrical-sod}: 2D cylindrical Sod shock) 2D cylindrical Sod shock problem using the 2D micro-macro solver with three different Knudsen numbers: $\varepsilon = 5\times 10^{-3}, 5\times 10^{-4}, 5\times 10^{-5}$, as well as a high-resolution, 1D cylindrically-symmetric Euler solution ($\varepsilon \rightarrow 0^+$) (see \eqref{eqn:cyl-euler-1d}). 
  The full domain is $\left(\vec{x},\vec{v}\right) \in \left[-1,1\right]^2 \times \left[-9,9 \right]^2$ with constant extrapolation boundary conditions in $x$ and $y$, the grid resolution is $320^2 \times 14^2$,
   the CFL number is $\text{CFL}\approx0.95$, and $\nu$ and $\tau$ are given by \eqref{eqn:nu-tau-Boscheri-Dimarco}.
    Shown in the panels are scatter plots of the following quantities for $0 < r < 1$: (a) density: $\rho$, (b) macroscopic velocity magnitude: $\|\vec{u}\|$, (c) scalar pressure: $p$, and 
    (d) scalar temperature: $T$. 
    \label{fig:cylindrical-euler}}
\end{figure}

In this example, we take $\left(x,y,v_1,v_2\right) \in \left[-1,1\right]^2 \times \left[-9, 9 \right]^2$ with constant extrapolation boundary conditions in $x$ and $y$ (see \eqref{eqn:2d_extrap_micro_1}--\eqref{eqn:2d_extrap_micro_5} and \eqref{eqn:2d_extrap_macro_1}--\eqref{eqn:2d_extrap_macro_4}), the grid resolution is $N_x \times N_y \times N_{v_1} \times N_{v_2} =320^2 \times 14^2$, the solution is computed out to $t=0.07$ with 213 time steps ($\text{CFL} \approx 0.94648$ and
$\Delta t \approx 3.2864\times10^{-4}$), and the collision parameters are chosen as
\begin{equation}
\label{eqn:nu-tau-Boscheri-Dimarco}
\nu = 0 \quad \text{and} \quad \tau = \rho
\end{equation}
to match Boscheri and Dimarco \cite{article:Boscheri2020}.
Scatter plots at the final time, $t=0.07$, for the density, $\rho$, velocity magnitude, $\| \vec{u} \|$, scalar pressure, $p$, and scalar temperature, $T$, in the range $0 < r < 1$ are shown in Figure \ref{fig:cylindrical-euler} with three different Knudsen numbers: $\varepsilon = 5\times 10^{-3}, 5\times 10^{-4}, 5\times 10^{-5}$. For each variable, we also include a 1D cylindrically symmetric compressible Euler solution (i.e., the $\varepsilon \rightarrow 0^+$ limit):
\begin{align}
\label{eqn:cyl-euler-1d}
	\frac{\partial}{\partial t} 
	\begin{bmatrix} \rho \\ \rho u_r \\ {\mathcal E} \end{bmatrix}
	+
	\frac{\partial}{\partial r}
	\begin{bmatrix} \rho u_r \\ \rho u_r^2 + p \\ u_r \left( {\mathcal E} + p \right) \end{bmatrix} = 
	-\frac{1}{r}
	\begin{bmatrix}
	\rho u_r \\ \rho u_r^2 \\ u_r \left( {\mathcal E} + p \right)
	\end{bmatrix}, \quad
	{\mathcal E} = p + \frac{1}{2} u^2_r,
\end{align}
computed with a fourth-order discontinuous Galerkin solver with 5000 elements using the {\sc DoGPack} software \cite{dogpack}, where we use a reflecting boundary at $r=0$ and extrapolation boundary conditions at $r=1$. The results in Figure \ref{fig:cylindrical-euler} are consistent with those in the literature (e.g., \cite{article:Boscheri2020}). This figure can be generated using the Python/C++ companion code \cite{code:RossmanithSar2025} by going to the code directory and typing:
\begin{tcolorbox}
\begin{verbatim}
cd 2d/example2-circular-shock/; python make_paper_figures.py
\end{verbatim}
\end{tcolorbox}

\subsection{Lid Driven Cavity Problem}
\label{subsec:liddriven}
Next, we apply it to a standard test case in fluid dynamics: the lid-driven cavity problem. In this setup, a gas is confined in a box with four impermeable walls (modeled here as diffusely reflecting), and the top wall moves to the right at a constant prescribed speed. The moving top wall induces a circular flow inside the box. Following the work of Rana, Torrilhon, and Struchtrup \cite{Rana2013ACavity},
we consider this problem with an intermediate Knudsen number, $\varepsilon=0.08$, for which dynamics beyond Navier-Stokes-Fourier occur. We demonstrate that the micro-macro approximation with a coarse mesh in velocity can produce the correct solution. The main challenge, just as in the R13 model considered by Rana et al. \cite{Rana2013ACavity}, is the implementation of the boundary conditions.

The domain we consider is a unit square: $(x,y) \in [0,1] \times [0,1]$. 
At the kinetic level, we assume that the four edges of the unit square are diffusely reflecting walls,
meaning that the walls are in thermodynamic equilibrium with the following Maxwell-Boltzmann distribution functions:
\begin{alignat}{2}
\label{eqn:lid_driven_bc_1}
f\left(t,x=0,y,v_1,v_2\right) &= \maxw^{\text{lft}}\left(t,y,v_1,v_2\right) = \frac{\rho_{\text{lft}}(t,y)}{{2\pi T_{\text{wall}}}} \exp\left[-\frac{v_1^2+v_2^2}{2 T_{\text{wall}}}\right] \quad &v_1>0&, \\
\label{eqn:lid_driven_bc_2}
f\left(t,x=1,y,v_1,v_2\right) &= \maxw^{\text{rgt}}\left(t,y,v_1,v_2\right) = \frac{\rho_{\text{rgt}}(t,y)}{{2\pi T_{\text{wall}}}} \exp\left[-\frac{v_1^2+v_2^2}{2 T_{\text{wall}}}\right] \quad &v_1<0&, \\
\label{eqn:lid_driven_bc_3}
f\left(t,x,y=0,v_1,v_2\right) &= \maxw^{\text{bot}}\left(t,x,v_1,v_2\right) = \frac{\rho_{\text{bot}}(t,x)}{{2\pi T_{\text{wall}}}} \exp\left[-\frac{v_1^2+v_2^2}{2 T_{\text{wall}}}\right] \quad &v_2>0&, \\
\label{eqn:lid_driven_bc_4}
f\left(t,x,y=1,v_1,v_2\right) &= \maxw^{\text{top}}\left(t,x,v_1,v_2\right) = \frac{\rho_{\text{top}}(t,x)}{{2\pi T_{\text{wall}}}} \exp\left[-\frac{\left(v_1 - U_{\text{lid}}\right)^2+v_2^2}{2 T_{\text{wall}}}\right] \quad &v_2<0&,
\end{alignat}
where $T_{\text{wall}}$ is the constant wall temperature, $U_{\text{lid}}>0$ is the constant rightward velocity of the top wall, and the densities, $\rho_{\text{lft}}(t,y)$, $\rho_{\text{rgt}}(t,y)$, $\rho_{\text{bot}}(t,x)$, and
$\rho_{\text{top}}(t,x)$ are chosen so that the normal velocity at the wall is zero. In practice, this means that we set these densities so that
\begin{equation}
\iint_{\vec{v} \cdot \vec{n} \, > \, 0} \left(\vec{v} \cdot \vec{n}\right)  {\mathcal M}\left(t,\vec{x},\vec{v}\right) \, d\vec{v} + 
\iint_{\vec{v} \cdot \vec{n} \, < \, 0} \left(\vec{v} \cdot \vec{n} \right)  {\mathcal G}\left(t,\vec{x},\vec{v}\right) \, d\vec{v} = 0,
\end{equation}
where ${\mathcal G}$ is the reconstructed Gaussian distribution based on the local fluid variables 
on the interior side of the wall, and $\vec{n}$ is the inward pointing unit normal at the wall.

We need several modifications to implement these boundary conditions in the micro-macro solver.
First, for the micro update, we modify \eqref{eqn:zmicrox} as follows:
\begin{align}
Z^n_{\left(1,j,k,\ell\right)} &= \left[v_1\right]_{k}^{-} \left( 
\frac{G^n_{\left(2, j, k, \ell\right)} - G^n_{\left(1, j, k, \ell\right)}}{\Delta x} \right)+
 \left[v_1\right]_{k}^+ \left( \frac{G^n_{\left(1, j, k, \ell\right)} - \inred{$0$}}{\Delta x} \right), \\
Z^n_{\left(N_x, j, k, \ell\right)} &= \left[v_1\right]_{k}^{-} \left( \frac{\inred{$0$} - G^n_{\left(N_x,  j, k, \ell\right)}}{\Delta x} \right)+ \left[v_1\right]_{k}^{+} \left( \frac{G^n_{\left(N_x, j, k, \ell\right)} - G^n_{\left(N_x-1, j, k, \ell\right)}}{\Delta x} \right),
\end{align}
and modify \eqref{eqn:zmicroy} as follows:
\begin{align}
Z^{\star}_{{\left(i, 1, k,\ell\right)}} &= \left[v_2\right]_{\ell}^{-} \left( \frac{G^{\star}_{\left(i, 2, k, \ell\right)} - G^{\star}_{\left(i, 1, k, \ell\right)}}{\Delta y} \right)+ \left[v_2\right]_{\ell}^{+} \left( \frac{G^{\star}_{\left(i, 1, k, \ell\right)} - \inred{$0$}}{\Delta y} \right), \\
Z^{\star}_{\left(i, N_y, k, \ell\right)} &= \left[v_2\right]_{\ell}^{-} \left( \frac{ \inred{$0$} - G^{\star}_{\left(i, N_y, k, \ell\right)}}{\Delta y} \right)+ \left[v_2\right]_{\ell}^{+} \left( \frac{G^{\star}_{\left(i, N_y, k, \ell\right)} - G^{\star}_{\left(i, N_y -1, k, \ell\right)}}{\Delta y} \right).
\end{align}

The most complex modification we make is that we
no longer use the analytic projection given by equation \eqref{eqn5-97} for the finite volume cells immediately adjacent to the walls. Instead, for the cells immediately adjacent to the wall, we directly compute the following upwind approximation of the term $\vec{v} \cdot \nabla_{\vec{x}} {\mathcal M}$:
\begin{equation}
\begin{split}
M_{ijk\ell} &= \left[v_1\right]_{k}^{-} \left( 
\frac{\maxw^n_{\left(i+1, j, k, \ell\right)} - \maxw^n_{\left(i, j, k, \ell\right)}}{\Delta x} \right) + \left[v_1\right]_{k}^{+} \left( \frac{\maxw^n_{\left(i, j, k, \ell\right)}-\maxw^n_{\left(i-1, j,k,\ell\right)}}{\Delta x} \right) \\
&+ \left[v_2\right]_{\ell}^{-} \left( \frac{\maxw^{n}_{\left(i, j+1, k, \ell\right)} - \maxw^{n}_{\left(i, j, k, \ell\right)}}{\Delta y} \right)  
+ \left[v_2\right]_{\ell}^{+} \left( \frac{\maxw^{n}_{\left(i, j, k, \ell\right)} - \maxw^n_{{\left(i, j-1, k,\ell\right)}} }{\Delta y} \right),
\end{split}
\end{equation}
where the ``ghost cell'' Maxwell-Boltzmann distributions are computed from the wall boundary conditions:
\begin{equation}
\left\{
\maxw^n_{\left(0, j,k,\ell\right)}, \,
\maxw^n_{\left(N_x+1, j,k,\ell\right)}, \,
\maxw^n_{\left(i, 0,k,\ell\right)}, \,
\maxw^n_{\left(i, N_y+1,k,\ell\right)} 
\right\} \leftarrow 
\left\{
\maxw^{\text{lft}}_{\left(j,k,\ell\right)}, \,
\maxw^{\text{rgt}}_{\left(j,k,\ell\right)}, \,
\maxw^{\text{bot}}_{\left(i,k,\ell\right)}, \,
\maxw^{\text{top}}_{\left(i,k,\ell\right)}
\right\}.
\end{equation}
In these expressions, the densities can be computed analytically as
\begin{gather}
\label{eqn:twodim_wall_densities_1}
R_{a}^{\pm} := \sqrt{\frac{{\mathbb T}_{aa}}{T_{\text{wall}}}} \exp\left[-\frac{u_a^2}{2 {\mathbb T}_{aa}} \right] + u_a \sqrt{\frac{\pi}{2T_{\text{wall}}}} \left( \text{erf}\left[\frac{u_a}{\sqrt{2 {\mathbb T}_{aa}}}\right] \pm 1 \right), \\
\label{eqn:twodim_wall_densities_2}
\left[\rho^{\text{lft}}\right]^n_{j} = \left[ \rho R^{-}_1 \right]^n_{\left(1, \, j\right)}, 
\quad
\left[\rho^{\text{rgt}}\right]^n_{j} = \left[ \rho R^{+}_1 \right]^n_{\left(N_x, \, j\right)}, 
\quad
\left[\rho^{\text{bot}}\right]^n_{i} = \left[ \rho R^{-}_2 \right]^n_{\left(i, \, 1\right)},
\quad
\left[\rho^{\text{top}}\right]^n_{i} = \left[ \rho R^{+}_2 \right]^n_{\left(i, \, N_y\right)}.
%
%
%
\end{gather}
%
%
Finally, at least for the micro portion of the numerical method, we need to compute the orthogonal projection onto the complement of the collision operator nullspace, which results in the following modification 
of \eqref{eqn5-97}:
\begin{equation}
\widehat{G}_{ijk\ell}^{n} = 
\begin{cases}
\text{right-hand side of \eqref{eqn5-97}} & \quad \text{if} \quad 1<i<N_x, \quad 1<j<N_y, \\
-\frac{1}{\tau^n_{ij}} \left(\mathcal{I} - \Pi_{\maxw^n}\right)\left[ M_{ijk\ell}\right] & \quad
\text{if} \quad i=1 \quad \text{or} \quad i=N_x \quad \text{or} \quad j=1 \quad \text{or} \quad j = N_y.
\end{cases}
\end{equation}

In the macro portion of the update, we need to know the states in each ``ghost cell'' to compute the wall fluxes. We take the following values:
\begin{align}
\Bigl( \,
\rho, \, \, u_1, \, \, u_2, \, \, {\mathbb T}_{11}, \, \, {\mathbb T}_{12}, \, \, {\mathbb T}_{22}  \,
\Bigr)^{\star}_{\left(0,j\right)} &\leftarrow  
\Bigl( \, \left[ \rho^{\text{lft}} \right]^{\star}_{j}, \, \, 0, \, \, 0, \, \, T_{\text{wall}}, \, \, 0,  \, \, T_{\text{wall}} \, \Bigr), \\
\Bigl( \,
\rho, \, \, u_1, \, \, u_2, \, \, {\mathbb T}_{11}, \, \, {\mathbb T}_{12}, \, \, {\mathbb T}_{22} 
\, \Bigr)^{\star}_{\left(N_x,j\right)} &\leftarrow  
\Bigl( \, \left[ \rho^{\text{rgt}} \right]^{\star}_{j}, \, \, 0, \, \, 0, \, \, T_{\text{wall}}, \, \, 0,  \, \, T_{\text{wall}} \, \Bigr), \\
\Bigl( \, 
\rho, \, \,  u_1, \, \,  u_2, \, \,  {\mathbb T}_{11}, \, \, {\mathbb T}_{12}, \, \, {\mathbb T}_{22} 
\, \Bigr)^{\star\star}_{\left(i,0\right)} &\leftarrow 
\Bigl( \, \left[ \rho^{\text{bot}} \right]^{\star\star}_{i}, \, \, 0, \, \, 0, \, \, T_{\text{wall}}, \, \, 0,  \, \, T_{\text{wall}} \, \Bigr), \\
\Bigl( \, 
\rho, \, \,  u_1, \, \,  u_2, \, \, {\mathbb T}_{11}, \, \, {\mathbb T}_{12}, \, \, {\mathbb T}_{22} 
\, \Bigr)^{\star\star}_{\left(i,N_y\right)} &\leftarrow 
\Bigl( \, \left[ \rho^{\text{top}} \right]^{\star\star}_{i}, \, \, U_{\text{lid}}, \, \, 0, \, \, T_{\text{wall}}, \, \,  0, \,  \, T_{\text{wall}} \, \Bigr),
\end{align}
where the wall densities are again defined via \eqref{eqn:twodim_wall_densities_1}--\eqref{eqn:twodim_wall_densities_2}, but this time with different time indeces.

Finally, to complete the boundary condition treatment in the macro update, 
we also need to specify the heat flux at the boundaries; for that, we take:
\begin{equation}
{\mathbb H}^{n+1}_{\left(1/2, \, j\right)} = \frac{1}{2} {\mathbb H}^{n+1}_{\left(1, \, j\right)} \, , \quad
{\mathbb H}^{n+1}_{\left(N_x + 1/2, \, j\right)} = \frac{1}{2} {\mathbb H}^{n+1}_{\left(N_x, \, j\right)} \, , \quad
{\mathbb H}^{n+1}_{\left(i, \, 1/2\right)} = \frac{1}{2} {\mathbb H}^{n+1}_{\left(i, \, 1\right)} \, , \quad
{\mathbb H}^{n+1}_{\left(i, \, N_y + 1/2\right)} = \frac{1}{2} {\mathbb H}^{n+1}_{\left(i, \, N_y\right)} \, .
\end{equation}

\begin{figure}[!th]
\begin{tabular}{cc}
    (a)\includegraphics[width=0.45\textwidth]{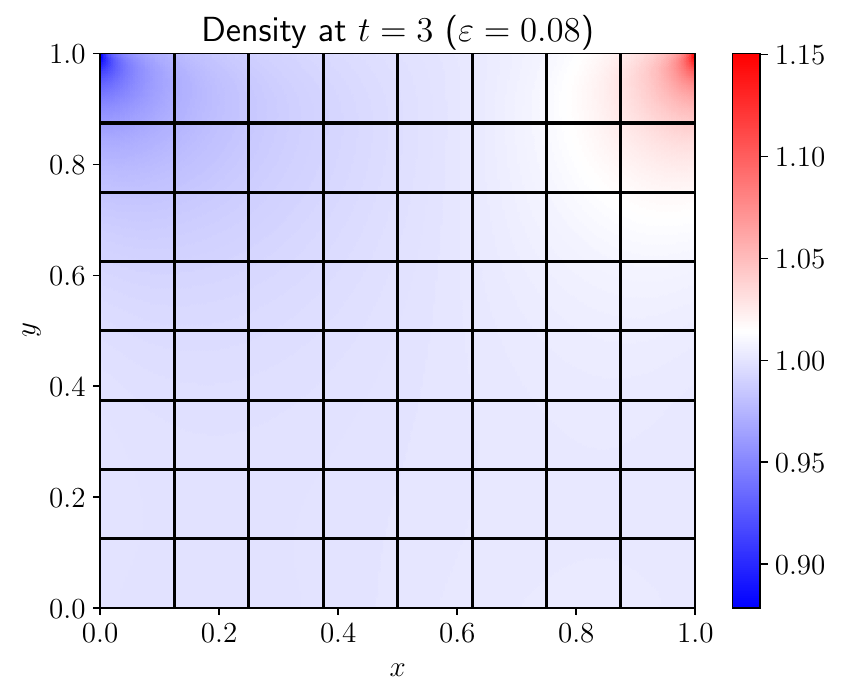} &
    (b)\includegraphics[width=0.46\textwidth]{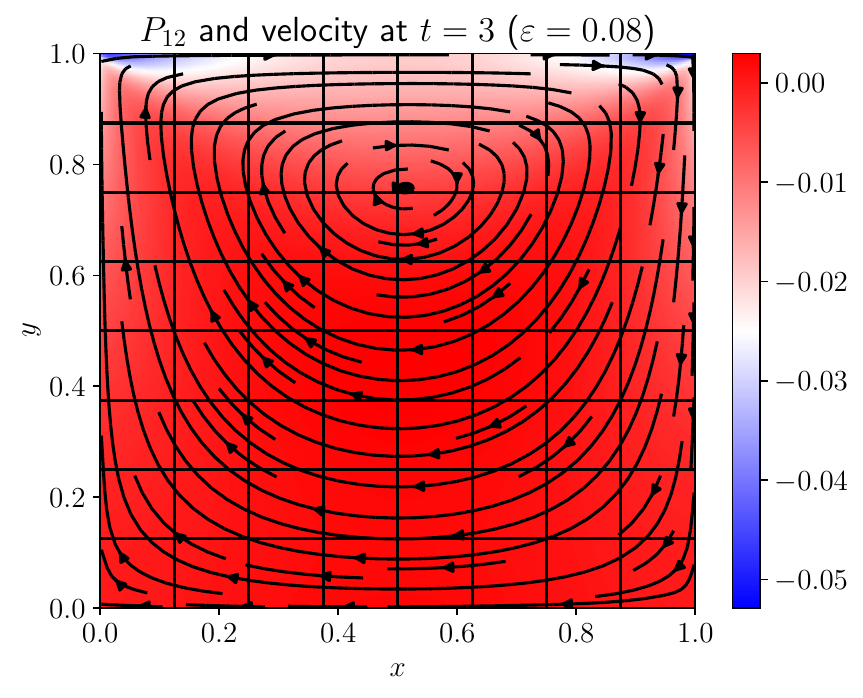} \\
    (c)\includegraphics[width=0.46\textwidth]{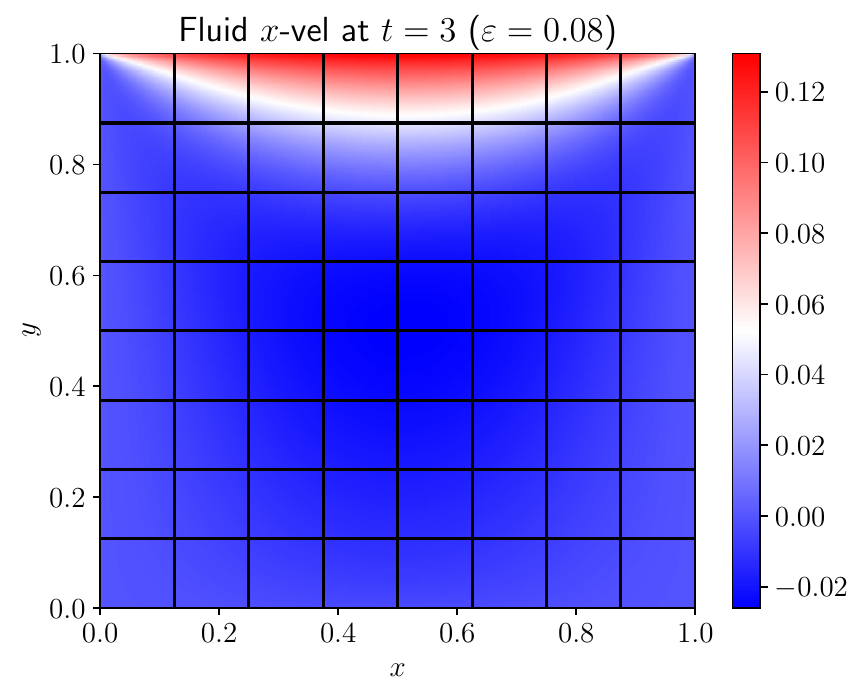} &
    (d)\includegraphics[width=0.46\textwidth]{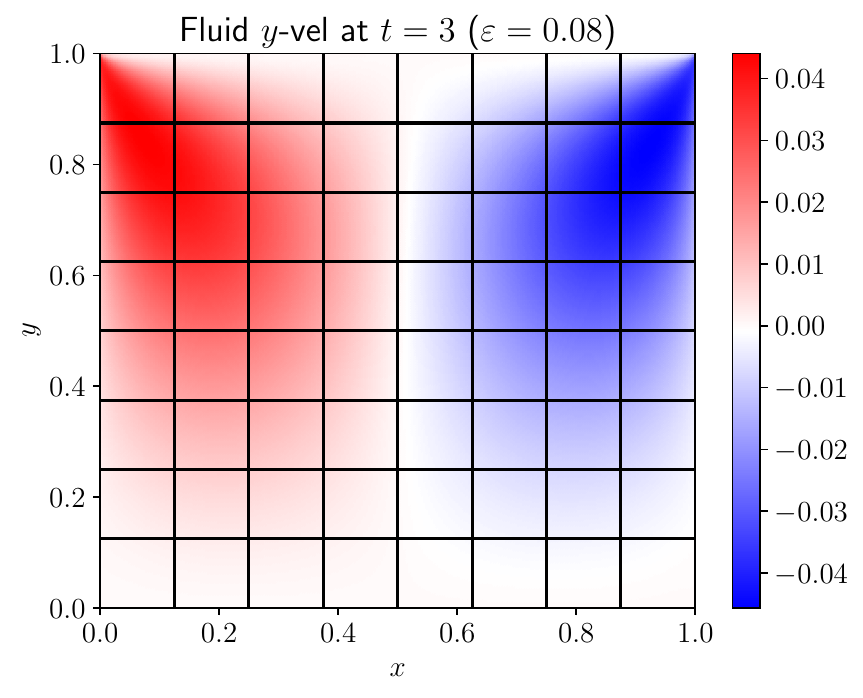} \\
    (e)\includegraphics[width=0.46\textwidth]{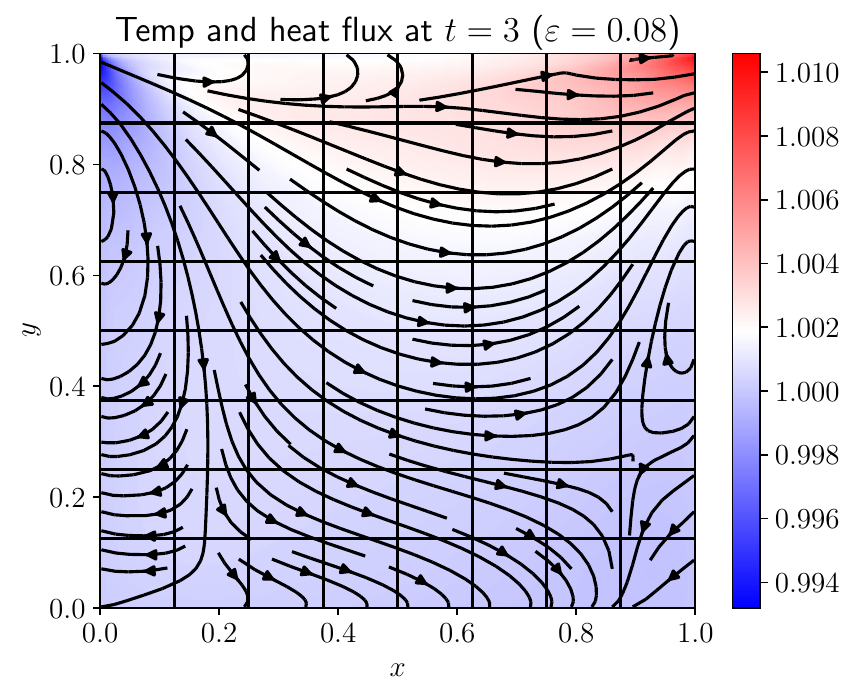} &
    (f)\includegraphics[width=0.46\textwidth]{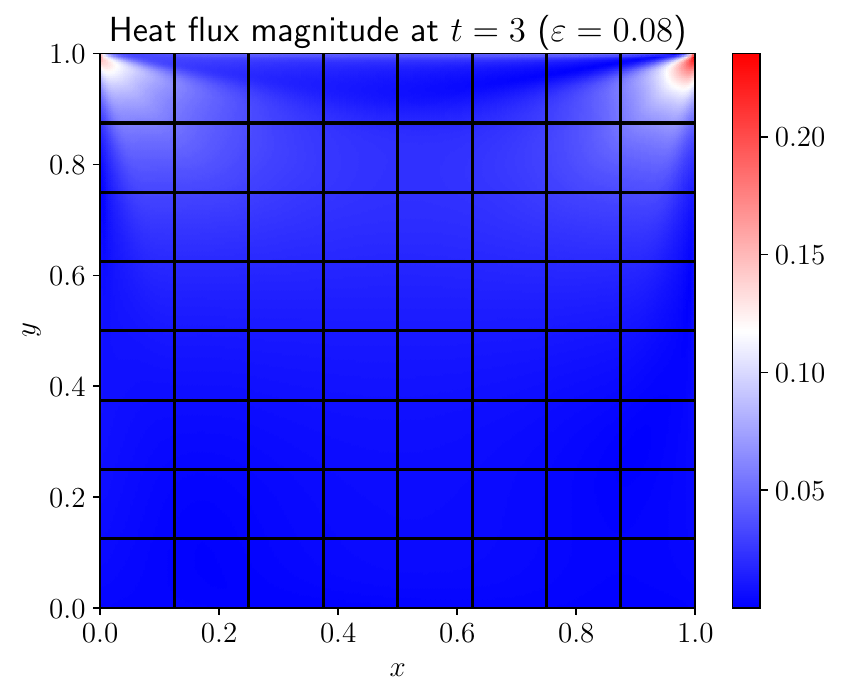}
\end{tabular}
    \caption{(\S\ref{subsec:liddriven}: lid-driven cavity) Simulation results for the lid-driven cavity problem with Knudsen number $\varepsilon=0.08$ on a mesh with $N_x=N_y=240$ and $N_{v_x}=N_{v_y}=14$ using 64 processors ($30\times30\times14\times14$ mesh on each processor). Shown in these panels are color plots at time $t=3$ (at this time the solution has settled down to a steady-state) of the
    (a) fluid density $\left(\rho\right)$, (b) off-diagonal pressure $\left({\mathbb P}_{12}\right)$ in color along with velocity streamlines $\left(\vec{u}\right)$ as contours with arrows, (c) $x$-component of the fluid velocity $\left(u_{1}\right)$, (d) $y$-component of the fluid velocity $\left(u_{2}\right)$, (e) scalar temperature $\left(T\right)$ in color along with heat flux streamlines $\left(\vec{h}\right)$ as contours with arrows, and (f) heat flux magnitude (which is a proxy for where microscopic effects are important).
    \label{fig5-3}}
\end{figure}

The results of a simulation for the lid-driven cavity problem with Knudsen number $\varepsilon=0.08$ on the domain
$\left(x,y,v_1,v_2\right) \in [0,1] \times [0,1] \times [-5,5] \times [-5,5]$ with initial condition:
\begin{equation}
\left( \, \rho, \, \, u_1, \, \, u_2, \, \, {\mathbb T}_{11}, \, \, {\mathbb T}_{12}, \, \, {\mathbb T}_{22} \, \right)\left(t=0,x,y\right) = \left( \, 1, \, \, 0, \, \, 0, \, \, 1, \, \, 0, \, \, 1 \, \right) \quad \text{and} \quad
g\left(t=0,x,y,v_1,v_2\right) = 0,
\end{equation}
is shown in Figure \ref{fig5-3}. 
We take the 2D ES-BGK collision parameter defined in expression \eqref{eqn:tau-for-2d-es-bgk}.
The mesh is $N_x=N_y=240$ and $N_{v_x}=N_{v_y}=14$, and the simulation is run using 64 processors, given a resolution of $30\times30\times14\times14$ on each processor. In this simulation, we take $T_{\text{wall}} = 1$, $U_{\text{lid}} = 0.16$, a CFL number of $0.95$, and run to a final time of $t=3$, at which point the solution has settled down to a steady-state. 
Shown in the panels of Figure \ref{fig5-3} are color plots at time $t=3$ of the (a) fluid density $\left(\rho\right)$, (b) off-diagonal pressure $\left({\mathbb P}_{12}\right)$ in color along with velocity streamlines $\left(\vec{u}\right)$ as contours with arrows, (c) $x$-component of the fluid velocity $\left(u_{1}\right)$, (d) $y$-component of the fluid velocity $\left(u_{2}\right)$, (e) scalar temperature $\left(T\right)$ in color along with heat flux streamlines $\left(\vec{h}\right)$ as contours with arrows, and (f) heat flux magnitude (which is a proxy for where microscopic effects are important). The heat flux vector, $\vec{h}$, mentioned in Panel \ref{fig5-3}(e) is defined as follows:
\begin{equation}
 \vec{h} := \frac{1}{2} \left( \, {\mathbb H}_{111} + {\mathbb H}_{122}, \, \, 
                      {\mathbb H}_{112} + {\mathbb H}_{222} \, \right).
\end{equation}

Figure \ref{fig5-3}(b) clearly shows the clockwise circular flow inside the unit square, which is typical for lid-driven cavity problems.
Figure \ref{fig5-3}(e) shows the tell-tale signs of the microscopic influence on the bulk flow; what we see in this panel so-called anti-Fourier heat flow, which occurs because heat is flowing from cold to hot (see temperature colors and heat flux streamlines), which violates Fourier's law of heat conduction:
\[
	\text{Fourier's law of heat conduction:} \quad \vec{h} \propto - \nabla_{\vec{x}} T.
\]
The violation of this law happens because at this Knudsen number, $\varepsilon=0.08$, microscopic effects are not negligible. See Rana et al. \cite{Rana2013ACavity} for a more detailed description. Our results are consistent with those in the literature, including Rana et al. \cite{Rana2013ACavity} and 
van der Woude et al. \cite{article:vanderWoude2024}. Remarkably, we can achieve these results with a relatively low resolution in velocity: $N_{v_1} \times N_{v_2} = 14 \times 14$.
Finally, we note that these figures can be generated using the Python/C++/MPI companion code \cite{code:RossmanithSar2025} by going to the code directory and typing:
\begin{tcolorbox}
\begin{verbatim}
cd MPI-2d/example3-lid-driven-cavity/; sbatch myjob
\end{verbatim}
\end{tcolorbox} 

\subsection{Scaling Analysis}
We end the numerical examples section with a scaling study of the proposed micro-macro finite volume method on the Nova Cluster at Iowa State University \cite{computer:nova-cluster}.

\begin{figure}[!th]
\begin{center}
    \begin{tabular}{cc}
    (a) \includegraphics[width=0.44\textwidth]{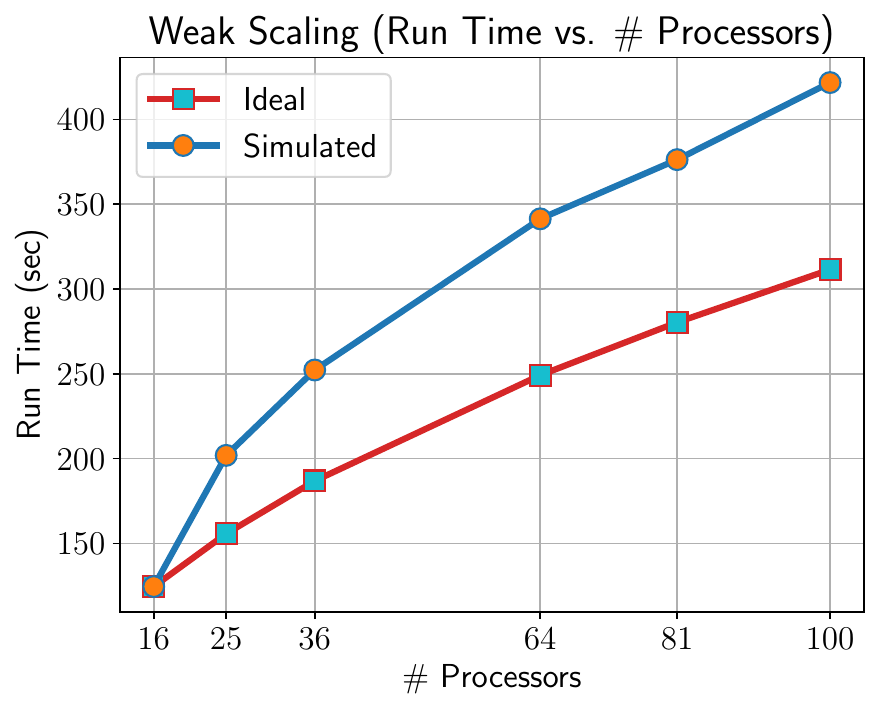} &
    (b) \includegraphics[width=0.44\textwidth]{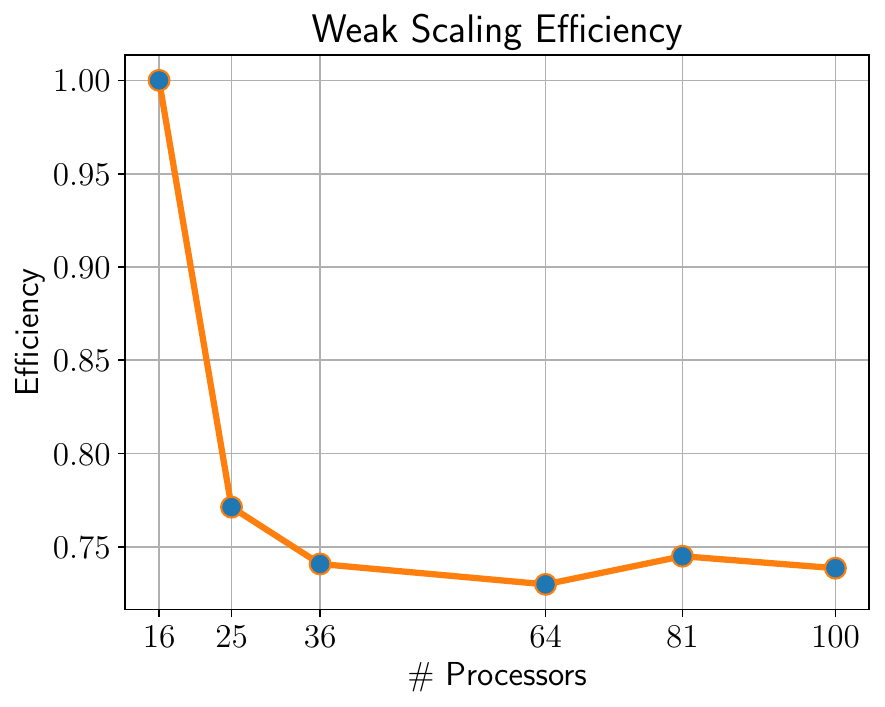}
    \end{tabular}
    \caption{(\S\ref{subsec:weak-scaling}: weak scaling analysis) Weak scaling analysis for the lid driven cavity problem. The simulation parameters are the same as those described in \S\ref{subsec:liddriven}, except that the total spatial mesh resolution varies depending on the number of processors: $N_x = N_y = 36\sqrt{\left(\text{\# Processors}\right)}$, where $\left(\text{\# Processors}\right) = 16, 25, 36, 64, 81, 100$. The velocity resolution is constant: $N_{v_1}=N_{v_2}=14$. Panel (a) shows the run time (measured in seconds) as a function of the number of processors; in addition to the experimentally computed run times, we show the ideal run times based on what the run time should be if the 16-processor run time is weakly scaled to more and more processors. In Panel (b), we show scaling efficiency, which is just the ratio of the ideal run time divided by the experimentally computed run time. \label{fig5-5}}
\end{center}
\end{figure}

\subsubsection{Weak Scaling}
\label{subsec:weak-scaling}
In weak scaling, the idea is to keep the total amount of work per processor constant as we increase the number of processors. The way we carry this out for the micro-macro method is that we keep the same simulation initial conditions, boundary conditions, and simulation parameters as those described in \S\ref{subsec:liddriven} (the lid-driven cavity problem), except that the total spatial mesh resolution varies depending on the number of processors: 
\begin{equation}
N_x = N_y = 36\sqrt{\left(\text{\# Processors}\right)}, \quad \text{where} \quad 
\left(\text{\# Processors}\right) = 16, \, 25, \, 36, \, 64, \, 81, \, 100.
\end{equation}
The velocity resolution is kept constant at $N_{{v_1}}=N_{{v_2}}=14$; and therefore,
the space-velocity mesh resolution on each processor is constant: $36 \times 36 \times 14 \times 14$. We again simulate out to time $t=3$. Note that the total number of time-steps will increase with increasing processors because the total mesh resolution increases with processor number; therefore, using 16 processors as the baseline, we define the following ideal run time:
\begin{equation}
\label{eqn:ideal-run-time-weak-scaling}
   \left(\text{ideal run time} \right) = \sqrt{\frac{\left(\text{\# Processors}\right)}{16}} 
   \times \left(\text{actual run time with 16 Processors} \right).
\end{equation}
The results of this experiment are shown in Figure \ref{fig5-5}. In Figure \ref{fig5-5}(a), we display the experimentally determined wall-clock run times (measured in seconds) as a function of the number of processors; additionally, we show the ideal run times as defined by \eqref{eqn:ideal-run-time-weak-scaling}. In Figure \ref{fig5-5}(b), we show scaling efficiency, which we define simply as the ratio of the ideal run time divided by the experimentally computed run time; this plot shows that as the number of processors increases, the efficiency decreases and saturates to a value of around $74\%$.

\begin{figure}[!th]
\begin{center}
    \begin{tabular}{cc}
    (a) \includegraphics[width=0.44\textwidth]{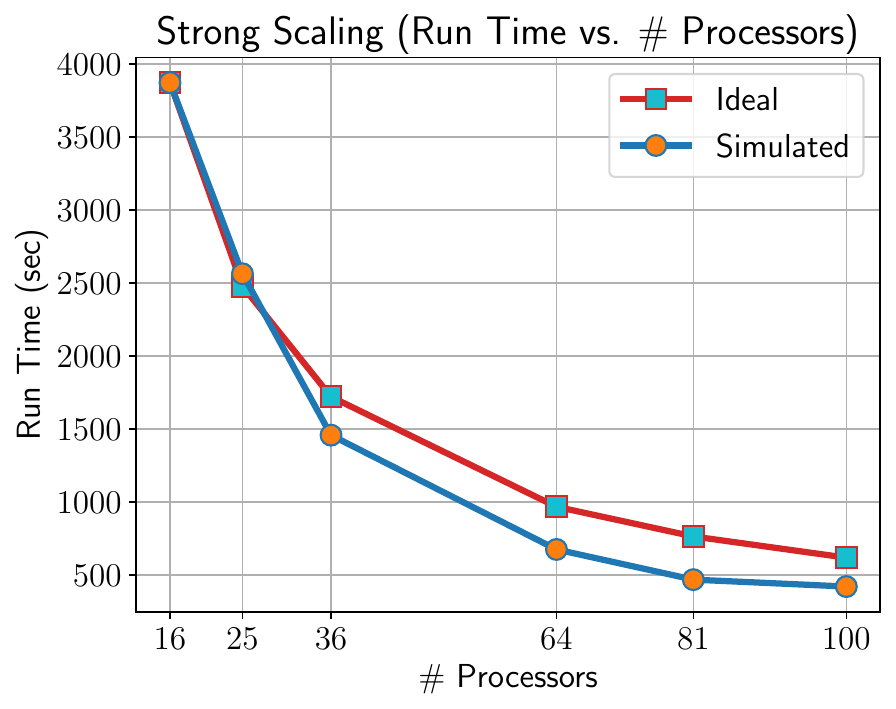} &
    (b) \includegraphics[width=0.44\textwidth]{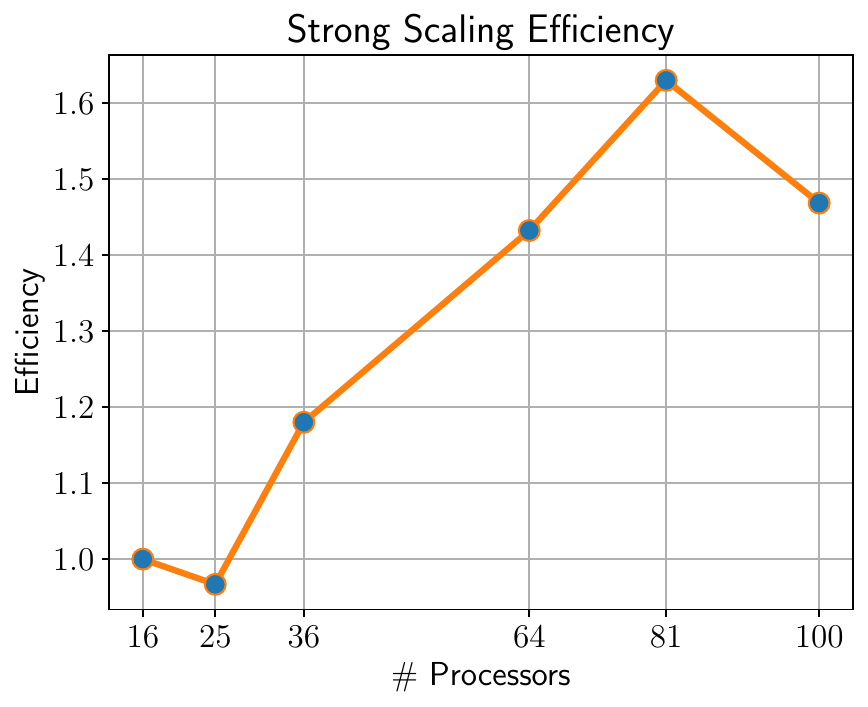}
    \end{tabular}
    \caption{(\S\ref{subsec:strong-scaling}: strong scaling analysis) Strong scaling analysis for the lid driven cavity problem. The simulation parameters are the same as those described in \S\ref{subsec:liddriven}, including
    the space and velocity resolutions: $N_x = N_y = 360$ and $N_{{v_1}}=N_{{v_2}}=14$.
     Panel (a) shows the run time (measured in seconds) as a function of the number of processors; in addition to the experimentally computed run times, we show the ideal run times based on what the run time should be if the 16-processor run time is strongly scaled to more and more processors. In Panel (b), we show scaling efficiency, which is just the ratio of the ideal run time divided by the experimentally computed run time. \label{fig5-6}}
\end{center}
\end{figure}

\subsubsection{Strong Scaling}
\label{subsec:strong-scaling}
In strong scaling, the idea is to keep the total problem size constant as we increase the number of processors. The way we carry this out for the micro-macro method is that we keep the same simulation initial conditions, boundary conditions, and simulation parameters as those described in \S\ref{subsec:liddriven} (the lid-driven cavity problem), and fix the total problem size to 
\begin{equation}
N_x \times N_y \times N_{{v_1}} \times N_{{v_2}} = 360 \times 360 \times 14 \times 14.
\end{equation}
The resolution on each processor is then given by
\begin{equation}
 \frac{360}{\sqrt{\left(\text{\# Processors}\right)}} \times \frac{360}{\sqrt{\left(\text{\# Processors}\right)}}
 \times 14 \times 14, \quad \text{where} \quad 
\left(\text{\# Processors}\right) = 16, \, 25, \, 36, \, 64, \, 81, \, 100.
\end{equation}
We again simulate out to time $t=3$. Using 16 processors as the baseline, we define the following ideal run time:
\begin{equation}
\label{eqn:ideal-run-time-strong-scaling}
   \left(\text{ideal run time} \right) = \frac{16 \times \left(\text{actual run time with 16 Processors} \right)}{\left(\text{\# Processors}\right)}.
\end{equation}
The results of this experiment are shown in Figure \ref{fig5-6}. In Figure \ref{fig5-6}(a), we display the experimentally determined wall-clock run times (measured in seconds) as a function of the number of processors; additionally, we show the ideal run times as defined by \eqref{eqn:ideal-run-time-strong-scaling}. In Figure \ref{fig5-6}(b), we show scaling efficiency, which we define simply as the ratio of the ideal run time divided by the experimentally computed run time. This plot shows that the experimentally efficiency is above one (i.e., the actual run time is better than what is predicted by the ideal value), because as the resolution per processor decreases, there are additional efficiency gains due to the amount of memory needed on each processor.

\section{Conclusion}
\label{sec:conclusions}
In this work, we developed a finite volume method for solving the kinetic Boltzmann equation using the
micro-macro decomposition approach first developed for the one-dimensional Boltzmann-BGK equation by Bennoune, Lemou, and Mieussens \cite{Bennoune2008UniformlyAsymptotics}. In particular, we extend the micro-macro approach to two spatial and two velocity dimensions, the ES-BGK collision operator, and examples with non-trivial boundary conditions. 
The scheme developed in this work was based on coupling a microscopic kinetic equation with a macroscopic fluid equation. A projection operator was used to separate the macroscopic and microscopic quantities. The updated microscopic quantity was coupled to the macroscopic update equation via the heat flux variable, which is obtained by taking the moment of the microscopic distribution function. The updated macroscopic variables were then used to calculate the projection operator, which affects the microscopic update.
As it appears in both the micro and macro equations, the collision operator was handled via L-stable implicit time discretizations. At the same time, the remaining transport terms were computed via kinetic flux vector splitting (for the macro equations) and upwind differencing (for the micro equation). Operator splitting was used both in the update's micro and macro portions to achieve optimal linear stability. The resulting scheme was applied to various test cases in 1D and 2D. In both 1D and 2D, we performed a convergence analysis using the method of manufactured solutions. We also tested the scheme on diffusely reflecting wall boundary conditions, showing how such boundary conditions can be incorporated into the micro-macro framework. The 2D version of the code was parallelized via MPI, and we presented weak and strong scaling studies with varying numbers of processors. 
Ongoing and future work will focus on extending the scheme to higher orders and problems in complex geometries.
The complete companion code to this paper, including all the presented examples, can be found in reference \cite{code:RossmanithSar2025}.




\section*{Acknowledgments}
We would like to thank the anonymous reviewers for their thoughtful comments and suggestions that helped to improve this paper.

\section*{Statements and Declarations}
\begin{description}
\item[{\bf Funding.}] This research was partially funded by US National Science Foundation Grants DMS--2012699 and DMS--2410538.

\medskip

\item[{\bf Competing interests.}] The authors have no
conflicts of interest to disclose.

\medskip

\item[{\bf Data availability statement.}] Data sharing does not apply to this article as no datasets were generated or analyzed during the current study.
\end{description}

\bibliographystyle{cas-model2-names}

\FloatBarrier
\clearpage
\bibliography{references}

\end{document}